\UseRawInputEncoding 
\documentclass[]{article}

\usepackage{amsmath,amssymb}

\DeclareFontFamily{U}{musix}{}%
\DeclareFontShape{U}{musix}{m}{n}{%
  <-12>   musix11
  <12-15> musix13
  <15-18> musix16
  <18-23> musix20
  <23->   musix29
}{}%
\newcommand*\musix{\usefont{U}{musix}{m}{n}\selectfont}
\DeclareTextFontCommand{\textmusix}{\musix}

\usepackage[dvipdfmx]{graphicx}

\usepackage[top=30truemm,bottom=30truemm,left=25truemm,right=25truemm]{geometry}

\usepackage{amsthm}
\usepackage{bm}
\usepackage{color}
\usepackage{here}
\usepackage{comment}

\theoremstyle{plain}

\theoremstyle{definition}

\theoremstyle{remark}

\begin{document}

\title{General Theory of Music by Icosahedron 1: 
\\
\Large
A bridge between ``artificial" scales and ``natural" scales, Duality between chromatic scale and Pythagorean chain, and Golden Major Minor Self-Duality}

\author{Yusuke Imai$^1$ \and  Sid C. Dellby \and Nobuaki Tanaka$^2$}
\date{%
\small
    $^1$Department of Physics, Nagoya University, Furo-cho, Chikusa-ku, Nagoya, Aichi 466-0804, Japan\\%
      $^2$Hochschule f\"ur Musik W\"urzburg, Hofstallstra$\beta$e 6-8, 97070 W\"urzburg, Germany\\[2ex]%
    \today
}

\maketitle

\small

\begin{center}
CONTACT: 93imaiyusuke@gmail.com
\end{center}

\begin{abstract}
Relations among various musical concepts are investigated through a new concept, musical icosahedron that is the regular icosahedron each of whose vertices has one of 12 tones, $C$, $C\sharp$, $D$, $E\flat$, $E$, $F$, $F\sharp$, $G$, $G\sharp$, $A$, $B\flat$, $B$.

First, we found that there exist four musical icosahedra that characterize the topology of the chromatic scale and one of the whole tone scales, and have the hexagon-icosahedron symmetry (an operation of raising all the tones of a given scale by two semitones corresponds to a symmetry transformation of the regular icosahedron): chromatic/whole tone musical icosahedra. The major triads or the minor triads are set on the golden triangles of these musical icosahedra. Also, various dualities between musical concepts are shown by these musical icosahedra: the major triads/scales and the minor triads/scales, the major/minor triads and the fundamental triads for the hexatonic major/minor scales, the major/minor scales and the Gregorian modes. Namely, these musical icosahedra connect ``natural" scales (Gregorian modes, major/minor scales, hexatonic major/minor scales) and ``artificial" scales (chromatic scale and whole tone scales).

Second, we derived duality relations between the chromatic scale and the Pythagorean chain that is a succession of the fifth based on $C$ by using musical icosahedra. We proposed Pythagorean/whole tone musical icosahedra that characterize the topology of the Pythagorean chain and one of the whole tone scales, and have the hexagon-icosahedron symmetry. The Pythagorean chain (chromatic scale) in the chromatic (Pythagorean)/whole tone musical icosahedron is constructed by ``middle" lines of the regular icosahedron. While some golden triangles correspond to the major/minor triads in the chromatic/whole tone musical icosahedra, in the Pythagorean/whole tone musical icosahedra, some golden gnomons correspond to the minor/major triads.

Third, we found four types of musical icosahedra other than the chromatic/whole tone musical icosahedra and the Pythagorean/whole tone musical icosahedra that have the hexagon-icosahedron symmetry. All the major triads and minor triads are represented by the golden triangles or the golden gnomons on each type. Then, these musical icosahedra may be applied to harmonic analysis by combining figures characterized by the golden ratio on the regular icosahedron.

All of these musical icosahedra naturally lead to generalizations of major/minor triads and scales through the symmetry of the regular icosahedron.

\end{abstract}

\newpage
\section{Introduction}
\begin{figure}
\centering
{
\resizebox*{8cm}{!}{\includegraphics{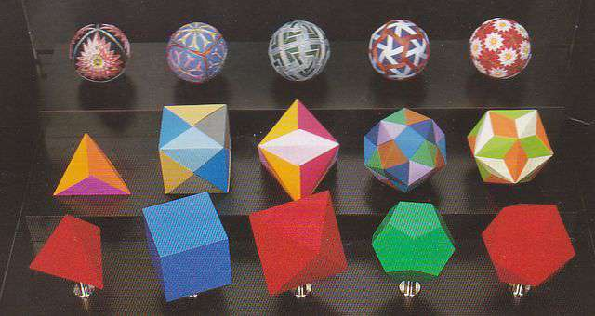}}}\hspace{5pt}
\caption{Temari balls on which polyhedra on a sphere are painted (Polyhedron Club, 2019a).} \label{sample-figure}
\end{figure}

\begin{figure}
\centering
{%
\resizebox*{10cm}{!}{\includegraphics{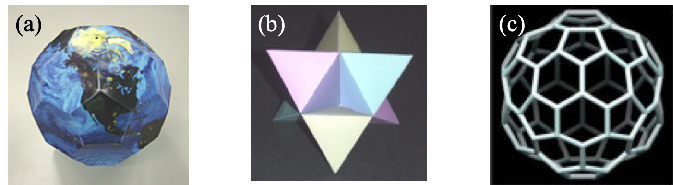}}}\hspace{5pt}
\caption{(a) Globe on a soccer-ball solid (Jaxa, 2019), (b) A folded paper of stellated octahedron (Polyhedron Club, 2019b) (c) One of the Goldberg polyhedra (Schein, 2014).} \label{Polyhedron2}
\end{figure}

The regular icosahedron\footnote{We will omit ``regular" in this paper to escape complexity.} is one of the Platonic solids that has 12 vertices. Platonic solids are a set of solids that are composed of one kind of the regular polygons and satisfy the convexity condition and the uniform condition. The name of Platonic solids originates from the story that Plato told that each of the classical four elements (fire, air, earth, water) is made of each of four Platonic solids (tetrahedron, octahedron, cube, icosahedron)\footnote{As for the dodecahedron, Plato said that ``God used it as a pattern for dividing the zodiac into its twelve signs" (Archer-Hind, 1888, 190).} in Timaeus (Archer-Hind, 1888, 190). In the reference of Cromwell (Cromwell, 1999), many historical works that include some polyhedra, such as pyramids at Giza in Egypt (square pyramid), the Sacrament of the Last Supper by Salvator Dali (the regular dodecahedron), the concept of ``bit" in the movie, Tron (octahedron and icosahedron), Shugakuin Imperial Villa in Kyoto (cuboctahedron) are shown. On the website, Polyhedron club, Temari balls, traditional Japanese toys, on which polyhedra on a sphere are painted, are shown (Fig.~\ref{sample-figure}). After Plato pointed out the importance of the Platonic solids, many ways of extension of the concept of the Platonic solids have been proposed. Archimedes proposed 13 solids that consist of polygons of many kinds and satisfy the convexity condition (Cromwell, 1999). A figure of a soccer ball is a typical example of the Archimedean solids (Fig.~\ref{Polyhedron2}(a)). Kepler tried to explain the orbits of planets in the solar system by using the Platonic solids (Kepler withdrew this idea later and found the famous Kepler's law (Kepler, 2009)), and found two Kepler-Poinsot polyhedra that are stellations (extension of edges or faces of a polyhedron) of the dodecahedron and the icosahedron in 1619, and Poinsot found two other Kepler–Poinsot polyhedra in 1809 (Coxeter \& Rouse Ball, 2010, 144-146). Cauchy proved there are no Kepler–Poinsot polyhedra other than the above four Kepler–Poinsot polyhedra (Coxeter \& Rouse Ball, 2010, 144-146). Figure~\ref{Polyhedron2}(b) shows one of the stellated polyhedra, the stellated octahedron. Recently, Schein and Gayed proposed the fourth type of polyhedra, Goldberg polyhedra that are equilateral and planar but not equiangular, and are composed of the regular triangles, rectangles, or pentagons, and hexagonal faces that are equilateral but not equiangular (see Fig.~\ref{Polyhedron2}(c), Schein \& Gayed, 2014, 2920-2925). These studies may have been driven by a simple question, ``What is a beautiful figure?".

 The number 12 is not only the number of vertices of an icosahedron, but also the number that divides an octave in European music. The development of European music is closely linked to the use of diatonic scales. From the early Middle Ages to the beginning of the Baroque period, music was conventionally built on the Gregorian modes. The development of polyphony, especially the various uses of consonance in the early Renaissance, provided an opportunity to discuss scales\footnote{We define ``a scale" as an ordered list of any tones in this paper.} in the context of chords and harmony. By the end of the 17th century, the commonly used musical scales among compositions were integrated into two types\footnote{The ionic and aeolian modes can be referred as major and minor scales in the modern sense. These Gregorian modes were added by Heinrich Glarean (1488-1563) to the conventional Gregorian modes, but they were in fact not differed from the transposed (`Transposition" is a transformation that shifts all the tones of a given scale by some semitones.) doric and Lydian church modes (Benary, 1995, 1594) in which B was conceived as B$\flat$, due to the practice of \emph{musica ficta}.},  namely the major and minor scales, which construct major and minor triads\footnote{We define ``a triad" as an ordered list of three tones in this paper.} on each first tone.

    In today's music, the major-minor tonality is still popular, but the composers in the late romantic period tried to deviate from the conventional tonal system. The music of that period certainly included the use of the Gregorian modes, but on the other hand it also aimed to deviate from the diatonic sound systems, namely through the use of the whole tone scale and the chromatic scale. It is possible to perceive some chromatic compositions in the earlier period, but their tonality was always diatonic, so the use of chromaticism has remained partial\footnote{e.g. some compositions by Carlo Gesualdo (1566-1613) and Johann Sebastian Bach (1685-1750). See also note examples 1 and 2.}.  In the compositional attempts based on the non-diatonic scales one can see the consciousness of the late romantic and modern generation that the musical settings with the major-minor sound system would soon reach the limit, because ``major and minor in the context of functional harmony since about 1850 get into a final stage"　(Benary, 1995, 1596). For example, there are some compositions after the late Romantic period that deliberately used the whole tone scale.\footnote{e.g. \emph{Voiles} by Claude Debussy (The second part of his \emph{24 Pr\'eludes}). See also note example 3.}  Doubts on diatonic settings of the generation were further proved by the attempt on using the chromatic scale, especially the twelve-tone technique of the Viennese school. From a music-historical point of view, the practice of un-diatonic settings can be understood thus as ``artificial", that apparently tried to deviate from the commonly used diatonic scales.

    Geometrical consideration on scales, however, indicates the opposite view; from the un-diatonic, ``artificial" scales one can derive the diatonic, ``natural" scales. In this process of proof, the icosahedron is assumed as a model for consideration. Since the age of antiquity, music has been widely discussed in the context of mathematics, but the icosahedron, which has twelve vertices and should therefore be appropriate for musical consideration, has been neglected in the field of music theory.

\begin{figure}[H]
\centering
{%
\resizebox*{10cm}{!}{\includegraphics{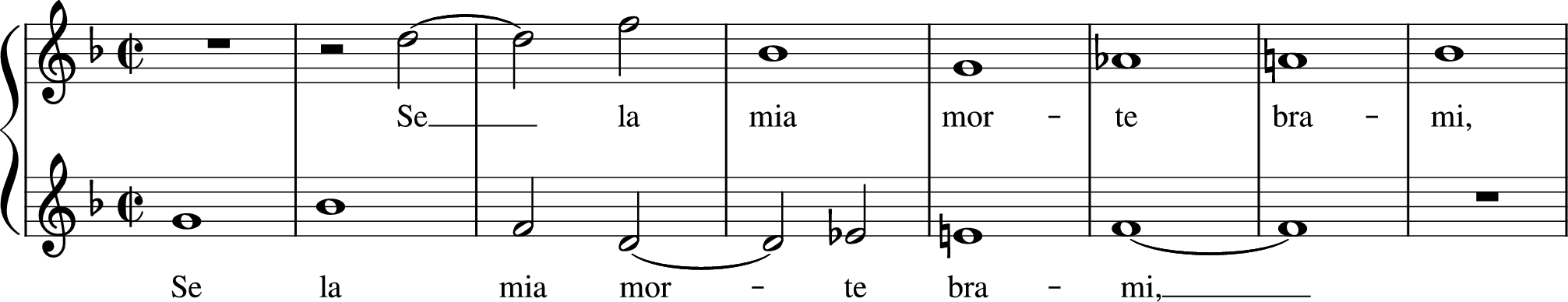}}}\hspace{5pt}
\caption{Example 1: The beginning of the madrigal \emph{Se la mia morte brami} by Carlo Gesualdo (The first piece of his sixth madrigal book. The collection was published in 1611).} \label{Exapmple_1}
\end{figure}

\begin{figure}[H]
\centering
{%
\resizebox*{10cm}{!}{\includegraphics{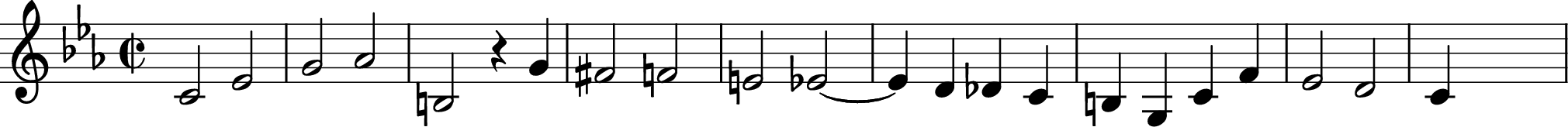}}}\hspace{5pt}
\caption{Example 2: The theme of \emph{Musical Offering} by Johann Sebastian Bach (BWV 1079, composed in 1747).} \label{Exapmple_2}
\end{figure}

\begin{figure}[H]
\centering
{%
\resizebox*{10cm}{!}{\includegraphics[clip]{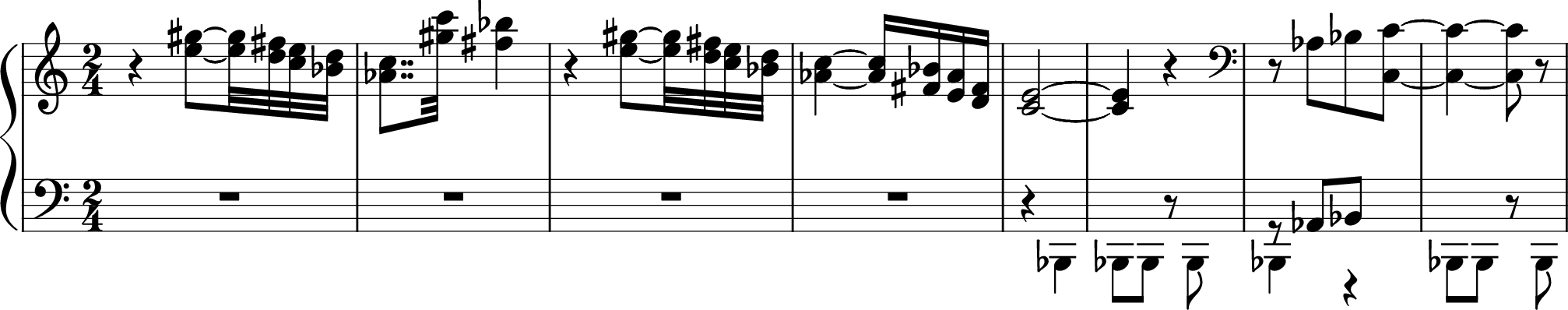}}}\hspace{5pt}
\caption{Example 3: The beginning of Claude Debussy's $Voiles$ (composed ca. 1910). Here dynamics and articulations are omitted.} \label{Exapmple_3}
\end{figure}

\begin{figure}
\centering
{%
\resizebox*{10cm}{!}{\includegraphics{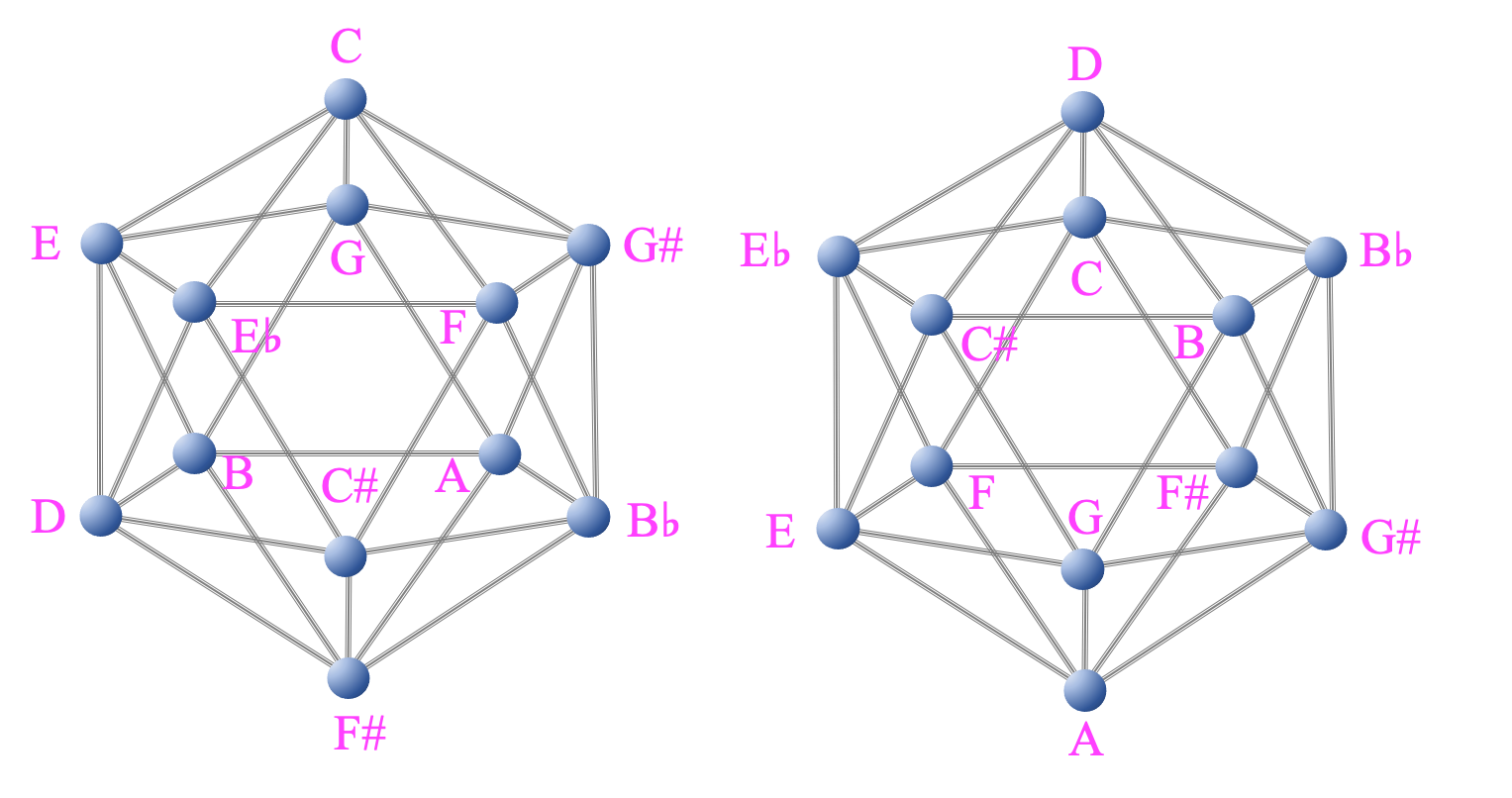}}}\hspace{5pt}
\caption{Examples of musical icosahedra.} \label{Icosahedron;sample}
\end{figure}

Our idea is simple. We attach 12 tones ($C$, $C\sharp$, $D$, $E\flat$ $E$, $F$, $F\sharp$, $G$, $G\sharp$, $A$, $B\flat$, $B$) to the vertices of the icosahedron as shown in Fig. \ref{Icosahedron;sample}. We call this ``musical icosahedron". First, we show how musical icosahedra connect ``artificial" scales and ``natural" scales through musical concepts with some beauty and without relying on any discussion of auditory frequencies. Concretely, we show the topology of the chromatic scale and the whole tone scales determines four types of musical icosahedra (chromatic/whole tone musical icosahedra) and these musical icosahedra have the hexagon-icosahedron symmetry (an operation of raising all the tones of a given scale by two semitones corresponds to a symmetry transformation of the regular icosahedron). These musical icosahedra lead to the uniqueness of the major/minor triads and the fundamental triads for the hexatonic major/minor scales through the golden triangles, and dualities between the major triads/scales and minor triads/scales, as well as between the major/minor scales and the Gregorian modes.

Then, we focus on the other musical icosahedra that have the hexagon-icosahedron symmetry.

We found there exist four musical icosahedra that reflect the topology of the Pythagorean chain (a succession of the fifth based on $C$: $C$, $G$, $D$, $A$, $E$, $B$, $F\sharp$, $C\sharp$, $G\sharp$, $E\flat$, $B\flat$, $F$) that appears in the Pythagorean intonation (Lindley, 2001), and one of the whole tone scales. We call them Pythagorean/whole tone musical icosahedra. We remark that the chromatic scale and the Pythagorean chain are only the tone sequences that are constructed by 12 tones and symmetric if a tone sequence is equivalent to the reversed tone sequence of it\footnote{The reversed tone sequence of $X_1$, $X_2$, $X_3$, $\cdots$, $X_n$ is $X_1$, $X_{11}$, $X_{10}$, $\cdots$ $X_2$.}. By the term, ``symmetric", we mean that the interval\footnote{We assume the interval is a non-negative integer less than 12. For a scale $C$, $G$, the interval between the first tone and the second tone is 7, not -5 or 5 or -7 or 19.} of neighboring two tones of a given tone sequence is constant. The interval of any neighboring two tones of the chromatic scale (Pythagorean chain) is 1 (7). The interval of any neighboring two tones of the reversed chain of the Pythagorean chain, $C$, $F$, $B\flat$, $E\flat$, $G\sharp$, $C\sharp$, $F\sharp$, $B$, $E$, $A$, $D$, $G$, is 5. Also, the interval of any neighboring two tones of the reversed scale of the chromatic scale, $C$, $B$, $B\flat$, $A$, $G\sharp$, $G$, $F\sharp$, $F$, $E$, $E\flat$, $D$, $C\sharp$, is 11. Note that 7, 5, 11 are the only natural numbers that are greater than 1 and less than 12, and coprime with 12. Note also that $C$, $E$, $G\sharp$ is constructed by 3 tones and a symmetric scale whose interval of any neighboring two tones is $4$, and $12/4=3$. Therefore, one can show the algebraic duality between the chromatic scale and the Pythagorean chain. We show the geometrical reflection of that duality through musical icosahedron, and find beautiful connections between them. These musical icosahedra also lead to the uniqueness of the major/minor triads and the fundamental triads for the hexatonic major/minor scales through golden gnomons, and dualities between the major triads/scales and minor triads/scales, as well as between the major/minor scales and the Gregorian modes.

In addition, we deal with four other musical icosahedra that have the hexagon-icosahedron symmetry. We call them exceptional musical icosahedra. Because the duality between the chromatic/whole tone musical icosahedra and the Pythagorean/whole tone musical icosahedra shown in this paper is beautiful, it may be strange that there exist other kinds of musical icosahedra that have the hexagon-icosahedron symmetry. In the chromatic/whole tone (Pythagorean/whole tone) musical icosahedra, all the major triads or minor triads are represented by the golden triangles (gnomons). Then, how are the major/minor triads represented in the exceptional musical icosahedra? We find that the self-duality with respect to the major/minor triads and the golden triangle/gnomon holds for them.

Also, symmetry transformations of the icosahedron naturally lead to generalization of the major/minor triads and scales through these musical icosahedra. This means that we obtain the third axis for the major/minor scale space; the first axis shows whether a given scale is major or minor, the second axis shows the key of a given scale, and then, the third axis shows what kind of the musical icosahedra correspond to a given scale. At this point, it is unclear that the generalized major/minor triads and scales naturally correspond to human emotion, just as the major (minor) triads/scales correspond to bright (dark) impressions.

\begin{figure}
\centering
\hspace*{3cm}
{%
\resizebox*{7cm}{!}{\includegraphics{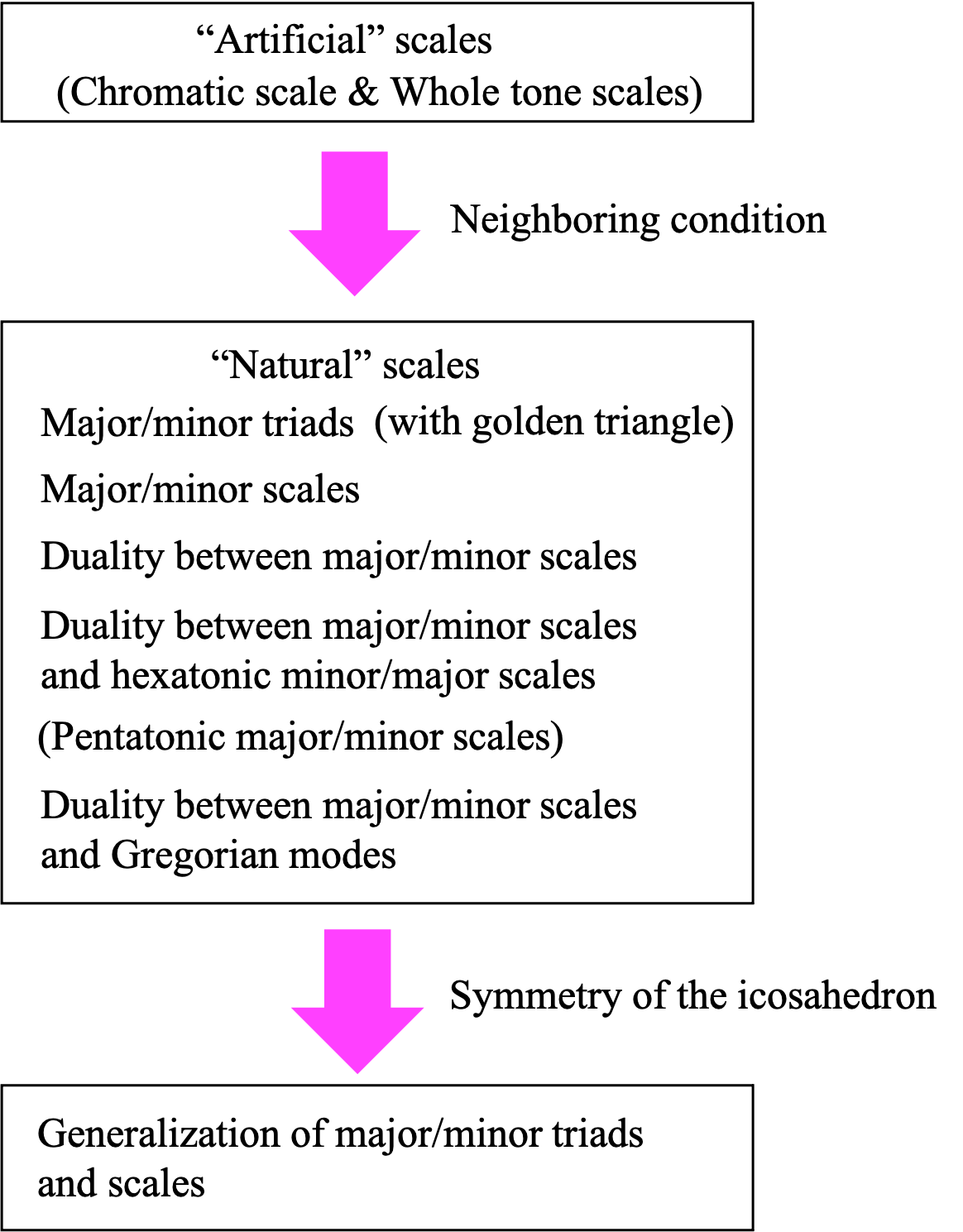}}}\hspace{5pt}
\caption{Summary of our study.} \label{summary}
\end{figure}

Figure~\ref{summary} shows the summary of our study.

This paper is organized as follows. In Section II, we review some properties of the icosahedron: symmetry and golden triangle. In Section III, we introduce the concept of a musical icosahedron and explain how the topology of scales are reflected on the musical icosahedra (neighboring/distant condition), and show links between ``artificial" scales and ``natural" scales. In Section IV, we introduce the musical icosahedra satisfying the neighboring condition for the Pythagorean chain and neighboring condition for one of the whole tone scales and distant condition for the other whole tone scale, and show some duality relations between the chromatic scale and the Pythagorean chain. In Section V, we show other kinds of musical icosahedra that have the same fundamental properties as the musical icosahedra shown in Section II and Section III, and their special characteristics (Golden Major Minor Self-Duality). In Appendix, we list all the scales we deal with in this paper. Because we would like to make this paper accessible to a variety of readers, we review musically well-known concepts in a simplified way for those who are not familiar with music and do not use mathematical-set-theory for those who are not familiar with mathematics (also for avoiding the confusion caused by the fact that the ``set" in music theory is different from the ``set" in mathematics).

This paper may be regarded as the fourth part of the following stories of strange connections among polyhedra, music, and astronomy; (i) Plato said that the four Platonic solids each make up one of the four classical elements, (ii) Plato paired astronomy (visual art) and scale (auditory art), (iii) Kepler tried to explain the orbits of planets in the solar system by using Platonic solids in the beginning, and then, he derived the famous Kepler's law stimulated by music of the spheres, (iv) We found that the Platonic solids link various musical concepts.

\section{Regular Icosahedron}
\begin{figure}
\centering
{%
\resizebox*{10cm}{!}{\includegraphics{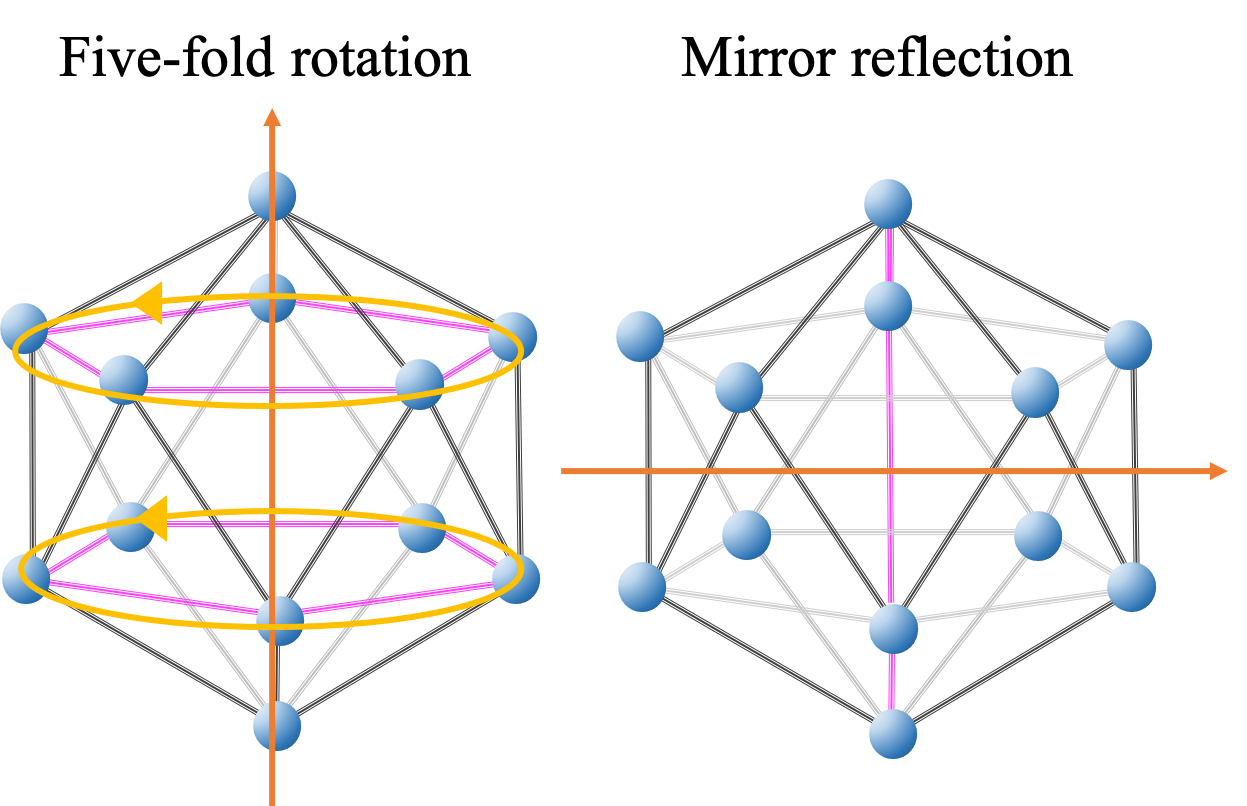}}}\hspace{5pt}
\caption{Five-fold rotation and mirror reflection.} \label{five_mirror}
\end{figure}

In this section, we review properties of the icosahedron with a focus on symmetry and golden triangle. First, we list all the 120 symmetric transformations of the icosahedron. 
\\
\\
\indent
(i) One of the five-fold rotations is given in Fig.~\ref{five_mirror}. 

\begin{figure}
\centering
{%
\resizebox*{10cm}{!}{\includegraphics{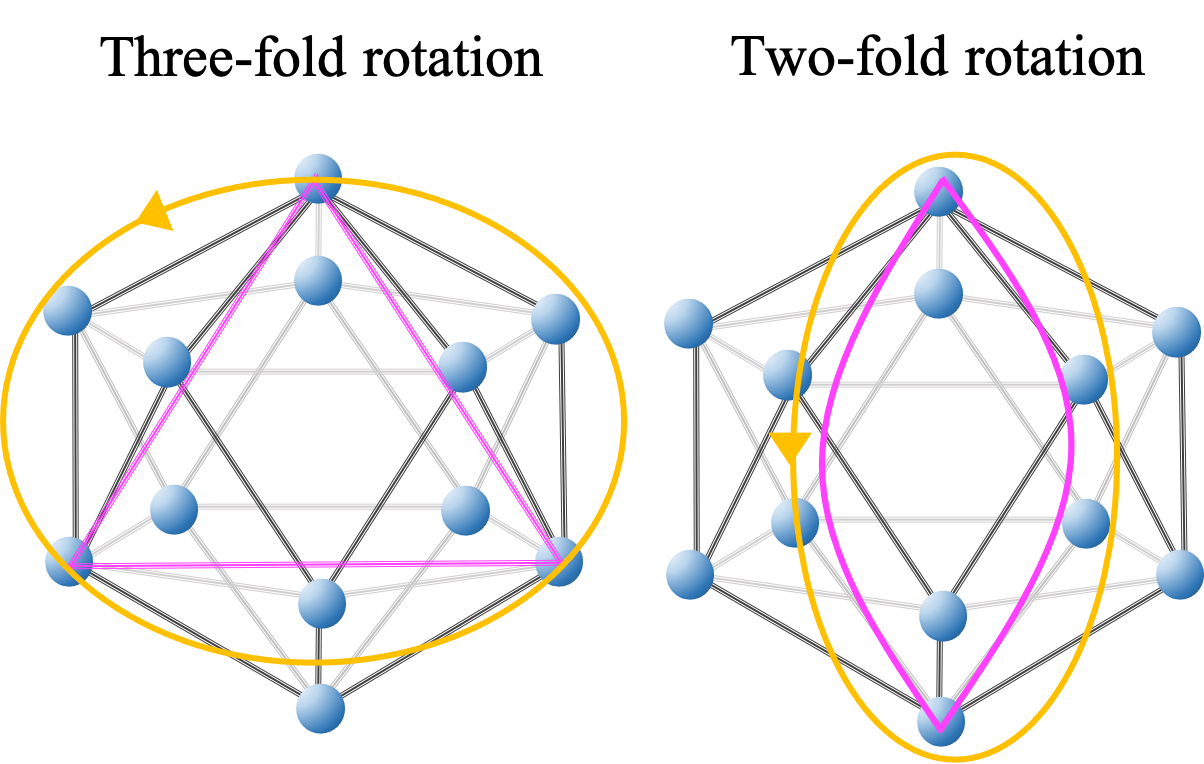}}}\hspace{5pt}
\caption{Three-fold rotation and two-fold rotation.} \label{two_three_rotation}
\end{figure}

(ii) One of the mirror reflections is given in Fig.~\ref{five_mirror}. 
\\
\\
\indent
The above two transformations fix the top vertex (and the bottom vertex). One can obtain all the transformations if considering transformations that do not fix the top vertex. In fact, the number of transformations that can be made by the above two transformations is 10, and the number of vertices is 12. Those transformations involve the five-fold rotation and (iii) a two-fold rotation and (iv) a three-fold rotation as shown in Fig. \ref{two_three_rotation}. By the three-fold rotation, one can move the two vertices to the top, and by the three-fold rotation after the two-fold rotation, the three vertices are moved to the top. By three-fold rotation after one time, three times, or four times of the five-fold rotation, the three vertices are moved to the top, and by the three-fold rotation after two-fold rotation after one time, two times, or four times of the five-fold rotation, the three vertices are moved to the top. Then, we moved all the vertices to the top.

Then, we can list all the symmetric transformations of the icosahedron. By defining $E$ as the elementary transformation, $C_5$ as the five-fold rotation, $M$ as the mirror reflection, $C_3$ as the three-fold rotation, $C_2$ as the two-fold rotation, one can represent those 120 transformations as follows.
\begin{align}
\{Y, \ M*Y\},
\end{align} 
where 
\begin{align}
Y &= \{X, \ C_3*X, \ C_3^2*X, \ C_2*X, \ C_2*C_3*X, \  C_2*C_3^2*X, 
\nonumber\\
&\qquad C_5*C_3*X, C_5^3*C_3*X, \ C_5^4*C_3*X, \ C_5*C_3*C_2*X, \ C_5^2*C_3*C_2*X, 
\nonumber\\
&\qquad C_5^4*C_3*C_2*X\},
\\
X &= \{E, C_5, C_5^2, C_5^3, C_5^4\}.
\end{align}

\begin{figure}
\centering
\hspace*{2.5cm}
{%
\resizebox*{5cm}{!}{\includegraphics{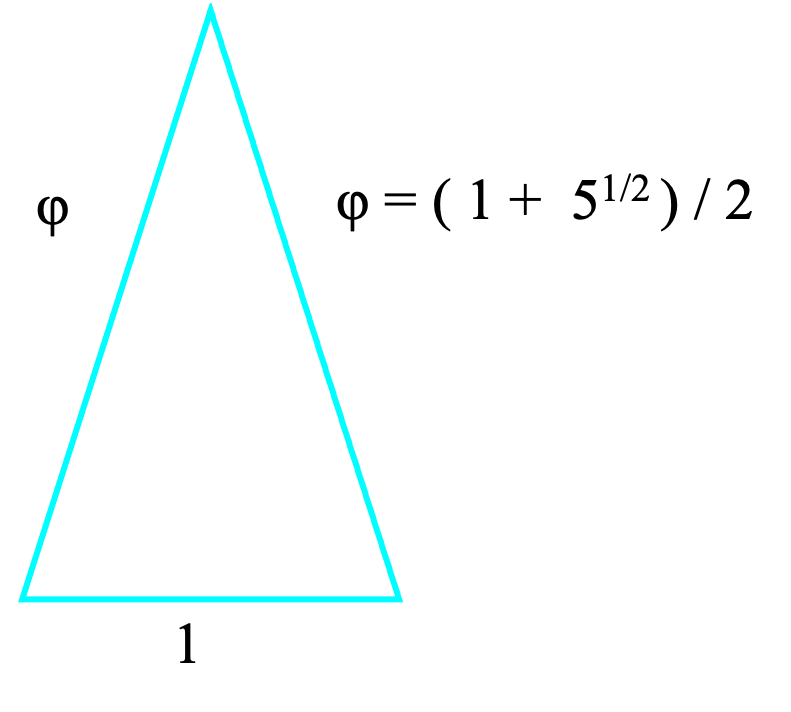}}}\hspace{5pt}
\caption{A golden triangle.} \label{golden_triangle}
\end{figure}

\begin{figure}
\centering
{%
\resizebox*{4cm}{!}{\includegraphics{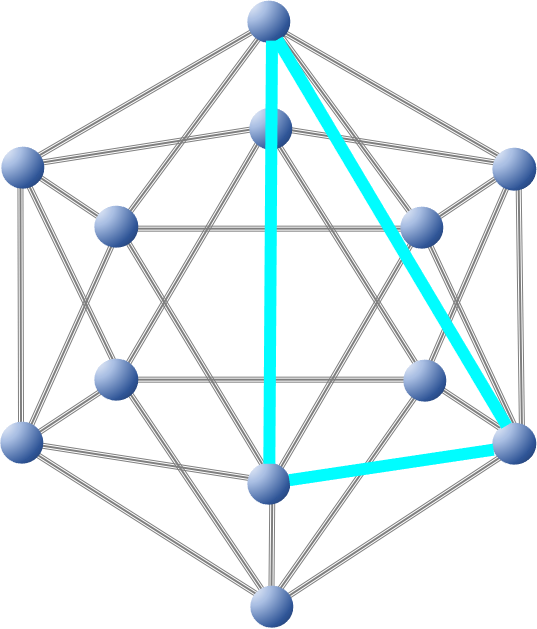}}}\hspace{5pt}
\caption{A golden triangle obtained by linking three vertices of the icosahedron.} \label{golden_icosahedron}
\end{figure}

Then, we also remark that there are 60 golden triangles in the icosahedron. The golden triangle (Fig.~\ref{golden_triangle}) is an isosceles triangle with two long edges and one short edge that satisfy the condition that the ratio of length of the long edge and the short edge is given by the golden ratio\footnote{The golden ratio appears in many things, e.g., seed heads, shells, aloe plants, ${\it Mona}$ ${\it Lisa}$ (Leonardo da Vinci), ${\it David}$ (Michelangelo), ${\it The}$ ${\it School}$ ${\it of}$ ${\it Athens}$ (Raphael), DNA, planets’ orbits, and black holes (Meisner, Araujo, 2018). This paper add a new example that has beauty of golden ratio to the list: the major/minor triads.}. Figure~\ref{golden_icosahedron} shows a golden triangle obtained by linking three vertices of the icosahedron.

\newpage
\section{Musical Icosahedron}

\subsection{Definition}
In this subsection, we introduce a new concept, ${\it a \ musical \  icosahedron}$. A musical icosahedron is the icosahedron each of whose vertices has one of the 12 tones, $C$, $C\sharp$, $D$, $E\flat$, $E$, $F$, $F\sharp$, $G$, $G\sharp$, $A$, $B\flat$, $B$ without duplication (Fig.~\ref{Icosahedron;sample}). There exist 12!=479001600 musical icosahedra. Therefore, this concept is useless unless one finds meaningful subsets of musical icosahedra. In the next section, we will find them by using topology of the chromatic scale and the whole tone scales.

In addition, identification of two musical icosahedra that can be transformed into one another by a symmetry transformation of the icosahedron is a way to impose a restriction on the set of musical icosahedra\footnote{Musical icosahedra should be transformed while retaining the connections.}. In this paper, we sometimes use this identification and we call this ``musical isomorphism".

Also, one can consider harmonies and scales on a musical icosahedron by connecting neighboring tones in a given scale on the musical icosahedron\footnote{We also use symmetry transformations of the icosahedron naturally for a figure representing a scale on the musical icosahedron.}. In this paper, we deal with the major/minor triads, the major/minor scales, the hexatonic major/minor scales, the Gregorian modes and their generalizations by this method.

Although we propose some ways to construct scales or find relations among various scales by using some musical icosahedra in this paper, they are not the only ways to use musical icosahedra, and there may be more meaningful ways to use musical icosahedra because there are over 400 million musical icosahedra and many fascinating relations of musical concepts must be hidden. Using other polyhedra may also be a fascinating approach to dealing with scales. For example, four Archimedean solids have 24 vertices and these are useful in dealing with the microtone found in the Arabic Maqam or enharmonic tetrachords.

\subsection{Neighboring Condition}
\begin{figure}
\centering
{%
\resizebox*{8cm}{!}{\includegraphics{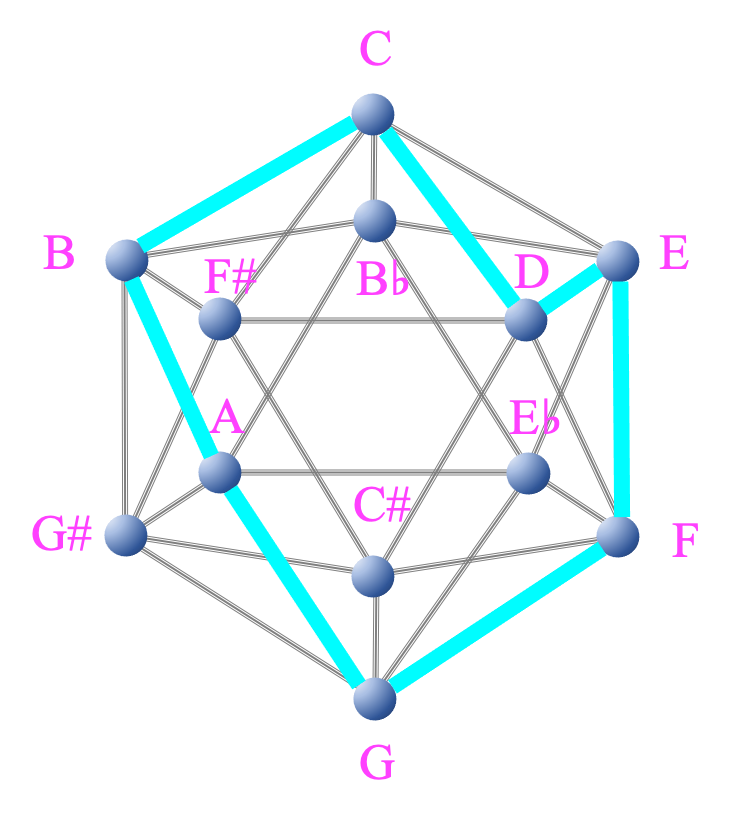}}}\hspace{5pt}
\caption{An example of a musical icosahedron that satisfies the neighboring condition for the $C$-major scale.} \label{c_major_neighbor}
\end{figure}

In this subsection, we introduce ``neighboring condition" in order to find an interesting subset of musical icosahedra. For convenience, we set the following algebra for 12 tones: $X+n$ is a tone higher than $X$ by $n$ semitones for any integer $n$ under the condition $X+12m=X$ for any integer $m$. We define the neighboring condition as follows.
\\
\\
$[{\bf Neighboring} \ {\bf Condition}]$
\\
\indent
If $X$ is a neighboring tone of $Y$ in a given scale, the vertex corresponding to $X$ must be a neighborhood of the vertex corresponding to $Y$ in a given musical icosahedron\footnote{The vertex $A$ is a ``neighborhood" of $B$ if $A$ is linked to $B$.}.
\\
\\
\indent
Mathematically, this condition can be constructed by using graph geometry. For example, the musical icosahedron shown in Fig.~\ref{c_major_neighbor} satisfies the neighboring condition for the $C$-major scale (Appendix A.1 (i)).

\begin{figure}
\centering
{%
\resizebox*{8cm}{!}{\includegraphics{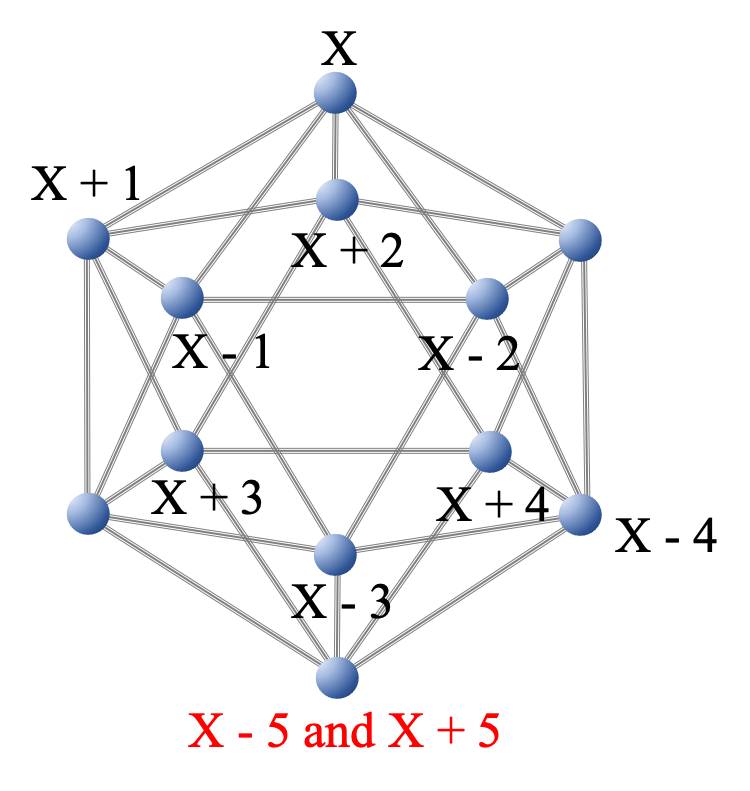}}}\hspace{5pt}
\caption{Proof of the prohibition lemma.} \label{proof_lemma}
\end{figure}

First, we tried to find a musical icosahedron that satisfies the neighboring condition for all the major or minor scales. However, we found this to be impossible by the following lemma.
\\
\\
$[{\bf Prohibition} \ {\bf Lemma}]$
\\
\indent
There are no musical icosahedra such that for any tone $X$, the vertex corresponding to $X+1$, $X-1$, $X+2$ and $X-2$ are included in the neighborhood of the vertex corresponding to $X$.
\\
\\
\indent
This lemma can be proved simply (Fig.~\ref{proof_lemma}). First, set a tone $X$ on the top. Then, $X-1$, $X-2$, $X+1$, $X+2$ must be set on the pentagon just below the top and (i) the vertex corresponding to $X-1$ must be linked with the vertex corresponding to $X+1$, $X-2$, (ii) the vertex corresponding to $X-2$ must be linked with the vertex corresponding to $X-1$, (iii) the vertex corresponding to $X+1$ must be linked with the vertex corresponding to $X+2$, the vertex corresponding to $X-1$, (iv) the vertex corresponding to $X+2$ must be linked with the vertex corresponding to $X+1$. Although there are many arrangements satisfying the condition, one can choose the arrangement in the figure (Fig.~\ref{proof_lemma}) without loss of generality. Then, a position of $X+3$ $(X-3)$ is uniquely determined by the position of $X+1$ $(X-1)$ and $X+2$ $(X-2)$. In addition, a position of $X+4$ $(X-4)$ is uniquely determined by the position of $X+2$ $(X-2)$ and $X+3$ $(X-3)$. Then, a position of $X+5$ $(X-5)$ is uniquely determined by the position of $X+3$ $(X-3)$ and $X+4$ $(X-4)$, and the position of $X+5$ is the same as the position of $X-5$. However, $X-5$ is not equal to $X+5$ because we deal with the 12 tone system.

\begin{figure}
\centering
{%
\resizebox*{14cm}{!}{\includegraphics{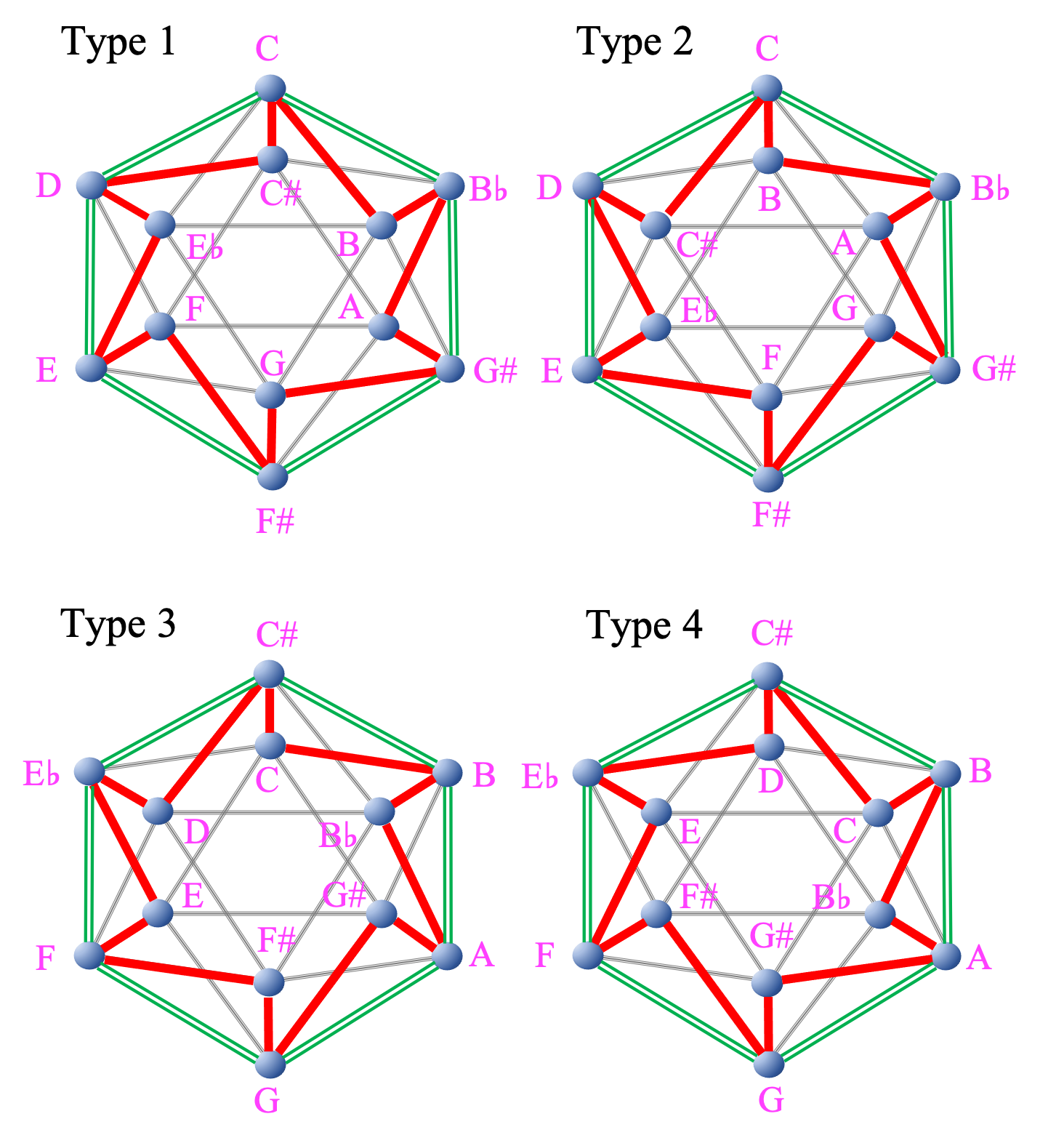}}}\hspace{5pt}
\caption{Four types of the chromatic/whole tone musical icosahedra. The red lines show the chromatic scale and the green lines show one of the whole tone scales.} \label{chromatic_whole tone}
\end{figure}

In order that a musical icosahedron satisfies the neighboring condition for all the major or minor scales, the musical icosahedron must satisfy the assumption of the lemma, i.e. for any tone $X$, the vertex corresponding to $X+1$, $X-1$, $X+2$ and $X-2$ are included in the neighborhood of the vertex corresponding to $X$. For example, the vertex corresponding to $C$ must be included in a neighborhood of the vertex corresponding to $B$, $D$ because one of them is included in a neighborhood of the others in $C$-major, and the vertex corresponding to $C$ must be included in a neighborhood of the vertex corresponding to $B\flat$, $C\sharp$ because one of them is included in a neighborhood of the others in $G\sharp$-major. Then, there are no musical icosahedra that satisfy the neighboring condition for all the major/minor scales.
\\
\indent
Then, we also introduce a condition reflecting topology of scales, distant condition as follows.
\\
\\
$[{\bf Distant} \ {\bf Condition}]$
\\
\indent
If $X$ is a neighboring tone of $Y$ in a given scale, the vertex corresponding to $X$ must not be a neighborhood of the vertex corresponding to $Y$ in a given musical icosahedron.
\\
\\
\indent
We found that there exist only four musical icosahedra that can satisfy the neighboring condition for the chromatic scale, and the neighboring condition for one of the whole tone scales and the distant condition for the other whole tone scale up to the musical isomorphism (Appendix A.8). We named them ``chromatic/whole tone musical icosahedra" (Fig.~\ref{chromatic_whole tone}). Considering the number of musical icosahedra is around 4 million up to the musical isomorphism, it is surprising that there are only four chromatic/whole tone musical icosahedra. Note that there are no musical icosahedra that can satisfy the neighboring condition for the chromatic scale and the two whole tone scales because of the prohibition lemma. 
\\
\indent
First, one can prove that there exist only two musical icosahedra satisfying the neighboring condition for the chromatic scale and the neighboring condition for the whole tone scale including $C$ and the distant condition for the whole tone scale including $C\sharp$, noting that all the vertices corresponding to the elements of $C$, $C\sharp$, $D$ must be linked with each other, and the same condition must be satisfied for $D$, $E\flat$, $E$, and $E$, $F$, $F\sharp$, and $F\sharp$, $G$, $G\sharp$, and $G\sharp$, $A$, $B\flat$, and $B\flat$, $B$, $C$. Note that all the triplets listed here begin with a tone constituting the whole tone scale including $C$. 
\\
\indent
Similarly, one can prove that there exist only two musical icosahedra satisfying the neighboring condition for the chromatic scale and the neighboring condition for the whole tone scale including $C\sharp$ and the distant condition for the whole tone scale including $C$.


Last, we analyze mathematical structure of transformations that change types of the chromatic/whole tone musical icosahedra. We define $Y_{nm}$ as a transformation from the type $n$ to $m$. $Y_{13(24)}$ raise all the tones in the type 1 (type 2) by one semitone, and $Y_{31(42)}$ is the inverse transformation of $Y_{13(24)}$. Because $Y_{13(24)}^2$ is a symmetry transformation of the icosahedron (in other words, $Y_{13(24)}$ is a ``root" of a symmetry transformation of the icosahedron that is not a symmetry transformation), if one adds $Y_{13}$, $Y_{31}$, $Y_{24}$, $Y_{42}$ to the symmetry transformations of the icosahedron, the resulting set has the group structure. One may consider the resulting set as the ``spin" group with respect to $A_5\times Z_2$\footnote{Spin is one of the fundamental properties of the electron in quantum physics, and usually discussed as a double-cover of $SO(n)$.}. In other words, the type 1 (type 4) chromatic/whole tone musical icosahedron can be considered as the ``spin up" (``spin down") chromatic/whole tone musical icosahedron.

On the other hand, if one adds $Y_{12}$, $Y_{21}$, $Y_{34}$, $Y_{43}$ to symmetry transformations of the icosahedron, the resulting set does not have the group structure. This is because acting $Y_{12}$ twice to the type 1 makes a musical icosahedron that does not satisfy the neighboring condition for the chromatic scale and one of the whole tone scales. Note that $Y_{12}$ is a ``six-fold rotation" for the inner 6 tones: $C\sharp \to E\flat \to F \to G \to A \to B \to C\sharp$. Analyzing the mathematical structure of $Y_{mn}$ may be an interesting topic in mathematics.

\subsection{Fundamental properties of Chromatic/Whole Tone Musical Icosahedra}

\begin{figure}
\centering
{%
\resizebox*{7cm}{!}{\includegraphics{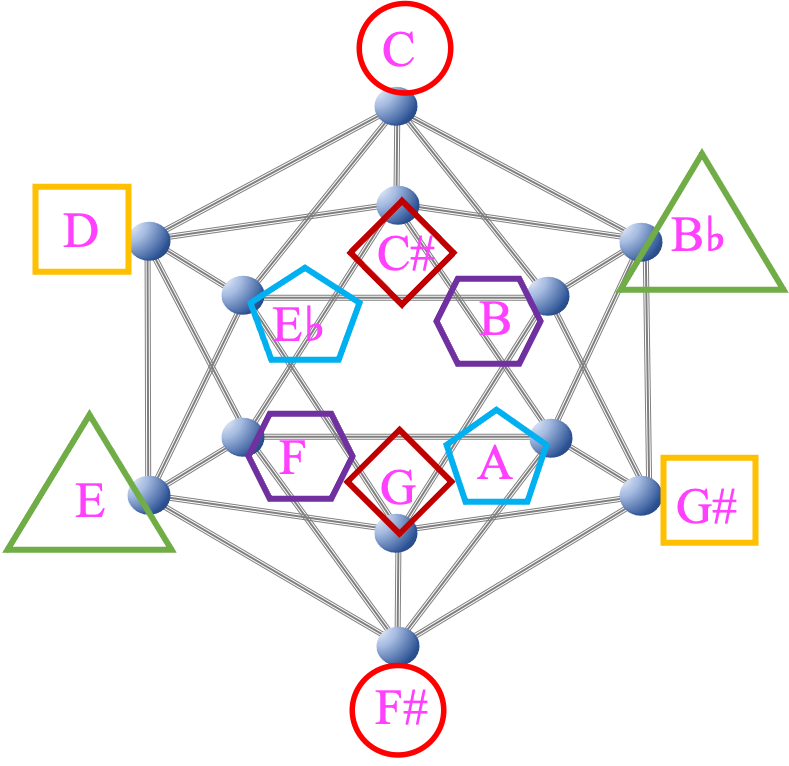}}}\hspace{5pt}
\caption{Tritones in the type 1 chromatic/whole tone musical icosahedron.} \label{tritone_chromatic_whole tone}
\end{figure}

Chromatic/whole tone musical icosahedra have many beautiful properties. First, let us see the properties common to all the four types.
\\
\\
\indent
(i) Any two tones that make tritones are set on the vertices facing each other. 
\\
\indent
(ii) One of the whole tone scales is set on the regular hexagon, and the other is set on the regular hexagram in the projection of the musical icosahedra onto the paper.
\\
\\
\indent
One can easily check (i) (Fig.~\ref{tritone_chromatic_whole tone}) and (ii) (Fig.~\ref{chromatic_whole tone}). The property (i) leads to the following strange correspondence. All the Messiaen's modes of limited transposition constructed by four tones ($C$, $E\flat$, $F\sharp$, $A$ (that can be transposed in 3 ways), or $C$, $C\sharp$, $F\sharp$, $G$ (that can be transposed in 6 ways), or $C$, $D$, $F\sharp$, $G\sharp$ (that can be transposed in 6 ways), or the modes obtained by raising the above three modes by some semitones) are set on the golden rectangle\footnote{A golden rectangle is a rectangle that satisfies the condition that ratio between length of the long edge and the short edge is given by the golden ratio.}.
 \\
\indent

Then, we see properties common to the type 1 and the type 2. 
\\
\\
\indent
(iii) The whole tone scale including $C$ is set on the regular hexagon in the projection of the musical icosahedra onto the paper, and the whole tone scale including $C\sharp$ is set on the regualr hexagram in the projection of the musical icosahedra onto the paper.
\\
\indent
(iv) $X$, $X+1$, $X+2$ ($X$ is $C$, $D$, $E$, $F\sharp$, $G\sharp$, or $B\flat$ that are tones constituting the whole tone scale including $C$) are set on a triangle made by edges of the icosahedron.
\\
\\
\indent
Next, we see properties common to the type 3 and the type 4. These properties are a kind of dual of the properties common to the type 1 and the type 2. 
\\
\\
\indent
(v) The whole tone scale including $C\sharp$ is set on the regular hexagon in the projection of the musical icosahedra onto the paper, and the whole tone scale including $C$ is set on the regular hexagram in the projection of the musical icosahedra onto the paper.
\\
\indent
(vi) $X$, $X+1$, $X+2$ ($X$ is $C\sharp$, $E\flat$, $F$, $G$, $A$, or $B$ that are tones constituting the whole tone scale including $C\sharp$) are set on a triangle made by edges of the icosahedron.
\\
\\
\indent
Then, we see properties common to the type 1 and the type 4.
\\
\\
\indent
(vii) The figure made by red lines (Fig.~\ref{chromatic_whole tone}).
\\
\\
\indent
Note that the type 4 is obtained by raising all the tones of the type 1 by one semitone.
\\
\\
\indent
Next, we see properties common to the type 2 and the type 3. These properties are a kind of dual of the properties common to the type 1 and the type 4. 
\\
\\
\indent
(viii) The figure made by red lines (Fig.~\ref{chromatic_whole tone}).
\\
\\
\indent
Note that the type 3 is obtained by raising all tones of the type 2 by one semitone.

\subsection{Major/Minor Triad and Golden Triangle}

\begin{figure}
\centering
{%
\resizebox*{16cm}{!}{\includegraphics{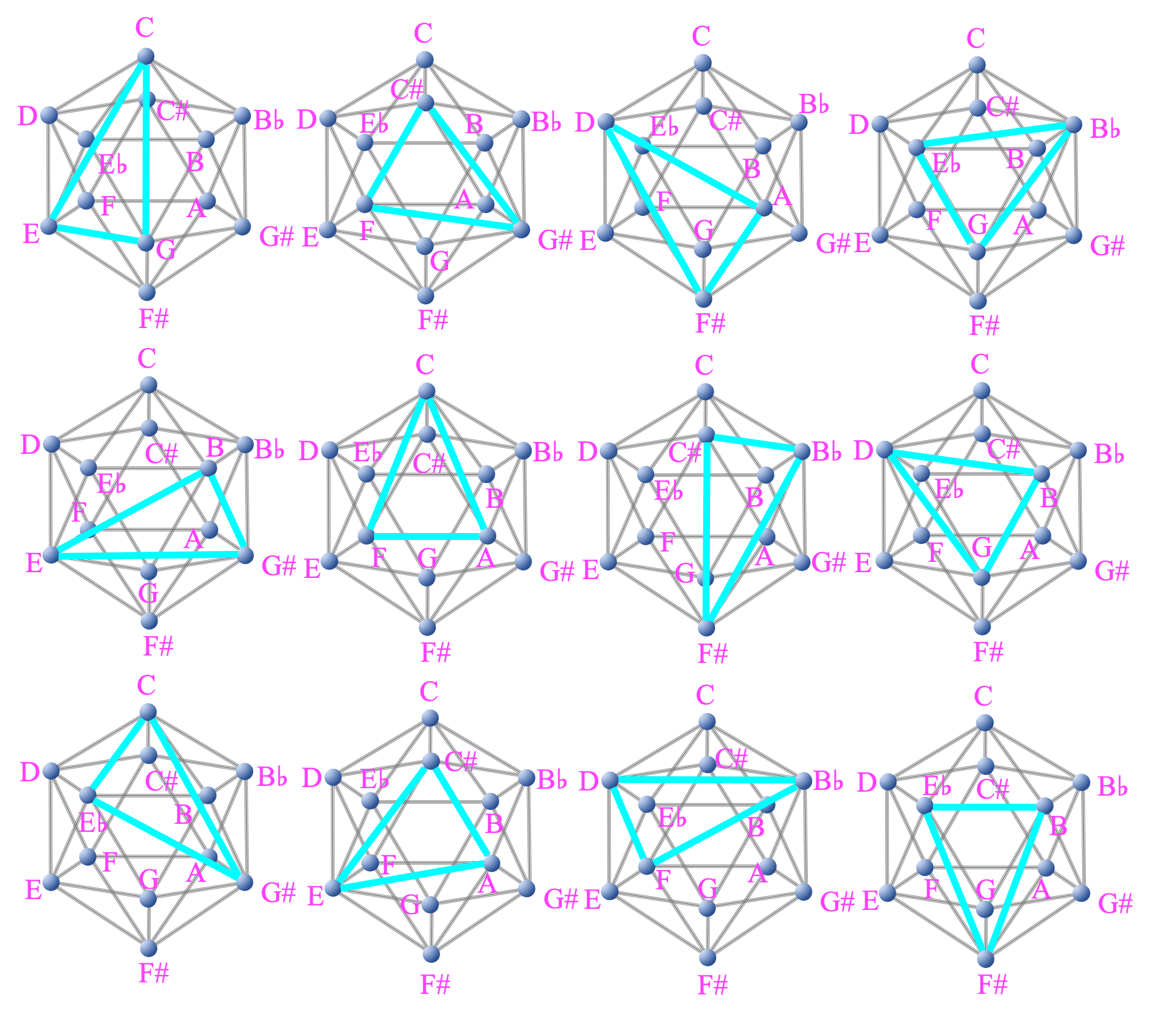}}}\hspace{5pt}
\caption{The golden triangles corresponding to the major triads in the type 1 of the chromatic/whole tone musical icosahedra.} \label{golden}
\end{figure}

Then, we move on to the main point of this section: arrangement of the major/minor triads on the four types of the chromatic/whole tone musical icosahedra. We discovered the following beautiful theorem.
\\
\\
$[{\bf Golden} \ {\bf Theorem}]$
\\
\indent
In the type 1 and the type 4, the major triads are set on the golden triangles.
\\
\indent
In the type 2 and the type 3, the minor triads are set on the golden triangles.
\\
\\
\indent
Figure~\ref{golden} shows the 12 golden triangles corresponding to the major triads in the type 1 of the chromatic/whole tone musical icosahedra. Similarly, one can check that the type 4 also includes the 12 golden triangles corresponding to the major triads, and the type 2 and the type 3 includes 12 golden triangles corresponding to the minor triads (Fig.~\ref{chromatic_whole tone}).

We remark that this discussion does not involve frequency of tones but uses only the golden triangles in the icosahedron and topology of the artificial scales (the chromatic scale and the whole tone scales) on musical icosahedra through the neighboring/distant condition. Also, the duality between the type 1, type 4 and the type 2, type 3 is none other than the duality between the major triads and the minor triads, i.e., we derived the duality between the major (bright) triads and the minor (dark) triads mathematically.
\\
\indent
Note that the number of golden triangles in the icosahedron is 60. Therefore, one can generalize major/minor triads. We discuss this point later. Also, one can uniquely determine the major/minor triads because all the major (minor) triads are set on the golden triangles on the type 1 and the type 4 (the type 2 and the type 3). If one transposes the other triads corresponding to the golden triangle, some of the resulting triads are not set on the golden triangle, e.g., $CGG\sharp$ is set on the golden triangle in the type 1 and $C\sharp G\sharp A$ is not set on the golden triangle in the type 1.
\\
\\
$[{\bf Uniqueness}$  ${\bf of}$ ${\bf Major/Minor}$ ${\bf Triads}]$
\\
\indent
For any tone $X$, triads $X$, $X+n$, $X+m$ $(m>n)$ are set on the golden triangles in the type 1 and the type 4 [the type 2 and the type 3] under the condition that the apex of the golden triangles corresponding to $X$, $X+n$, $X+m$ is set on the regular hexagon in the projection of the musical icosahedra onto the paper, if and only if $n=4$, $m=7$ [$n=3$, $m=7$].
\\
\\
Replacing the condition ``the apex of the golden triangles corresponding to $X$, $X+n$, $X+m$ is set on the regular hexagon in the projection of the musical icosahedra onto the paper" by the condition ``the apex of the golden triangle corresponding to $X$, $X+n$, $X+m$ is set on the regular hexagram in the projection of the musical icosahedra onto the paper", one has the following triads uniquely: $X$, $X+7$, $X+10$ by the type 1 and the type 4, $X$, $X+7$, $X+9$ by the type 2 and the type 3. The former/latter triads should be regarded as the dual triads of the major/minor triads. Note that these triads can be obtained by replacing the middle tone of the major/minor triads by its tritone.

We found $X$, $X+7$, $X+10$ ($X$, $X+7$, $X+9$) are the fundamental triads for the hexatonic minor/major scale (Appendix A.3, A.4). Recall that the major/minor scale is obtained by combining three major/minor triads, e.g., $C$, $E$, $G$ ($C$-major triad), and $G$, $B$, $D$ ($G$-major triad), and $F$, $A$, $C$ ($F$-major triad) make $C$, $D$, $E$, $F$, $G$, $A$, $B$, $C$ ($C$-major scale), $C$, $E\flat$, $G$ ($C$-minor triad), and $G$, $B\flat$, $D$ ($G$-minor triad), and $F$, $A\flat$, $C$ ($F$-minor triad) make $C$, $D$, $E\flat$, $F$, $G$, $A\flat$, $B\flat$, $C$ ($C$-minor scale). By combining $C$, $G$, $B\flat$ (the fundamental hexatonic minor triad whose base tone is $C$), and $G$, $D$, $F$ (the fundamental hexatonic minor triad whose base tone is $G$), and $F$, $C$, $E\flat$ (the fundamental hexatonic minor triad whose base tone is $F$), one has $C$, $D$, $E\flat$, $F$, $G$, $B\flat$ that is nothing but the hexatonic minor scale whose base tone is $C$. Also, by combining $C$, $G$, $A$ (the fundamental hexatonic major triad whose base tone is $C$), and $G$, $D$, $E$ (the fundamental hexatonic major triad whose base tone is $G$), and $F$, $C$, $D$ (the fundamental hexatonic major triad whose base tone is $F$), one has $C$, $D$, $E$, $F$, $G$, $A$ that is the hexatonic major scale whose base tone is $C$. By combining $C$, $G$, $B\flat$ and $F$, $C$, $E\flat$, one has $C$, $E\flat$, $F$, $G$, $B\flat$ that is the pentatonic minor scale whose base tone is $C$. By combining $C$, $G$, $A$ and $G$, $D$, $E$, one has $C$, $D$, $E$, $G$, $A$ that is the pentatonic major scale whose base tone is $C$ (Appendix A.5, A.6).
\\
\\
\indent
Note that this derivation of the pentatonic major/minor scales may be not beautiful because there may be no reasons why we choose $C$, $G$, $B\flat$ ($C$, $G$, $A$) and $F$, $C$, $E\flat$ ($G$, $D$, $E$), and do not choose $C$, $G$, $B\flat$ ($C$, $G$, $A$), and $G$, $D$, $F$ ($F$, $C$, $D$). The latter combination makes a scale, $C$ $D$, $F$, $G$, $B\flat$ ($C$ $D$, $F$, $G$, $A$). Although these pentatonic scales and the major/minor pentatonic scales correspond to four of the Chinese musical scales ($Gong$, $Shang$, $Zhi$, $Yu$), we cannnot list all the five Chinese scales ($Gong$, $Shang$, $Jue$, $Zhi$, $Yu$). The remarkable point is that $Jue$ is given by $C$, $E\flat$, $F$, $A\flat$, $B\flat$ and does not include $G$ (Shi C. Y. 2002). One may create a new theory that deals with these scales and/or Japanese pentatonic scales, etc. by using musical icosahedra in a different way.

\subsection{Dualities between Major Scales and Minor Scales}
We found the following general theorem.
\\
\\
$[{\bf Hexagon-Icosahedron}$ ${\bf Symmetry}]$\footnote{The regular hexagons represent the whole tone scales in a musical dodecagon that satisfies the neighboring condition for the chromatic scale.}
\\
\indent
Raising all the tones of the type $n$ chromatic/whole tone musical icosahedron by any number of whole tones is equivalent to a symmetry transformation of the icosahedron.
\\
\\
\indent
By using the above theorem, one can show that for a figure of $C_X^n$ representing any scale $X$ on the type $n$ chromatic/whole tone musical icosahedron and the figure of $C_Y^n$ obtained by raising all the tones of $X$ by some whole tones, the figure of $C_X^n$ can be superposed on the figure of $C_Y^n$ by a symmetry transformation of the icosahedron to the figure of $C_X^n$. If one draws the $X$-major/minor scale on the type $n$ chromatic/whole tone musical icosahedron, one can superpose the figure of the $X+2$-, $X+4$-, $X+6$-, $X+8$-, $X+10$-major/minor scale to the figure of the $X$-major/minor scale by a symmetry transformation to the $X$-major/minor scale. Also, one can obtain the $X+1$-, $X+3$-, $X+5$-, $X+7$-, $X+9$-, $X+11$-major/minor scale by transforming a figure representing the $X$-, $X+2$-, $X+4$-, $X+6$-, $X+8$-, $X+10$-major/minor scale from the type 1, type 2, type 3, type 4 to the type 4, type 3, type 2, type 1.

\begin{figure}
\centering
{%
\resizebox*{14cm}{!}{\includegraphics{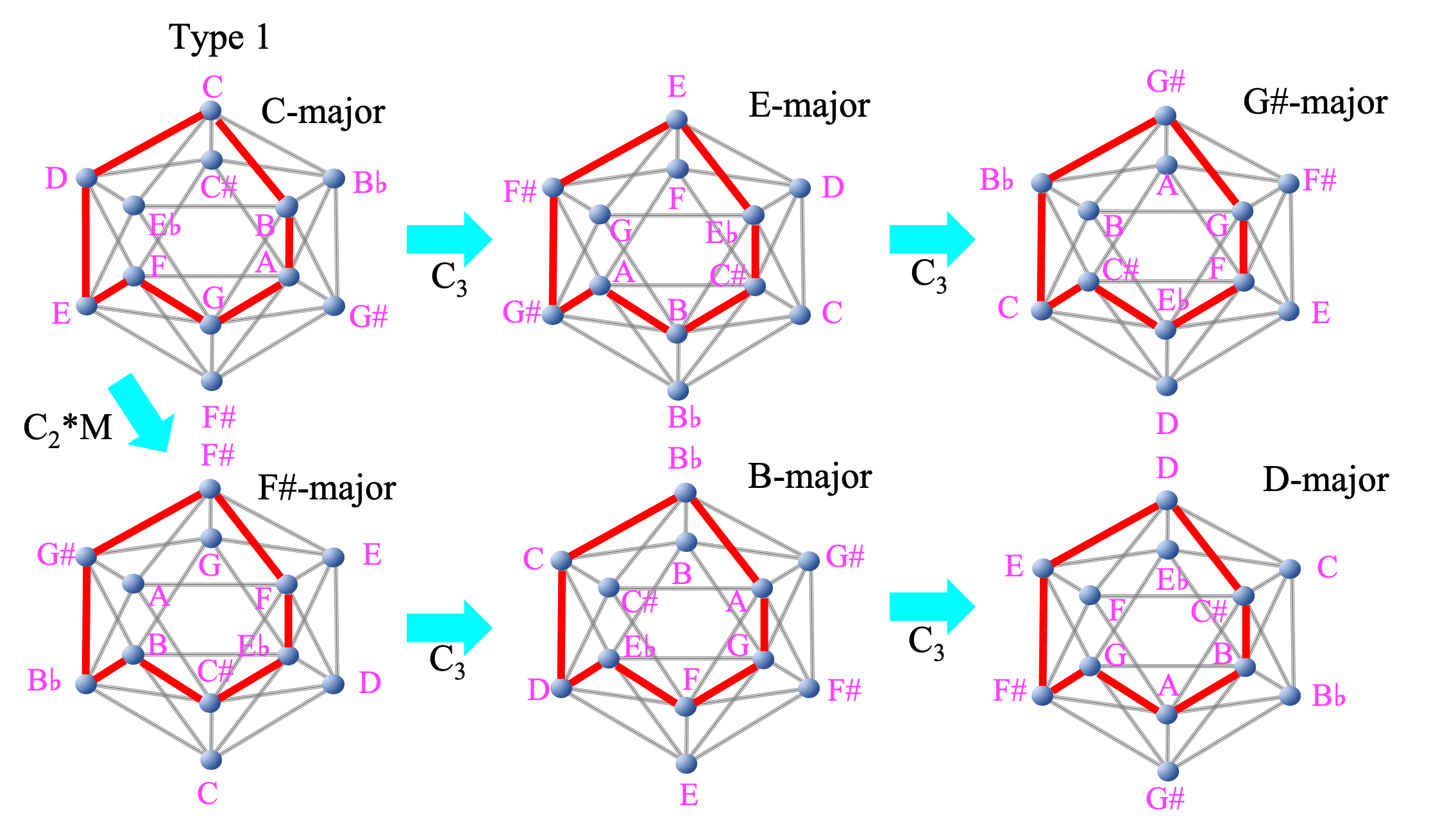}}}\hspace{5pt}
\caption{The $D$-major, $E$-major, $F\sharp$-major, $G\sharp$-major, $B\flat$-major scale are obtained by a symmetry transformation of the icosahedon for the $C$-major scale.} \label{major_type1}
\end{figure}

For example, the $C$-major scale is transformed to the $D$-major, $E$-major, $F\sharp$-major, $G\sharp$-major, $B\flat$-major scale by using a symmetry transformation of the icosahedron as shown in Fig.~\ref{major_type1}.

\begin{figure}
\centering
{%
\resizebox*{14cm}{!}{\includegraphics{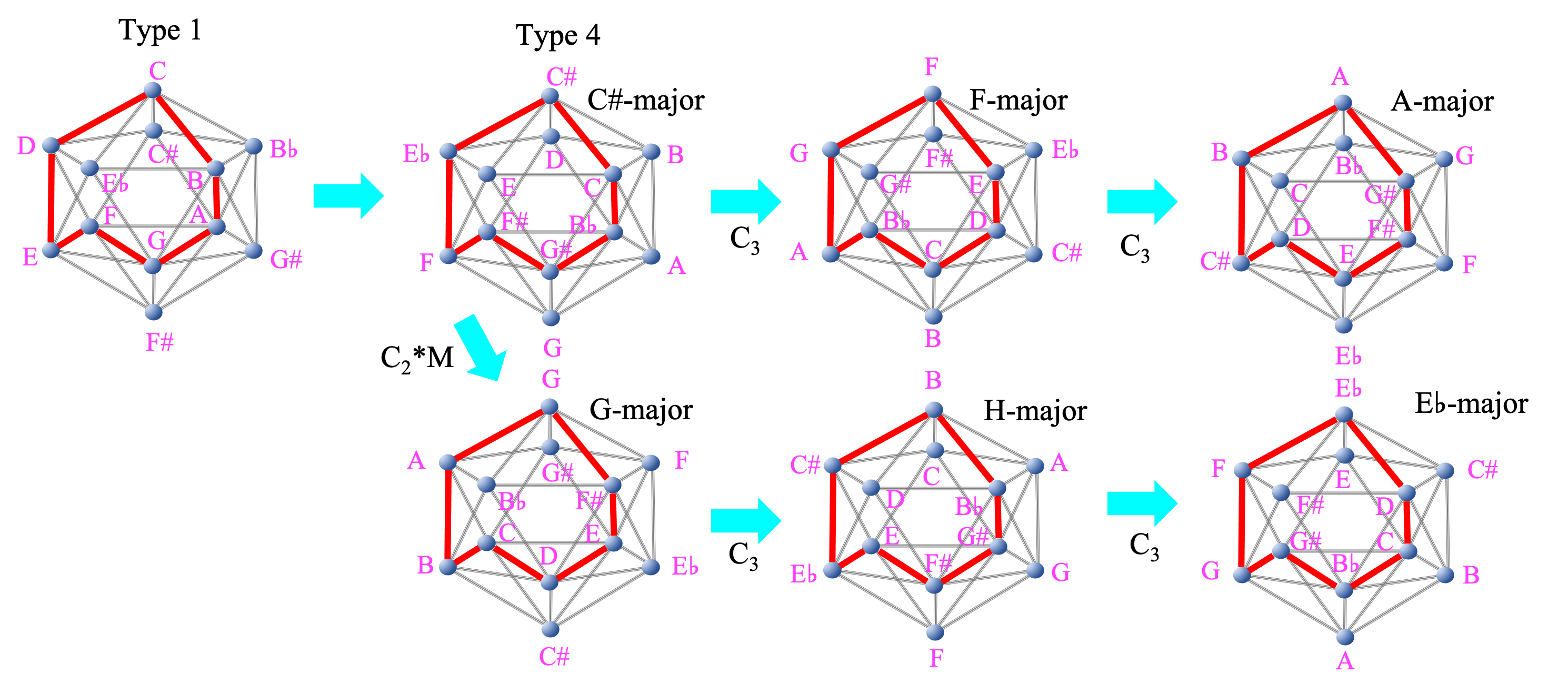}}}\hspace{5pt}
\caption{The $C\sharp$-major scale is obtained by the type change from the type 1 to the type 4 for the $C$-major scale on the type1, and the $E\flat$-major, $F$-major, $G$-major, $A$-major, and $B$-major scales are obtained by a symmetry transformation of the icosahedron for the $C\sharp$-major scale.} \label{major_type2}
\end{figure}

In addition, the $C$-major scale is transformed to the $C\sharp$-major by the type change from the type 1 to the type 4, and the $C\sharp$-major scale is transformed to the $E\flat$-major, $F$-major, $G$-major, $A$-major, and $B$-major scale by a symmetry transformation of the icosahedron as shown in Fig.~\ref{major_type2}.

\begin{figure}
\centering
{%
\resizebox*{14cm}{!}{\includegraphics{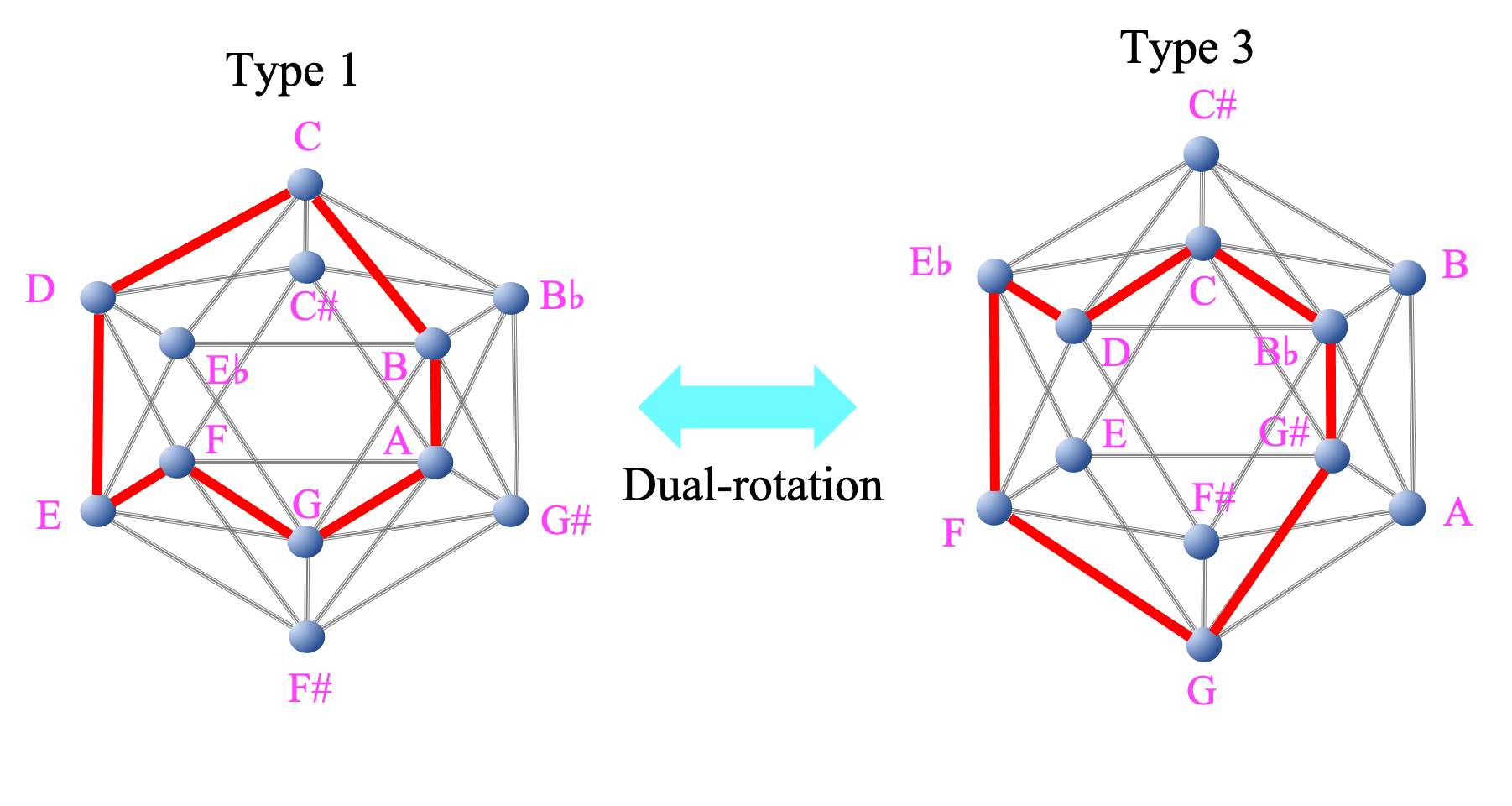}}}\hspace{5pt}
\caption{The $C$-minor scale is obtained by transforming the $C$-major scale from the type 1 to the type 3 and applying the two-fold rotation.} \label{two_fold_dual}
\end{figure}

Furthermore, the $X$-major/minor scale on the type 1 (type 2) is obtained by $A_{3(4)}C_2A_{1(2)}^{-1}$ to the $X$-minor/major on the type 3 (type 4) where $A_n$ is the symmetry transformation that changes the $C$-major/minor scale to the $X$-major/minor scale on the type $n$. Figure~\ref{two_fold_dual} shows the transformation from the $C$-major scale on the type 1 to the $C$-minor scale on the type 3. This shows the duality between the major (bright) scales and the minor (dark) scales.

\subsection{Duality between Gregorian modes and Major/Minor Scales}
\begin{figure}
\centering
{%
\resizebox*{14cm}{!}{\includegraphics{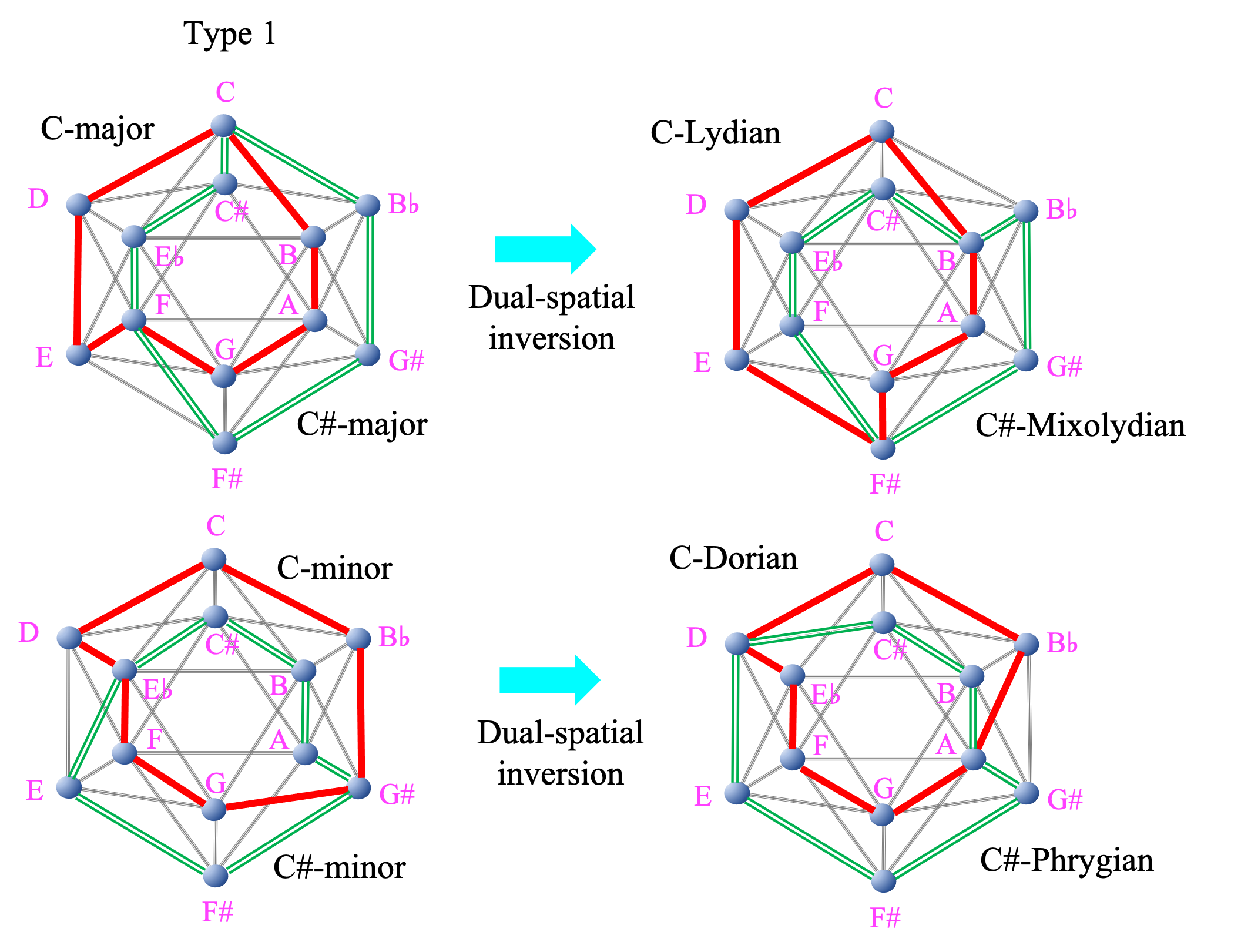}}}\hspace{5pt}
\caption{An example of the correspondence between the major/minor scales and the Gregorian modes.} \label{c_major_gregorian}
\end{figure}

One may perceive that this theory is not beautiful because transposition of the major scales is not closed in one type of chromatic/whole tone musical icosahedra, and the figure of the $C$-major scale on the type 1 is not the same as the $C\sharp$-major scale on the type 1. However, the $C$-major scale on the type 1 is similar to the $C\sharp$-major scale on the type 1, and this slight difference is nothing but the slight difference between the major/minor scales and the Gregorian modes (Appendix 7). That is to say, the fact that there exist four chromatic/whole tone musical icosahedra indicates not only the duality between the major triads/scales and minor triads/scales, but also the duality between the major/minor scales and the Gregorian modes, by lack of beauty with respect to the transposition of the major/minor scales. We found that all the Gregorian modes can be obtained by using the spatial inversion\footnote{The spatial inversion for a scale $X$ transforms all the edges $a$ constructing $X$ to the edge facing to $a$.} for all the major scales and minor scales on the type 1 and the type 4 (or the type 2 and the type 3) of chromatic/whole tone musical icosahedra. If one chooses the type 1, then, one has the $C\sharp$-Myxolydian, $C$-Lydian, $C\sharp$-Phrygian, and $C$-Dorian by the spatial inversion to the $C$-major, $C\sharp$-major, $C$-minor, and $C\sharp$-minor scale (Fig.~\ref{c_major_gregorian}). Note that the spatial inversion is one of the most elementary operations in mathematics or physics, e.g. the electric field changes its sign by the spatial inversion while the magnetic field does not change its sign by the spatial inversion.
\\
\\
\indent

$[{\bf Duality}$ ${\bf between}$ ${\bf Major/Minor}$ ${\bf Scales}$ ${\bf and}$ ${\bf Gregorian}$ ${\bf Modes}]$
\\
\indent
All the Gregorian modes can be obtained by acting the spatial inversion to all the major/minor scales on the type 1 and the type 4 (or the type 2 and the type 3) of the chromatic/whole tone musical icosahedra.

\subsection{Generalization of Major/Minor Triads and Scales}
As said before, some of the golden triangles in the chromatic/whole tone musical icosahedra do not correspond to the major/minor triads but other triads. By using this fact, one can generalize major/minor triads in two ways. Note that there are two types of the golden triangles in the icosahedron. One is a set of the golden triangles whose apex is included in the regular hexagon in the projection of the musical icosahedra onto the paper. The other is a set of the golden triangles whose apex is included in the regular hexagram in the projection of the musical icosahedra onto the paper. Then, the first generalization of the major/minor triads is obtained by considering the former case. The second generalization is obtained by considering all the golden triangles.

\begin{figure}
\centering
{%
\resizebox*{14cm}{!}{\includegraphics{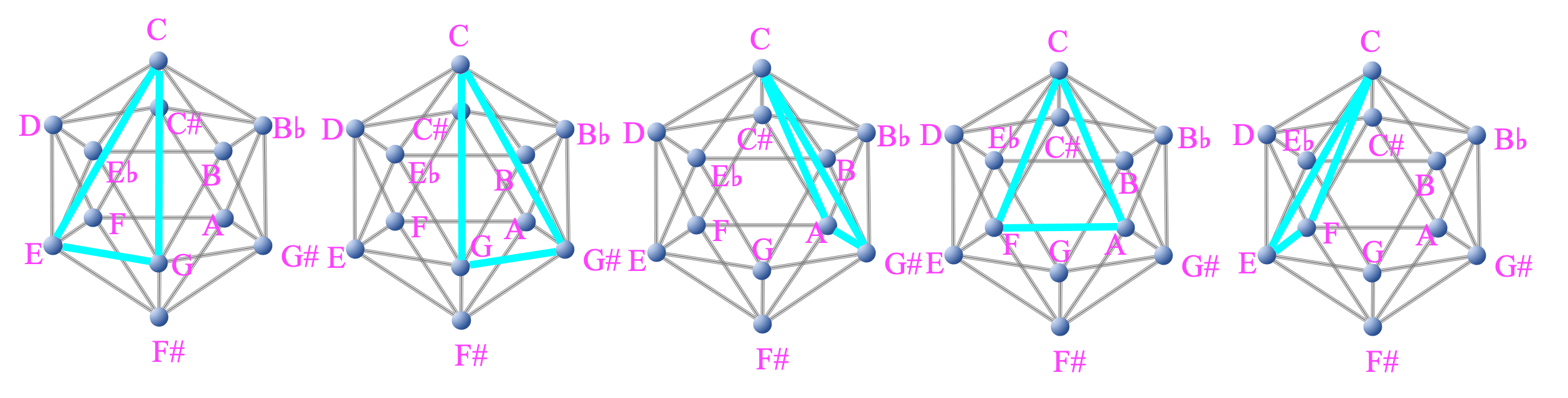}}}\hspace{5pt}
\caption{The first generalization of $C$-major triad and the corresponding golden triangles.} \label{general_major_1}
\end{figure}

Now, we list all the triads included in the first generalization of the major triads (Fig.~\ref{general_major_1}).
\\
\\
\indent
$CEG$, $CGG\sharp$, $CG\sharp A$, $CFA$, $CEF$,  
\\
\indent
and triads obtained by raising all the tones of each of the above triads by $n$ semitones ($1\leq n\leq11$).
\\
\indent

\begin{figure}
\centering
{%
\resizebox*{15cm}{!}{\includegraphics{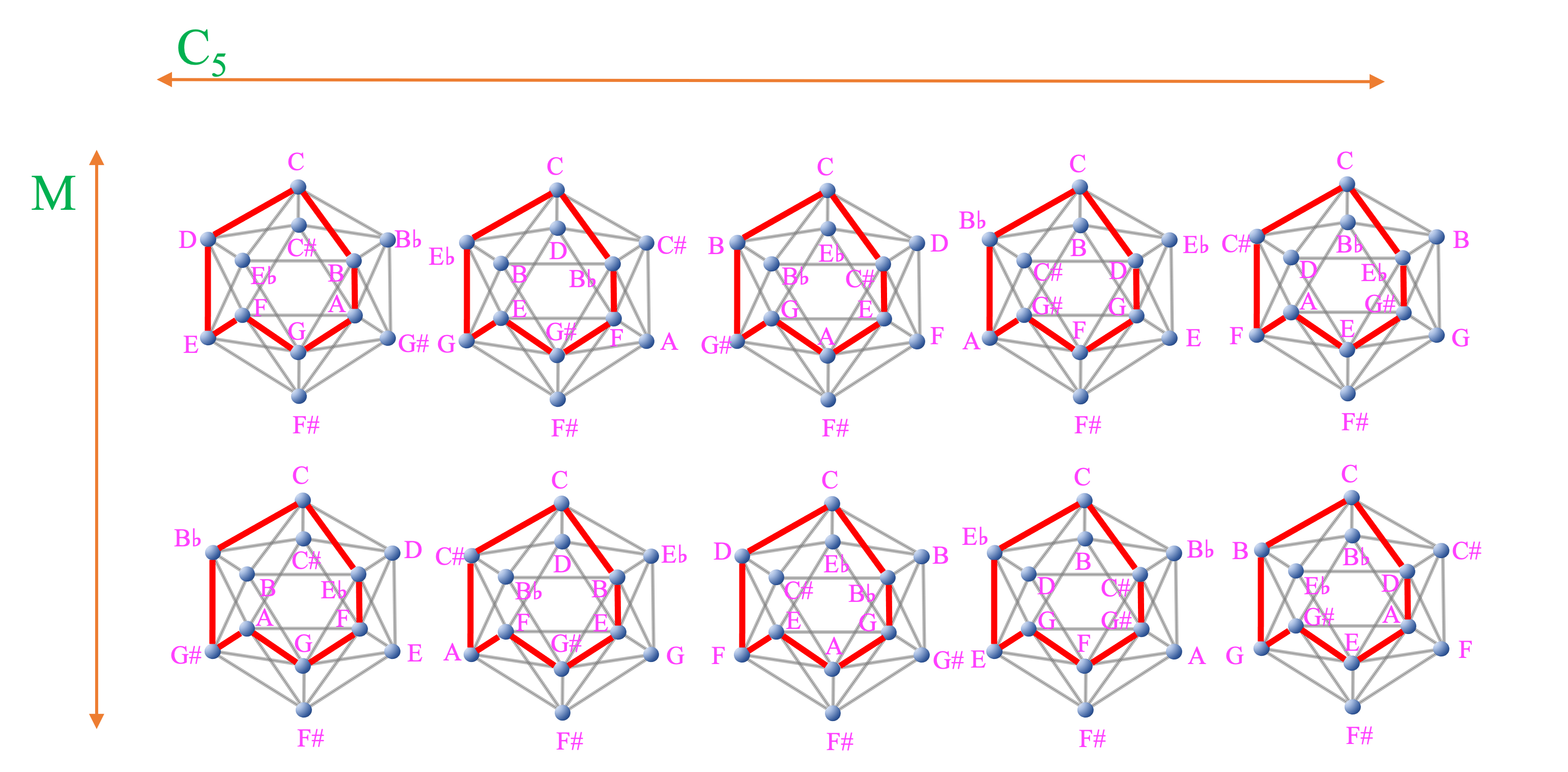}}}\hspace{5pt}
\caption{The first generalization of the major scales.} \label{general_major_scale}
\end{figure}

The first triad is a usual major triad ($CEG$), the second triad is a major seventh without the fifth ($G\sharp C (E) G$), the third triad is a minor/major seventh without the fifth ($AC(E)G\sharp$), the fourth triad is a usual major triad ($FAC$), the fifth triad is a major seventh without the third ($F(A)CE$). Next, we explain the second generalization.
\\
\\
\indent
$CEG$, $CGG\sharp$, $CG\sharp A$, $CFA$, $CEF$, $CDE\flat$, $CE\flat F$, $CFG$, $CGB\flat$, $CDB\flat$, 
\\
\indent
and triads obtained by raising all the tones of each of the above triads by $n$ semitones ($1\leq n\leq11$).
\\
\\
\indent
The sixth triad is a minor ninth without the fifth ($CE\flat G D$), the seventh triad is a major/minor seventh without the third ($F(A)CE\flat$), the eighth triad is Csus4, the ninth triad is a major/minor seventh without the triad ($C(E)GB\flat$), the tenth triad is an added-ninth without the third and the fifth ($C(EG)B\flat D$).
\\
\\
\indent
Similarly, one has the following first and second generalizations of the minor triad.
\\
\\
\indent
$CE\flat G$, $CGG\sharp$, $CFG\sharp$, $CEF$, $CE\flat E$
\\
\indent
and triads obtained by raising all the tones of each of the above triads by $n$ semitones ($1\leq n\leq12$).
\\
\\
\indent
$CE\flat G$, $CGG\sharp$, $CFG\sharp$, $CEF$, $CE\flat E$, $CDF$, $CFG$, $CGA$, $CAB\flat$, $CDB\flat$
\\
\indent
and triads obtained by raising all the tones of each of the above triads by $n$ semitones ($1\leq n\leq12$).
\\
\\
\indent
The first triad is a usual minor triad, and the third triad is a usual minor triad, the fifth triad is an augmented ninth without the fifth ($CE(G)E\flat$), the sixth triad is F6 without the third ($F(A)CD$), the eighth triad is C6 without the third ($C(E)GA$), the ninth triad is a minor ninth without the fifth ($AC(E)B\flat$). We remark that $CGG\sharp$ and $CEF$ are common to the first generalization of the minor triads, and the major triads, $CFG$, and $CDB\flat$ are common to the second generalization of the minor triads and the major triads. 
\\
\\
\indent
Next, we generalize the major scales. By using the five-fold rotation ($C_5$) and the mirror transformation ($M$), one obtains the first generalization of the major scales (see Fig.~\ref{general_major_scale} and Appendix A.9). The generalization has two meanings. The first meaning is a generalization of the major scales in terms of symmetry transformation of the icosahedron. Because we obtained the major scales by using a part of symmetry transformations of the icosahedron ($C_2*M$ and $C_3$), using $C_5$ and $M$ is a simple generalization of the major scales. The second meaning is the correspondence between the generalization of the major triad and the generalization of the major scale. First, note that the major scales are each composed of three major triads, e.g., $C$-major is composed of $CEG$ (tonic), $FAC$ (sub-dominant), $GBD$ (dominant). The generalized major scales, without a usual major scale, are also composed of two or three kinds of three generalized major triads each. Note that although a usual major scale is composed of only one kind of major triad\footnote{$GBD$ and $FAC$ is obtained by transposing $CEG$.}.

\begin{figure}
\centering
{%
\resizebox*{15cm}{!}{\includegraphics{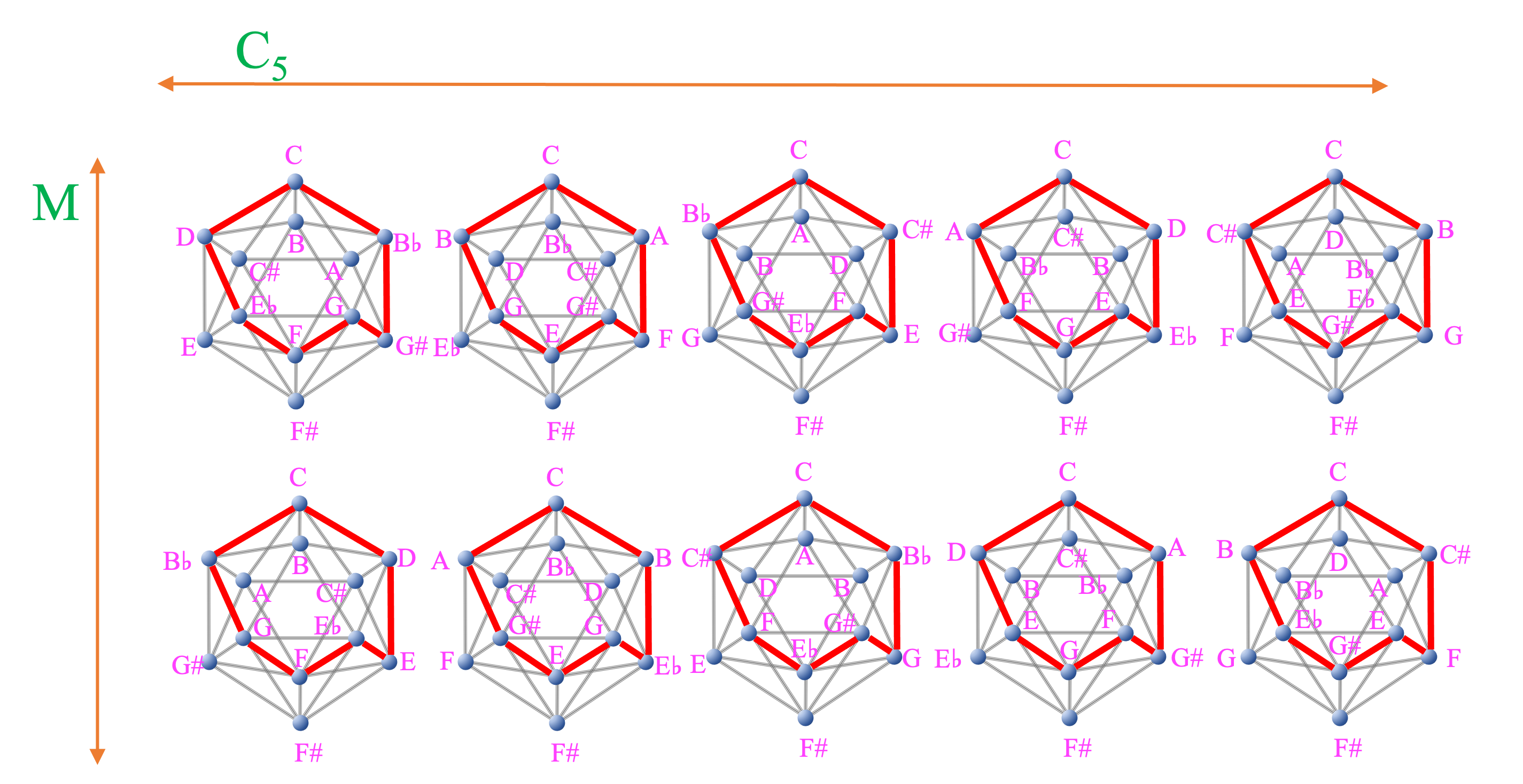}}}\hspace{5pt}
\caption{The first generalization of the minor scales.} \label{general_minor_scale}
\end{figure}

Similarly, one can generalize the minor scales (see Fig.~\ref{general_minor_scale} and Appendix A.10).

\begin{figure}
\centering
{%
\resizebox*{15cm}{!}{\includegraphics{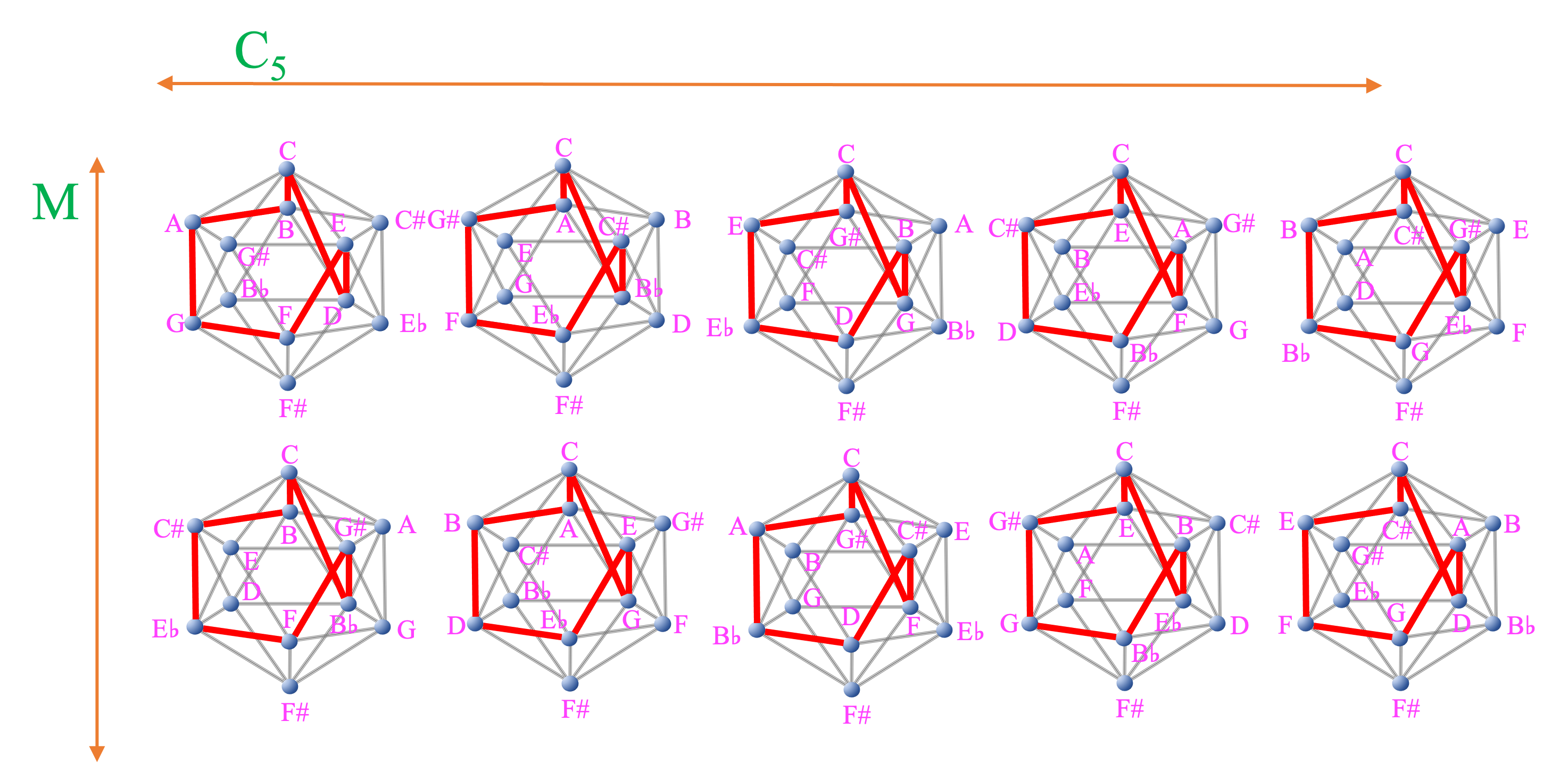}}}\hspace{5pt}
\caption{The second generalization of the major scales.} \label{general_major_scale_2}
\end{figure}

By using the transformation that moves $C$ to the top in the type 4 chromatic/whole tone musical icosahedron and doing a similar procedure as before, one has the second generalization of the major scales (see Fig.~\ref{general_major_scale_2} and Appendix A.11).

\begin{figure}
\centering
{%
\resizebox*{15cm}{!}{\includegraphics{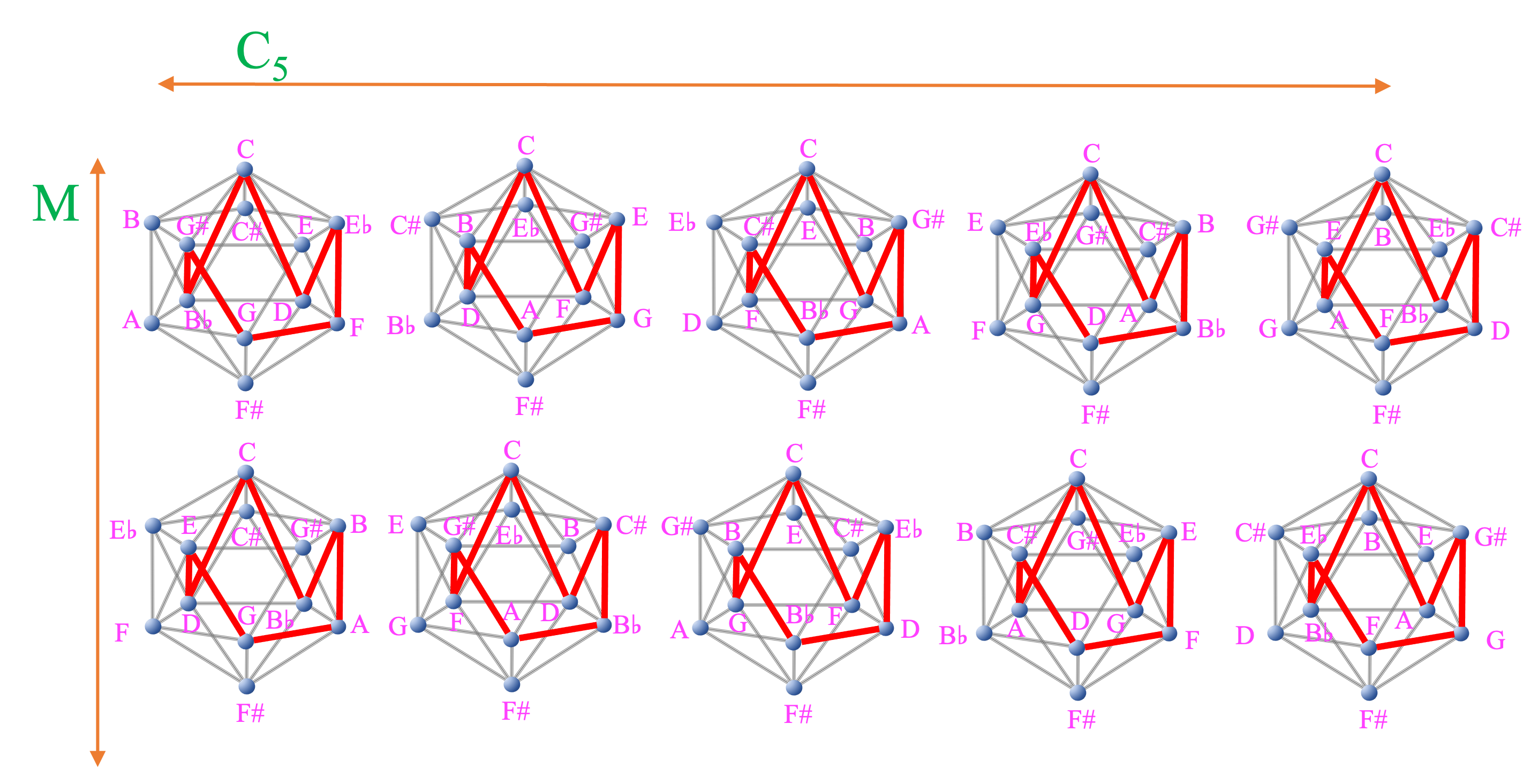}}}\hspace{5pt}
\caption{The second generalization of the minor scales.} \label{general_minor_scale_2}
\end{figure}

One also obtains the second generalization of the minor scales as follows (see Fig.~\ref{general_minor_scale_2} and Appendix A.12).

By drawing the chromatic scale on the left-top figure in Fig.~\ref{general_major_scale}, \ref{general_major_scale_2}, (Fig.~\ref{general_minor_scale}, \ref{general_minor_scale_2}) and the same figure on the other figures in Fig.~\ref{general_major_scale},  \ref{general_major_scale_2}, (Fig.~\ref{general_minor_scale}, \ref{general_minor_scale_2}), one has the generalized chromatic scales. By extracting the first tone ($C$), the third tone, the fifth tone (the fourth tone), the sixth tone, the eighth tone, the tenth tone ( the ninth tone), and the twelfth (the eleventh) tone in each generalized chromatic scale, one has the generalized major (minor) scales whose base tone is $C$\footnote{Note that a usual major scale is also constructed in the same way.}.

One can also obtain the generalization of the Gregorian modes, hexatonic major/minor scales, pentatonic major/minor scales, and any other scales in the same way.

\newpage
\section{Duality between Chromatic Scale and Pythagorean chain}
In this section, we introduce the musical icosahedra that satisfy the neighboring condition for the Pythagorean chain and the neighboring condition for one of the whole tone scales and the distant condition for the other whole tone scale: Pythagorean/whole tone musical icosahedra. Although there are many common properties between the chromatic/whole tone musical icosahedra and the Pythagorean/whole tone musical icosahedra, there exist properties that the chromatic/whole tone musical icosahedra have and the Pythagorean/whole tone musical icosahedra do not have. Instead, it is shown that there are beautiful duality relations between the chromatic scale and the Pythagorean chain through mathematical properties of the regular icosahedron.

\subsection{Four types of the Pythagorean/whole tone musical icosahedra}
\begin{figure}[H]
\centering
{%
\resizebox*{14cm}{!}{\includegraphics{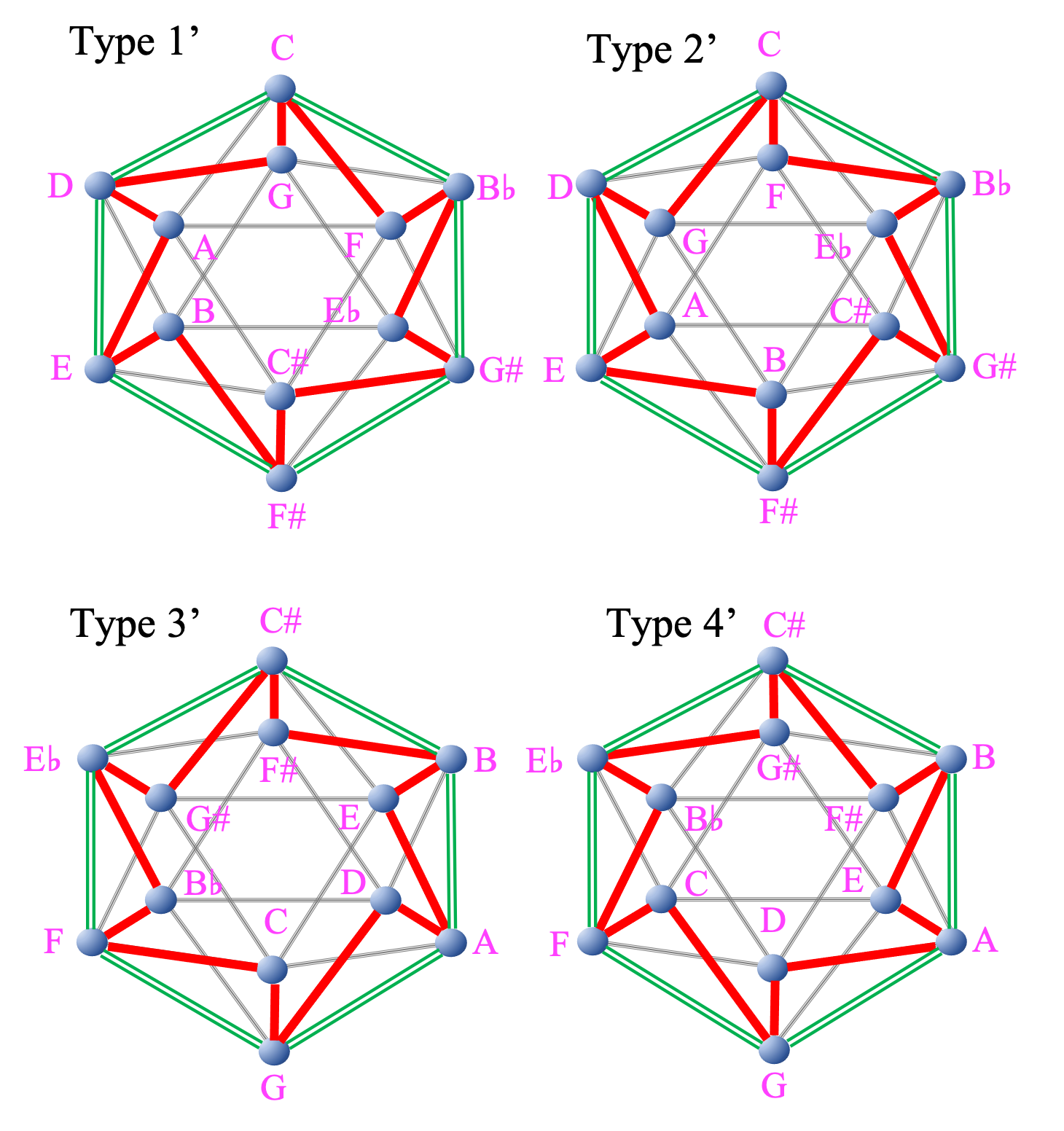}}}\hspace{5pt}
\caption{Four types of the Pythagorean/whole tone musical icosahedra.} \label{2_type}
\end{figure}

Similar to the proof of the theorem that there exist only four chromatic/whole tone musical icosahedra, one can easily show that there exist only four types of the Pythagorean/whole tone musical icosahedra (Fig.~\ref{2_type}). The red lines correspond to the Pythagorean chain and the green lines correspond to one of the whole tone scales. The type 1' and type 2' (type 3' and type 4') satisfy the neighboring condition for the Pythagorean chain and the neighboring condition for one of the whole tone scales including $C$ ($C\sharp$) and the distant condition for the whole tone scales including $C\sharp$ ($C$).

\begin{figure}[t]
\centering
{%
\resizebox*{14cm}{!}{\includegraphics{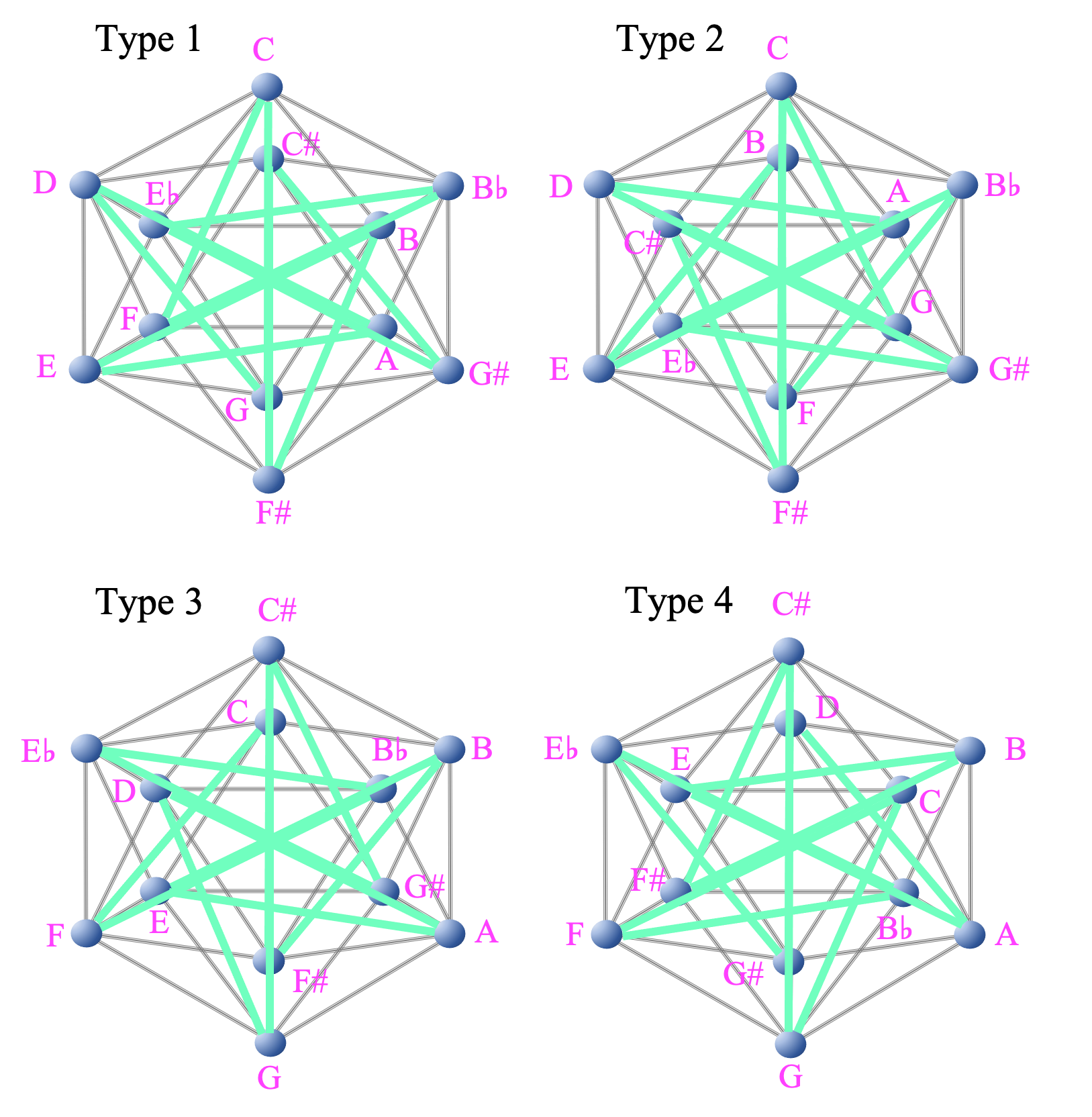}}}\hspace{5pt}
\caption{The Pythagorean chain in the chromatic/whole tone musical icosahedra.} \label{2_pythagorean}
\end{figure}

\begin{figure}[t]
\centering
{%
\resizebox*{14cm}{!}{\includegraphics{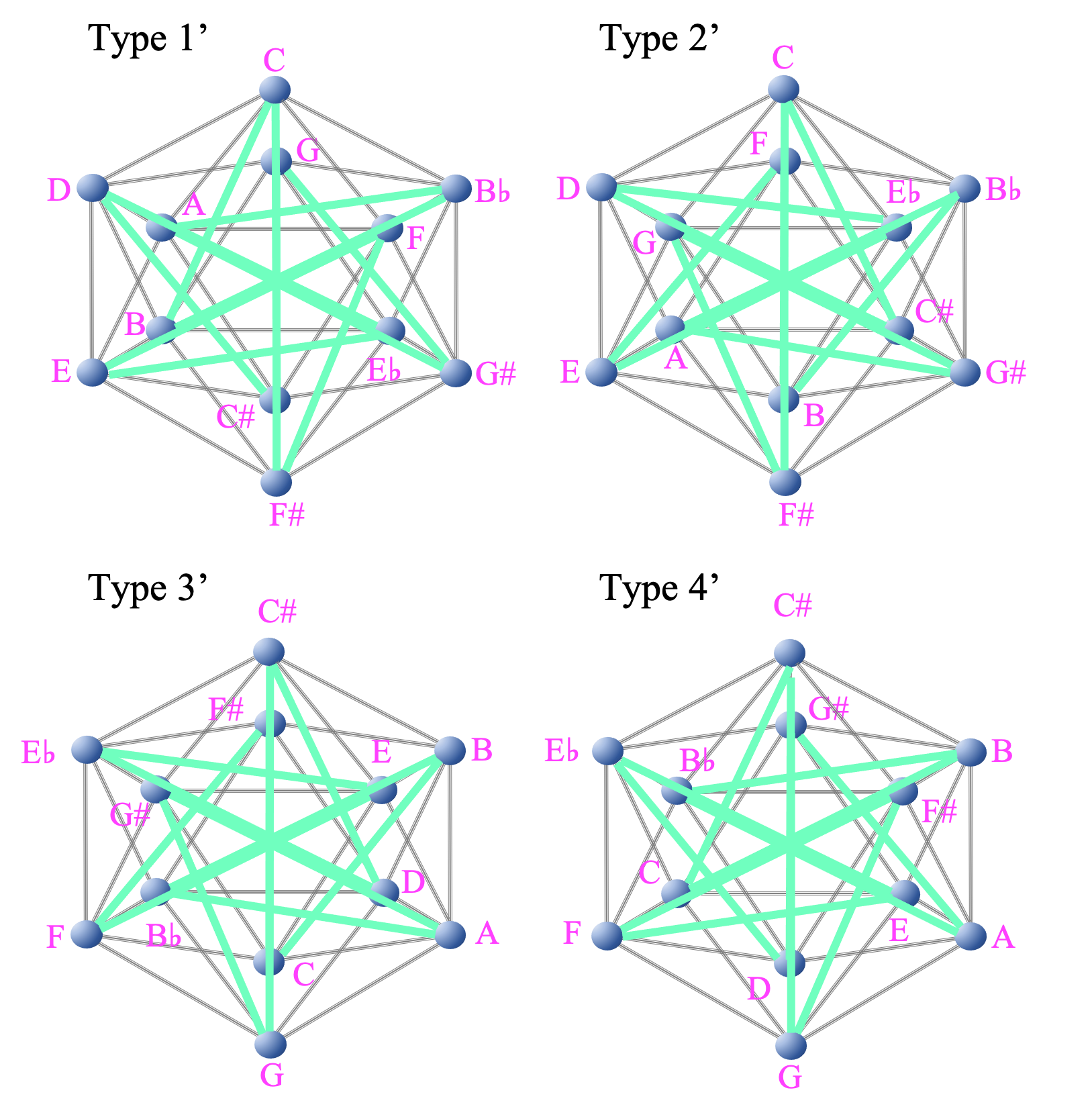}}}\hspace{5pt}
\caption{The chromatic scale in the Pythagorean/whole tone musical icosahedra.} \label{2_chromatic}
\end{figure}

Focusing on the middle line (obtained by connecting two vertices, which do not neighbor each other and are not opposite to each other), one can obtain a duality between the Pythagorean/whole tone musical icosahedra and the chromatic/whole tone musical icosahedra. The Pythagorean chain is constructed by only the middle lines in the chromatic/whole tone musical icosahedra (Fig.~\ref{2_pythagorean}) and the chromatic scale is constructed by only the middle lines in the Pythagorean/whole tone musical icosahedra (Fig.~\ref{2_chromatic}). We call this duality Chromatic Pythagorean Duality.
\\
\\
$[{\bf Chromatic} \ {\bf Pythagorean} \ {\bf Duality}]$
\\
\indent
The Pythagorean chain (chromatic scale) is constructed by only the middle lines in the chromatic (Pythagorean)/whole tone musical icosahedra.

\subsection{Common properties between the chromatic/whole tone musical icosahedra and the Pythagorean/whole tone musical icosahedra}
There are some properties that are common to the chromatic/whole tone musical icosahedra and the Pythagorean/whole tone musical icosahedra. For example, tone pairs constructing tritones are set oppositely. Also, some transpositions of the major/minor scales can be done by a symmetry operation of the regular icosahedron. In addition, one can show a duality between major scales and minor scales, and a duality between major/minor scales and the Gregorian modes in the same manner as we dealt with the chromatic/whole tone musical icosahedra. In the next subsection, we show the difference between the chromatic/whole tone musical icosahedra and Pythagorean/whole tone musical icosahedra.

\subsection{Golden gnomon and major/minor triads}
\begin{figure}[H]
\centering
{%
\resizebox*{14cm}{!}{\includegraphics{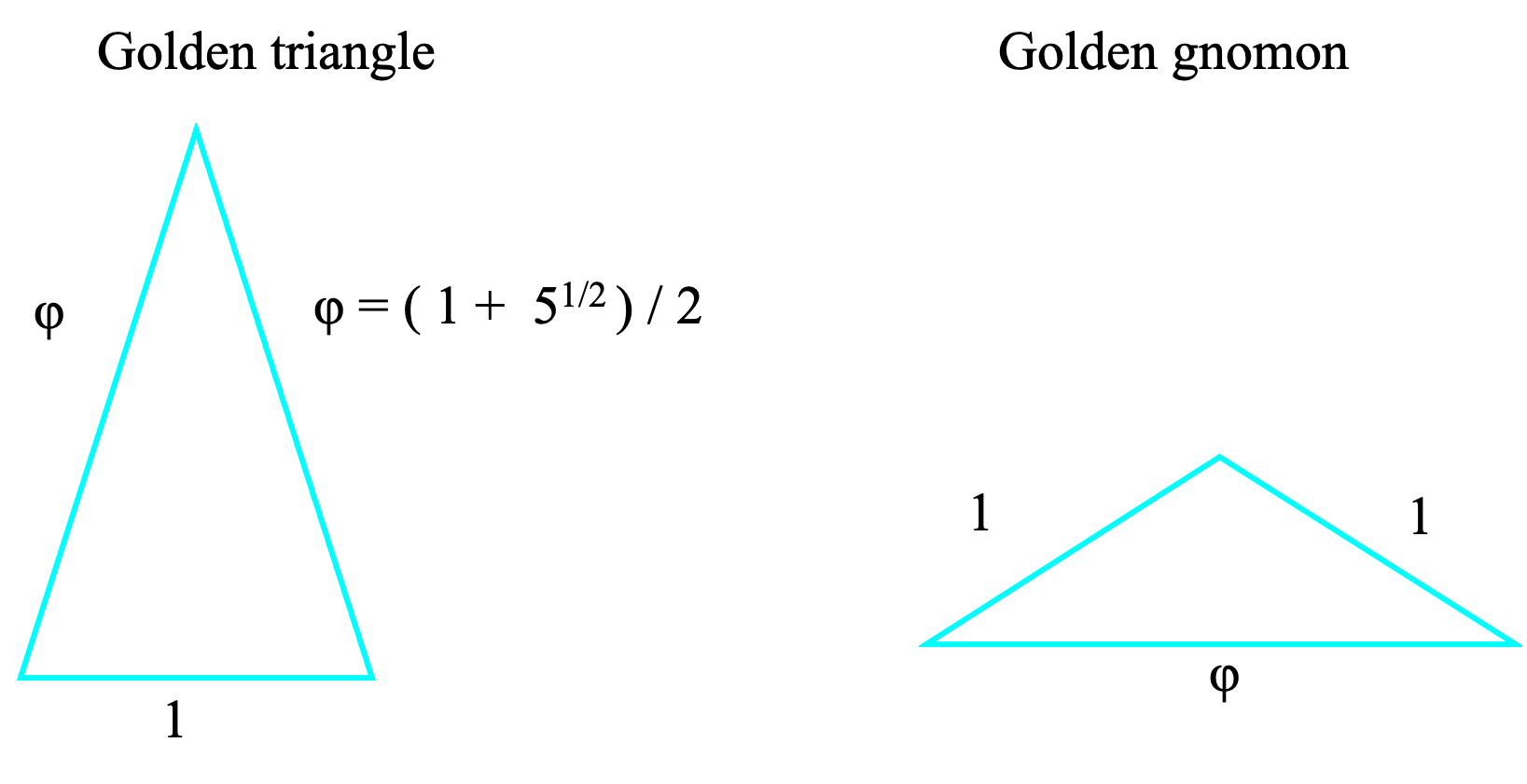}}}\hspace{5pt}
\caption{The golden gnomon and the golden triangle.} \label{2_gnomon}
\end{figure}

\begin{figure}[t]
\centering
{%
\resizebox*{16cm}{!}{\includegraphics{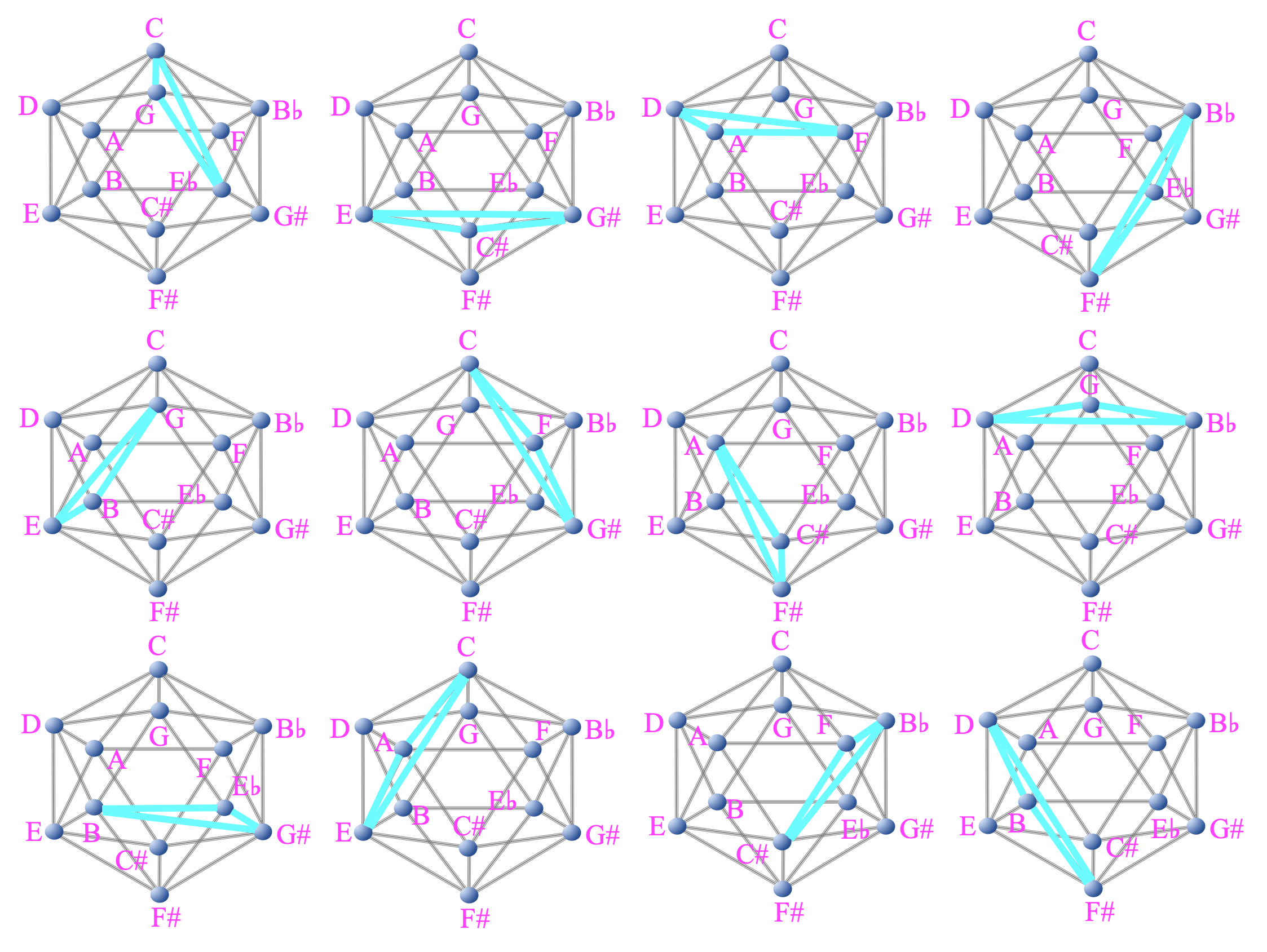}}}\hspace{5pt}
\caption{All the minor triads correspond to the golden gnomons in the type 1' Pythagorean/whole tone musical icosahedron.} \label{2_golden}
\end{figure}

In the subsection, ``Four types of Pythagorean/whole tone musical icosahedron", we showed the Chromatic Pythagorean Duality, and in this subsection, we find another duality: the duality of major/minor triads through the golden triangle and the golden gnomon (Fig.~\ref{2_gnomon}). First, we find the following theorem (Fig.~\ref{2_golden}).
$[{\bf Golden} \ {\bf theorem} \ {\bf 2}]$
\\
\indent
All the minor (major) triads correspond to the golden gnomons in the type 1', type 4' (type2', type 3') Pythagorean/whole tone musical icosahedron.
\\
\\
\indent
The above theorem is the dual theorem of the Golden theorem of the chromatic/whole tone musical icosahedra because the Golden theorem can be obtained by replacing ``major", ``minor" ,``gnomon" and ``Pythagorean" in Golden theorem 2 by `minor", ``major" ,``triangle" and ``chromatic". Therefore, we have the following theorem.
\\
\\
\indent
$[{\bf Golden} \ {\bf Duality}]$

All the major (minor) triads correspond to the golden triangles (gnomons) in the type 1 and type 4 (type 1' and type 4') of the chromatic (Pythagorean)/whole tone musical icosahedra and all the minor (major) triads correspond to golden triangles (gnomons) in the type 2 and type 3 (type 2' and type 3') of the chromatic (Pythagorean)/whole tone musical icosahedra. 
\\
\\
\indent
Also, the fundamental triads for the hexatonic major/minor scales are represented by golden gnomons, and both major/minor triads and hexatonic major/minor scales are uniquely determined in the same manner as we dealt with the chromatic/whole tone musical icosahedra.

\subsection{Generalization of major/minor triads and scales}
\begin{figure}[H]
\centering
{%
\resizebox*{16cm}{!}{\includegraphics{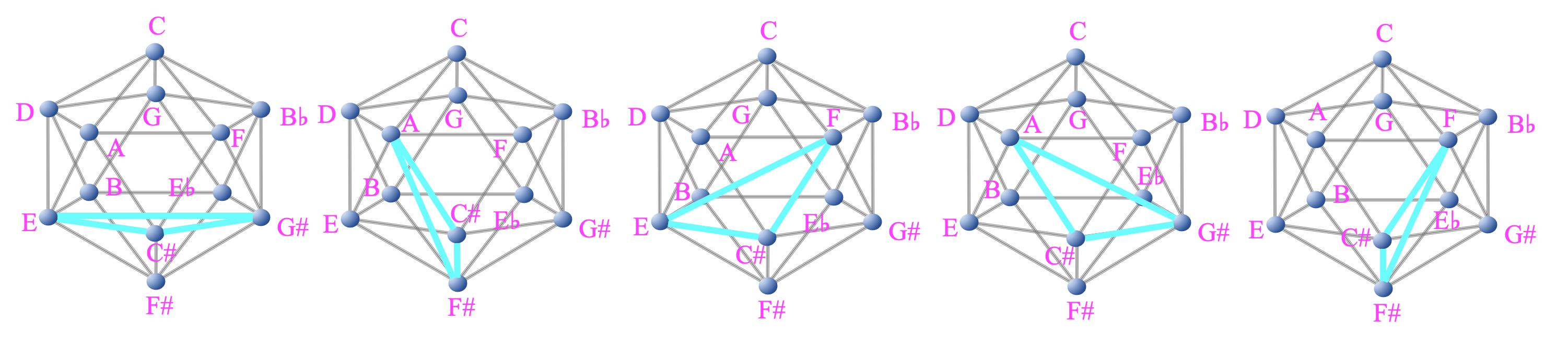}}}\hspace{5pt}
\caption{The first generalization of the $C\sharp$-minor triad.} 
\label{2_general_minor_1}
\end{figure}

\begin{figure}[H]
\centering
{%
\resizebox*{16cm}{!}{\includegraphics{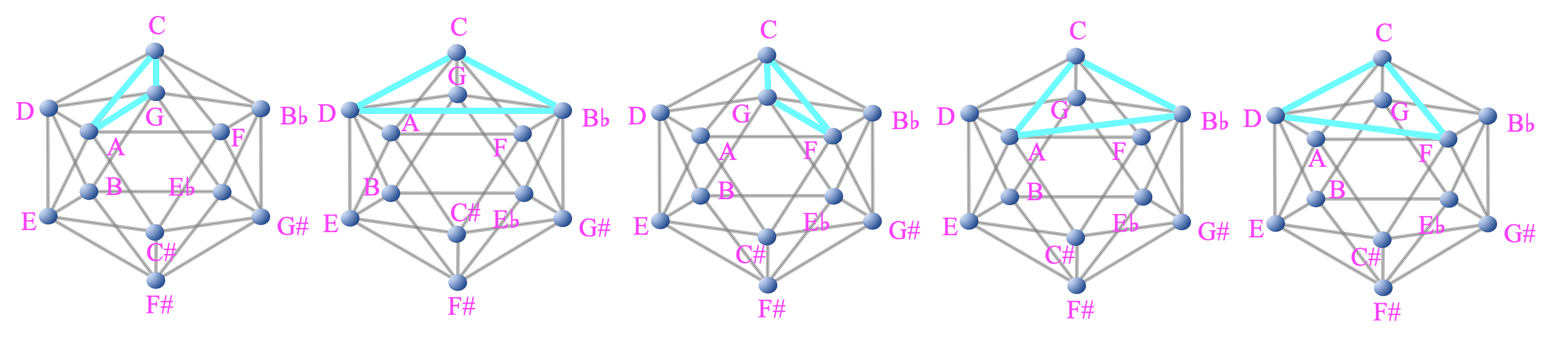}}}\hspace{5pt}
\caption{The second generalization of the $C$-minor triad.} 
\label{2_general_minor_2}
\end{figure}

\begin{figure}[H]
\centering
{%
\resizebox*{16cm}{!}{\includegraphics{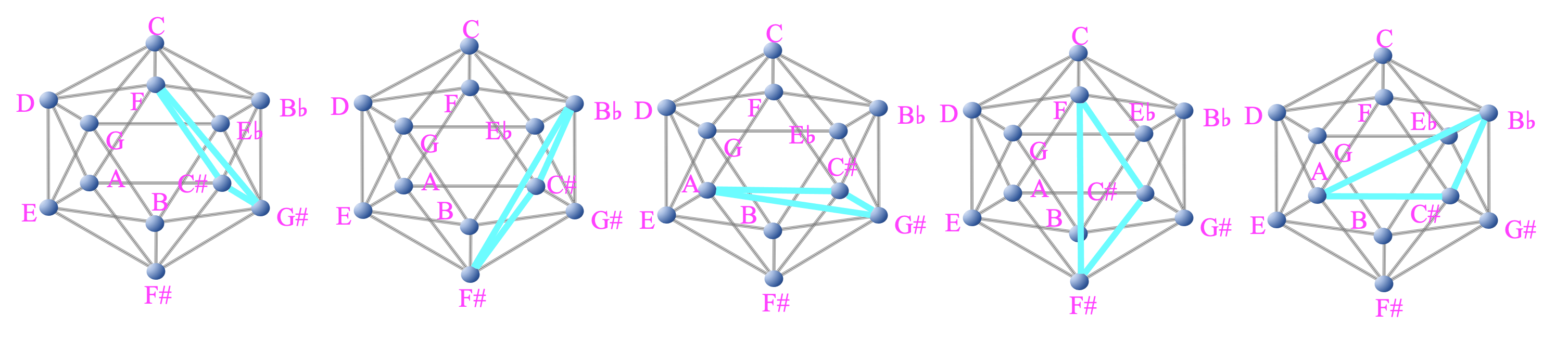}}}\hspace{5pt}
\caption{The first generalization of the $C\sharp$-major triad.} \label{2_general_major_1}
\end{figure}

\begin{figure}[H]
\centering
{%
\resizebox*{16cm}{!}{\includegraphics{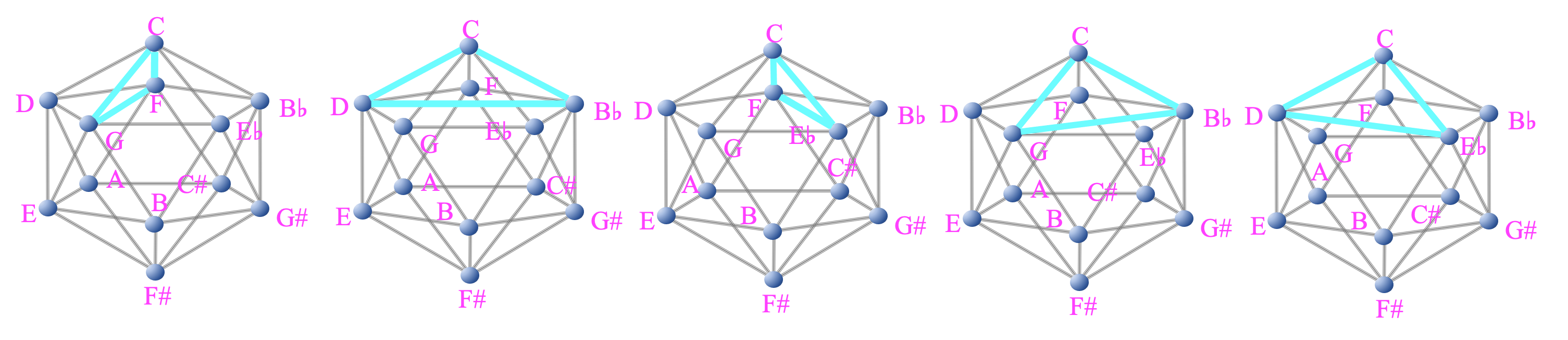}}}\hspace{5pt}
\caption{The second generalization of the  $C$-major triad.} 
\label{2_general_major_2}
\end{figure}

As with the generalizations of the major/minor triads and scales by the chromatic/whole tone musical icosahedra, one has generalizations of the major/minor triads and scales through the icosahedron symmetry and the Pythagorean/whole tone musical icosahedra. The general major/minor triads can be obtained by considering all the golden gnomons in all the types of the Pythagorean/whole tone musical icosahedra (Fig.~\ref{2_general_minor_1}, \ref{2_general_minor_2}, \ref{2_general_major_1}, \ref{2_general_major_2}). We remark that in the type 1', the apex of the triangle made by $CE\flat G$ corresponds to $G$, not $C$, and the apex of the triangle made by $C\sharp E G\sharp$ corresponds to $C\sharp$. Therefore, one can obtain the generalization of the $C\sharp$-minor triad naturally. One may say the type 1' Pythagorean/whole tone musical icosahedron is compatible with the $C\sharp$-minor triad, not with the $C$-minor triad. Note that in the type 1 chromatic/whole tone musical icosahedron, the apex of the triangle made by $CE\flat G$ corresponds to $C$, and then, $CEG$ is naturally generalized. Therefore, the type 1 chromatic/whole tone musical icosahedron is compatible with the $C$-major triad. It is strange that although both the type 1' Pythagorean/whole tone musical icosahedron and the type 1 chromatic/whole tone musical icosahedron satisfy the neighboring condition for the whole tone scale including $C$, the compatibility with the major/minor triads is different. The fact that the apex of the triangle made by $CE\flat G$ corresponds to $G$ in the type 1' Pythagorean/whole tone musical icosahedron may suggest that the main tone of $CEG$ needs not to be $C$ and that $G$ can be the main tone of $CEG$. Because we know that $CEG$ is the $C$-major triad, we tend to consider $C$ to be the main tone of $CEG$. However, $C$ is originally equivalent to $G$ in $CEG$. This point will be discussed further in the next section.

\begin{figure}[H]
\centering
{%
\resizebox*{16cm}{!}{\includegraphics{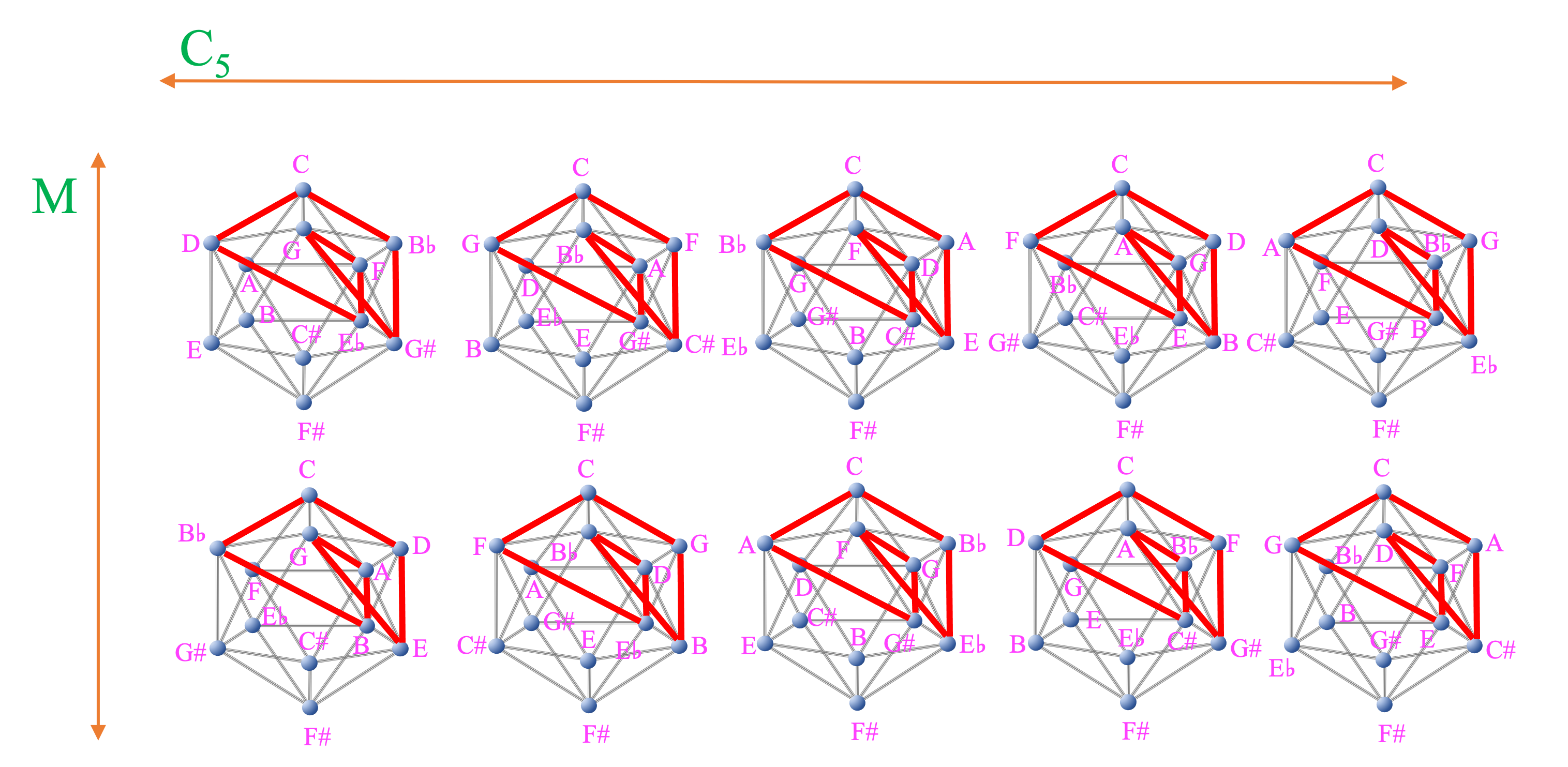}}}\hspace{5pt}
\caption{The first generalization of the $C$-minor scale.} \label{2_general_minor_scale_1}
\end{figure}

\begin{figure}[H]
\centering
{%
\resizebox*{16cm}{!}{\includegraphics{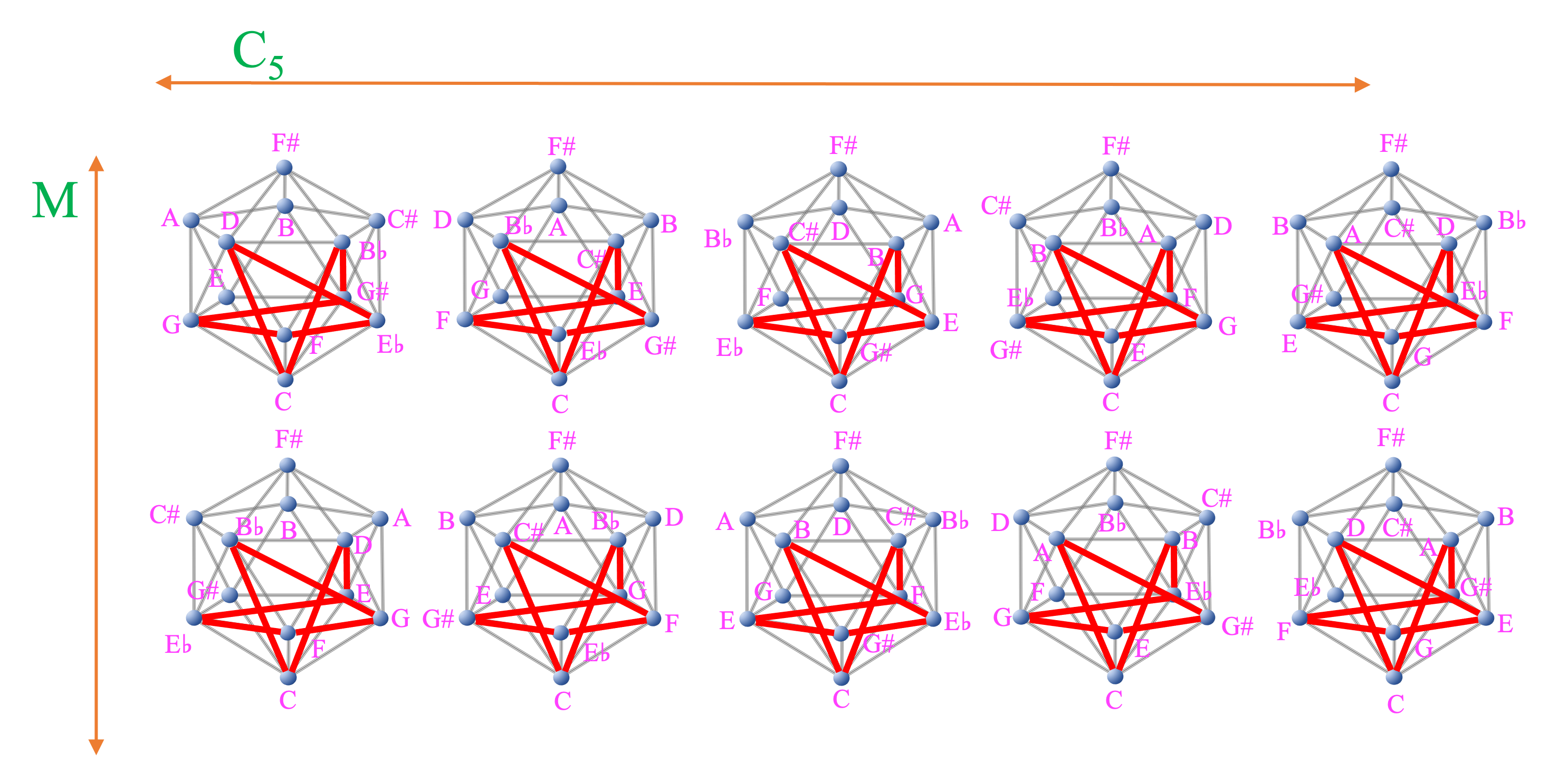}}}\hspace{5pt}
\caption{The second generalization of the $C$-minor scale.} \label{2_general_minor_scale_2}
\end{figure}

\begin{figure}[H]
\centering
{%
\resizebox*{16cm}{!}{\includegraphics{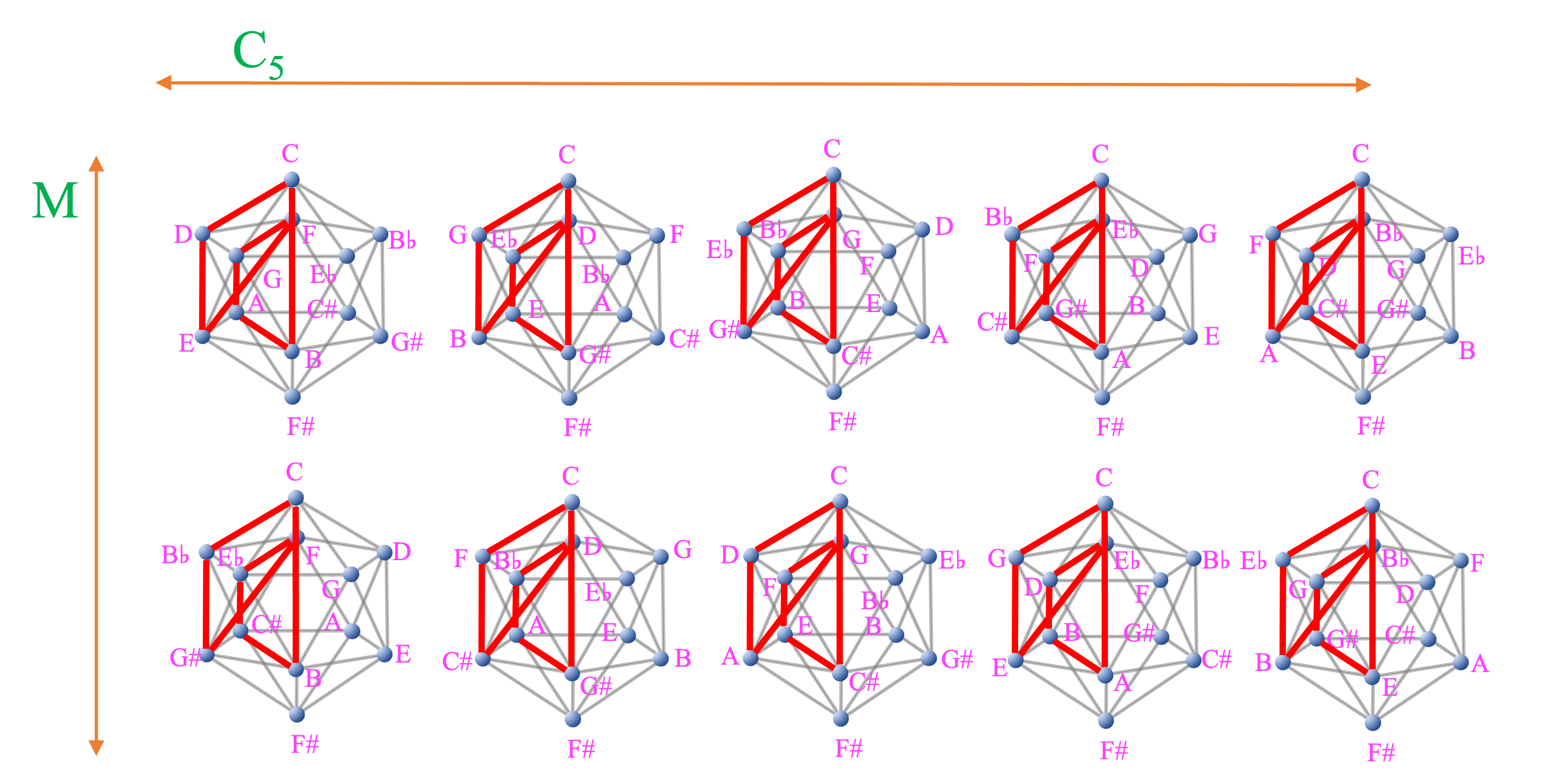}}}\hspace{5pt}
\caption{The first generalization of the $C$-major scale.} \label{2_general_major_scale_1}
\end{figure}

\begin{figure}[H]
\centering
{%
\resizebox*{16cm}{!}{\includegraphics{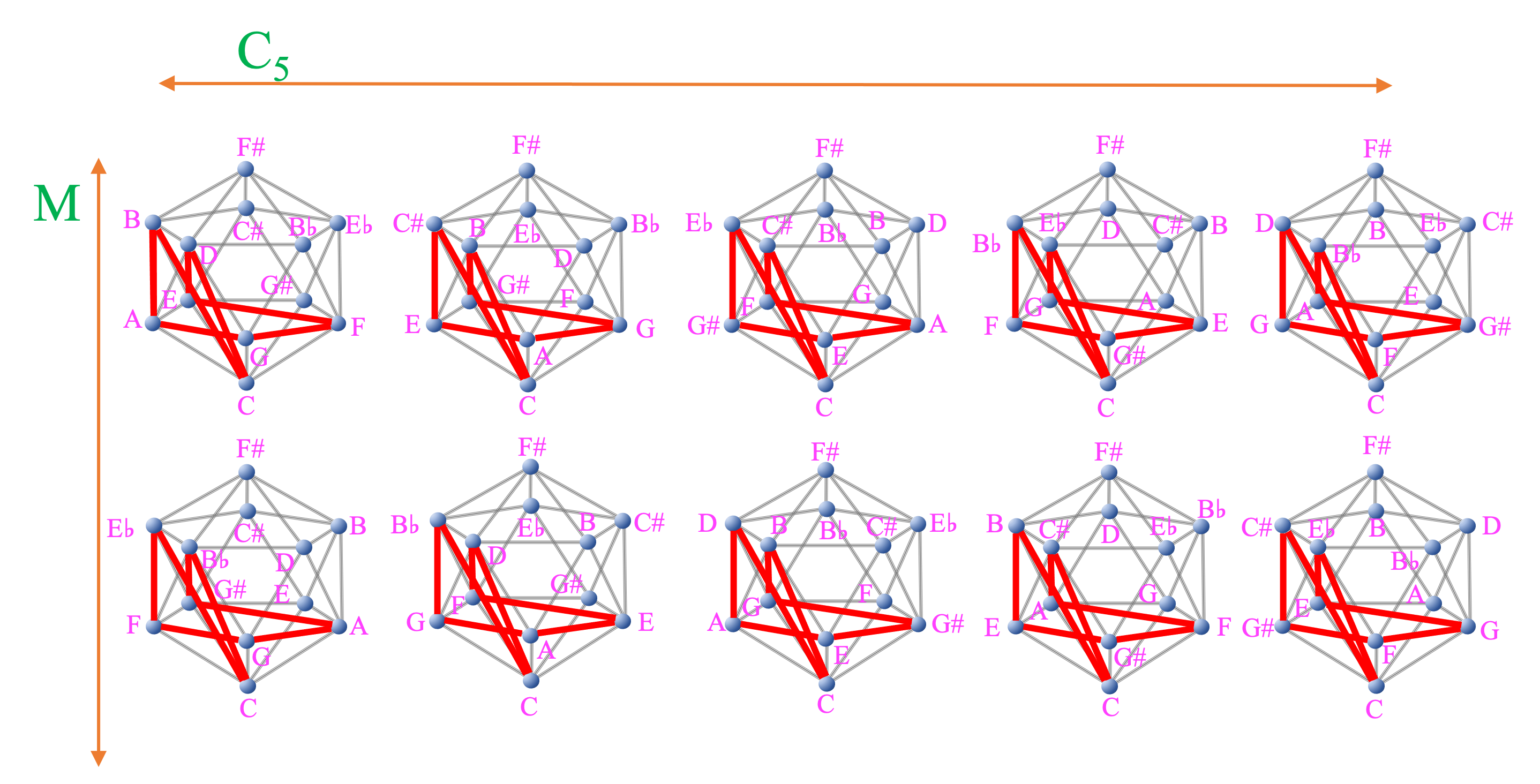}}}\hspace{5pt}
\caption{The second generalization of the $C$-major scale.} \label{2_general_major_scale_2}
\end{figure}

In addition, the general major/minor scales can be obtained by considering all the symmetry transformations of the regular icosahedron to the $C$-major/minor scale in all the types of Pythagorean/whole tone musical icosahedra (Fig.~\ref{2_general_minor_scale_1}, \ref{2_general_minor_scale_2}, \ref{2_general_major_scale_1}, \ref{2_general_major_scale_2}).

\newpage
\section{Golden Major Minor Self-Duality}
In this section, we introduce other kinds of musical icosahedra that have the hexagon-icosahedron symmetry: exceptional musical icosahedra. We show they have the special characteristics with respect to the golden ratio and the major/minor triads (Golden Major Minor Self-Duality).

\subsection{Four Types of Exceptional Musical Icosahedron}
\begin{figure}[H]
\centering
{%
\resizebox*{14cm}{!}{\includegraphics{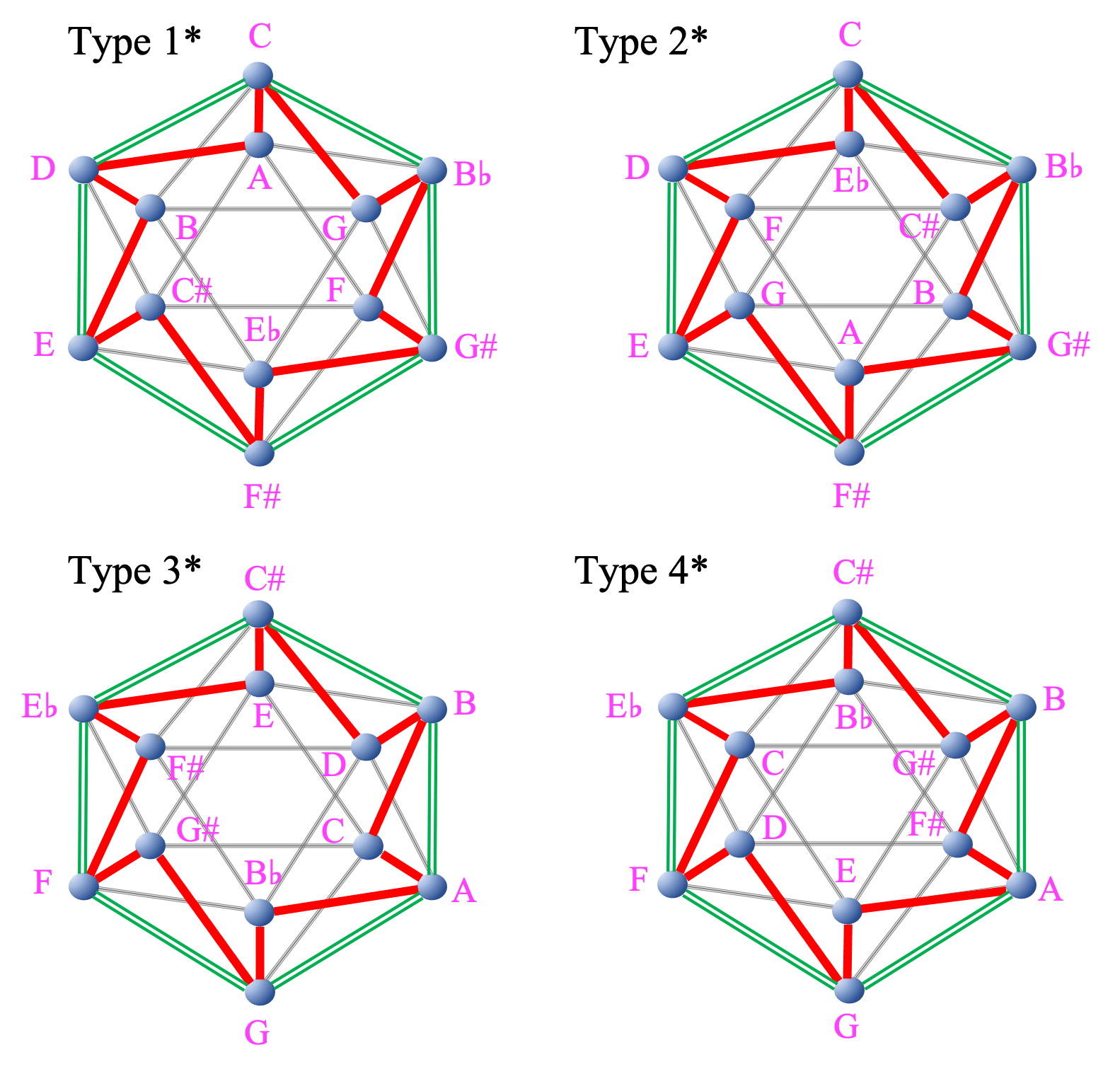}}}\hspace{5pt}
\caption{Type ${\rm 1^*}$, type ${\rm 2^*}$, type ${\rm 3^*}$', and type ${\rm 4^*}$ of the exceptional musical icosahedra.} \label{3_type}
\end{figure}

In order that a musical icosahedron has the hexagon-icosahedron symmetry, one of the whole tone scales is set on the regular hexagon in the projection of the musical icosahedra onto the paper in turn and the other whole tone scale is set on the regular hexagram in the projection of the musical icosahedra onto the paper in turn. Then, one has four types of musical icosahedra other than the chromatic/whole tone musical icosahedra and the Pythagorean/whole tone musical icosahedra with the hexagon-icosahedron symmetry and we call them type ${\rm 1^*}$, type ${\rm 2^*}$, type ${\rm 3^*}$, and type ${\rm 4^*}$ exceptional musical icosahedron (Fig.~\ref{3_type}). A tone just below the $C$ in the type 1, type 2, type 1', type 2', type ${\rm 1^*}$, type ${\rm 2^*}$ is $C\sharp$, $B$, $G$, $F$, $A$, $E\flat$, and they make the whole tone scale including $C\sharp$.

We remark that one has the type ${\rm 4^*}$ (type ${\rm 3^*}$) by raising the all the tones of the type ${\rm 1^*}$ (type ${\rm 2^*}$)  by one semitone. This property is common to the choromatic/whole tone musical icosahedra and the Pythagorean/whole tone musical icosahedra. In the type ${\rm 1^*}$ and type ${\rm 2^*}$, the whole tone scale including $C$ is set on the regular hexagon in the projection of the musical icosahedra onto the paper while in the type ${\rm 3^*}$ and the type ${\rm 4^*}$, the whole tone scale including $C\sharp$ is set on the regular hexagram in the projection of the musical icosahedra onto the paper.

The red lines on the figures do not represent the chromatic scale or the Pythagorean chain. The red lines on the type ${\rm 1^*}$ represent $C$, $A$, $D$, $B$, $E$, $C\sharp$, $F\sharp$, $E\flat$, $G\sharp$, $F$ $B\flat$, $G$. This scale is symmetric in ``modulo 4" because $A - C =9=4\times 2+1$, $D - A=4+1$. Also, the red lines on the type ${\rm 2^*}$ represent $C$, $E\flat$, $D$, $F$, $E$, $G$, $F\sharp$, $A$, $G\sharp$, $B$, $B\flat$, $C\sharp$. This scale is symmetric in ``modulo 4" because $E\flat$ - $C$ = 3, $D - E\flat=11=4\times 2+3$. These scales can be obtained by operating the mirror reflection to the red lines representing the chromatic scale in the type 1 of the chromatic/whole tone musical icosahedra and the type 1' of the Pythagorean/whole tone musical icosahedra. Also, the red lines on the type ${\rm 3^*}$ represent $C$, $B$, $D$, $C\sharp$, $E$, $E\flat$, $F\sharp$, $F$, $G\sharp$, $G$, $B\flat$, $A$. This scale is symmetric in ``modulo 4" because $B - C=11=4\times 2+3$, $D-B=3$. Also, the red lines on the type ${\rm 4^*}$ represent $C$, $F$, $D$, $G$, $E$, $A$, $F\sharp$, $B$, $G\sharp$, $C\sharp$ $B\flat$, $E\flat$. This scale is symmetric in ``modulo 4" because $F-C=5=4+1$, $D-F=9=4\times 2+1$. These scales can be obtained by operating the mirror reflection to the red lines representing the chromatic scale in the type 1 of the chromatic/whole tone musical icosahedra and the type 1' of the Pythagorean/whole tone musical icosahedra. Note that the above four scales are the only scales that are not symmetric in ``modeulo 12" and are symmetric in ``modulo 4". The properties common to the chromatic/whole tone musical icosahedra and the Pythagorean/whole tone musical icosahedra is also common to the exceptional musical icosahedra other than the major-minor duality. Instead of the major-minor duality, the exceptional musical icosahedra have another beautiful property.

\subsection{Golden Major Minor Self-Duality}
\begin{figure}[H]
\centering
{%
\resizebox*{14cm}{!}{\includegraphics{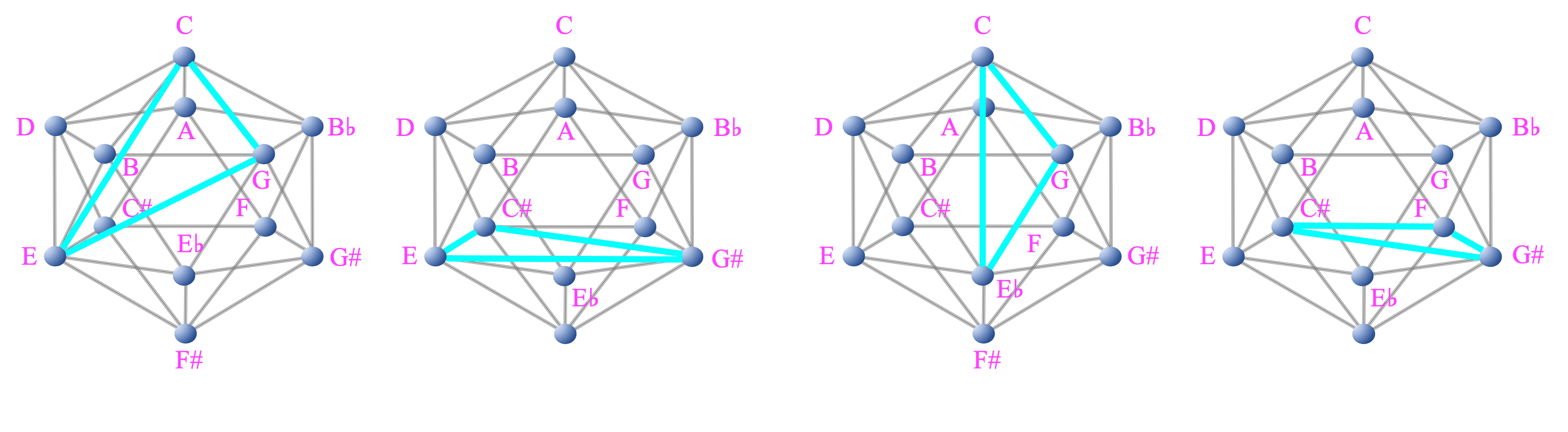}}}\hspace{5pt}
\caption{The $C$-major triad and $C\sharp$-minor triad (The $C$-minor triad and $C\sharp$-major triad) correspond to the golden triangle (golden gnomon) on type ${\rm 1^*}$. } \label{3_harmony1}
\end{figure}

\begin{figure}[H]
\centering
{%
\resizebox*{14cm}{!}{\includegraphics{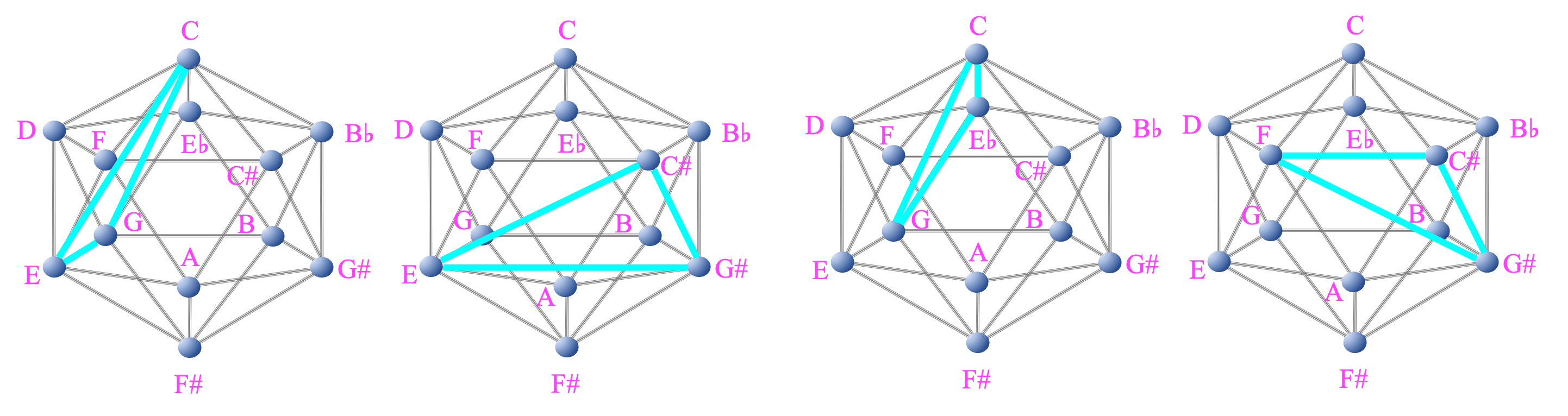}}}\hspace{5pt}
\caption{The $C$-major triad and $C\sharp$-minor triad (The $C$-minor triad and $C\sharp$-major triad) correspond to the golden triangle (golden gnomon) on type ${\rm 2^*}$.} \label{3_harmony2}
\end{figure}

In each type of the exceptional musical icosahedra, not all the major/minor triads correspond to the golden triangles or not all the major/minor triads correspond to the golden gnomons. Instead, one obtains the Golden Major Minor Self-Duality summarized below and in Fig.~\ref{3_harmony1}, \ref{3_harmony2}.
\\
\\
$[{\bf Golden} \ {\bf Major}  \ {\bf Minor}   \ {\bf Self}-{\bf Duality}]$
The major triad whose base tone is $C$, $D$, $E$, $F\sharp$, $G\sharp$, or $B\flat$ ($C\sharp$, $E\flat$, $F$, $G$, $A$, or $B$) corresponds to the golden triangle and the minor triad whose base tone is $C\sharp$, $E\flat$, $F$, $G$, $A$, or $B$ ($C$, $D$, $E$, $F\sharp$, $G\sharp$, or $B\flat$) corresponds to the golden gnomon on the type ${\rm 1^*}$ and the type ${\rm 4^*}$ (the type ${\rm 2^*}$ and the type ${\rm 3^*}$).
\\
\\
\indent
Because all the major/minor triads correspond to the golden triangles or the golden gnomons on each type of the exceptional musical icosahedra, then, the exceptional musical icosahedra may be applied to harmonic analysis because a combination of some of the major triads and the minor triads correspond to a combination of some of the golden triangles and the golden gnomons.

\subsection{Generalization of major/minor triads and scales}
As with the generalizations of the major/minor triads and scales by the chromatic/whole tone musical icosahedra and by the Pythagorean/whole tone musical icosahedra, one has generalizations of the major/minor triads and scales by the exceptional musical icosahedra. General major/minor triads can be obtained by considering all the golden triangles and gnomons in all the types of the exceptional musical icosahedra. In the previous section, we discussed the compatibility of musical icosahedra with the major/minor triads. This problem appears again in this situation. In the type ${\rm 1^*}$, the apex of the triangle made by $CEG$, $C\sharp EG\sharp$, $CE\flat G$, $C\sharp FG\sharp$ is not $C$, $C\sharp$, $C$, $C\sharp$ but $E$, $G\sharp$, $G$, $F$ (Fig.~\ref{3_harmony2}). In the type ${\rm 2^*}$, the apex of the triangle made by $C\sharp EG\sharp$, $CE\flat G$ is not $C\sharp$, $C$ but $E$, $E\flat$ (Fig.~\ref{3_harmony2}). Then, in the exceptional musical icosahedra, the generalized major/minor triad $X$ should not be regarded as the generalization of the major/minor triad $Y$.

\begin{figure}[H]
\centering
{%
\resizebox*{14cm}{!}{\includegraphics{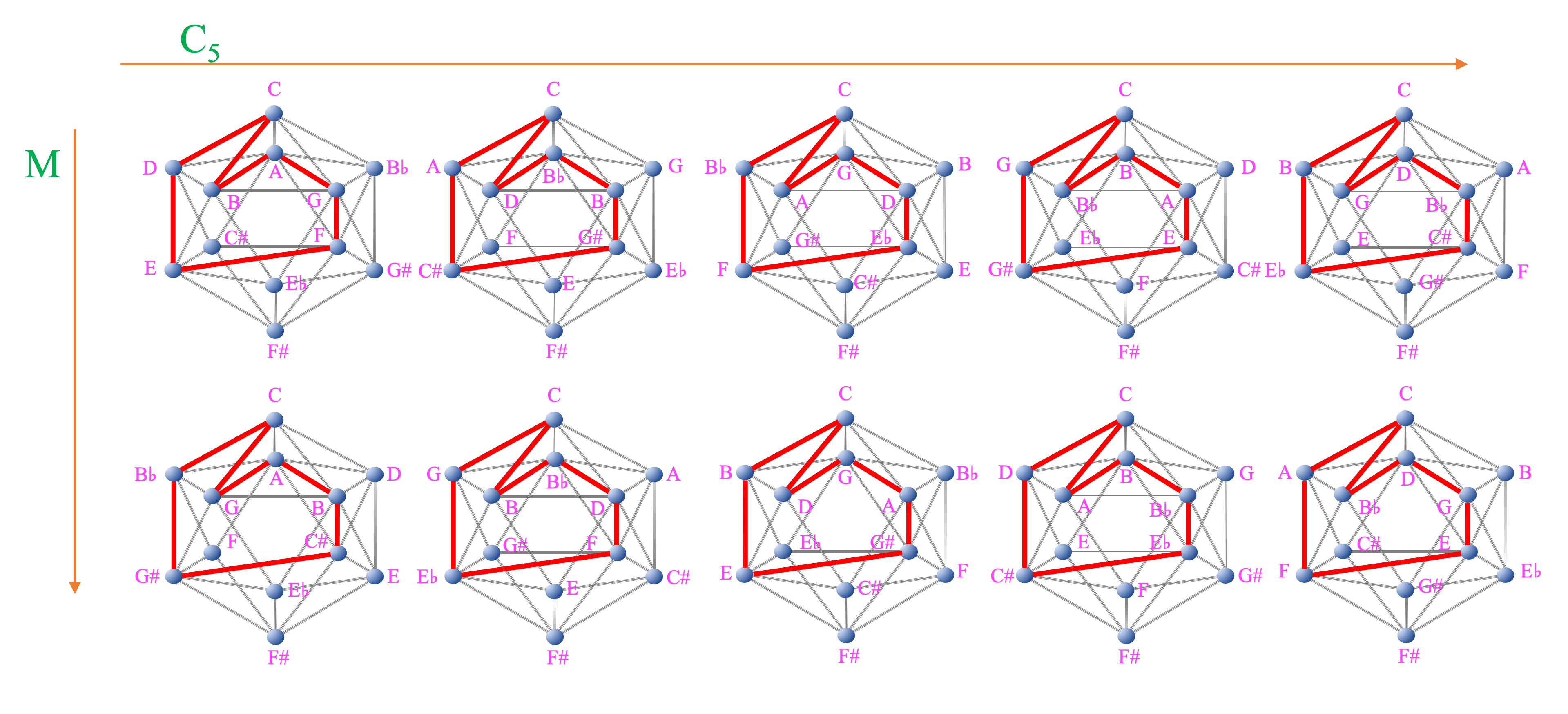}}}\hspace{5pt}
\caption{The first generalization of the $C$ major scale by the type ${\rm 1^*}$ exceptional musical icosahedron.} \label{3_generalization_major1_1}
\end{figure}

\begin{figure}[H]
\centering
{%
\resizebox*{14cm}{!}{\includegraphics{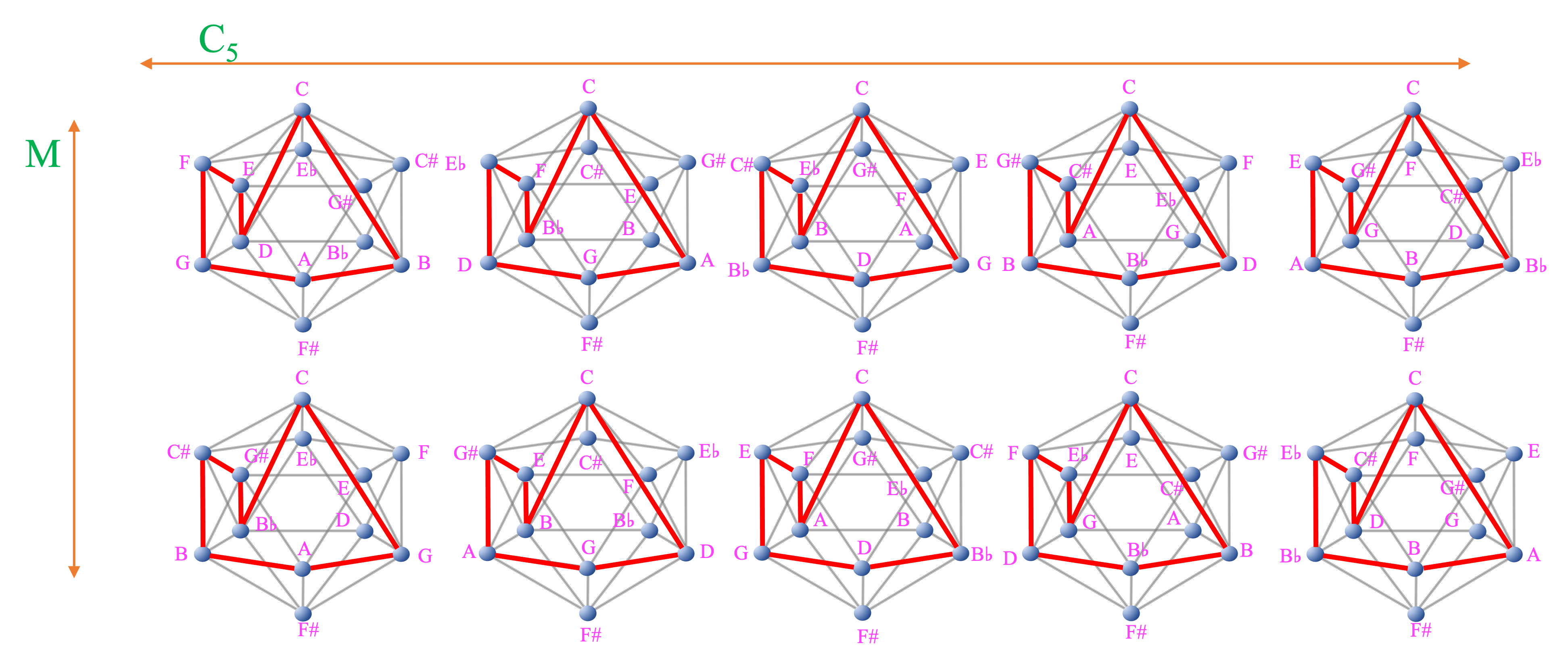}}}\hspace{5pt}
\caption{The second generalization of the $C$ major scale by the type ${\rm 1^*}$ exceptional musical icosahedron.} \label{3_generalization_major1_2}
\end{figure}

\begin{figure}[H]
\centering
{%
\resizebox*{14cm}{!}{\includegraphics{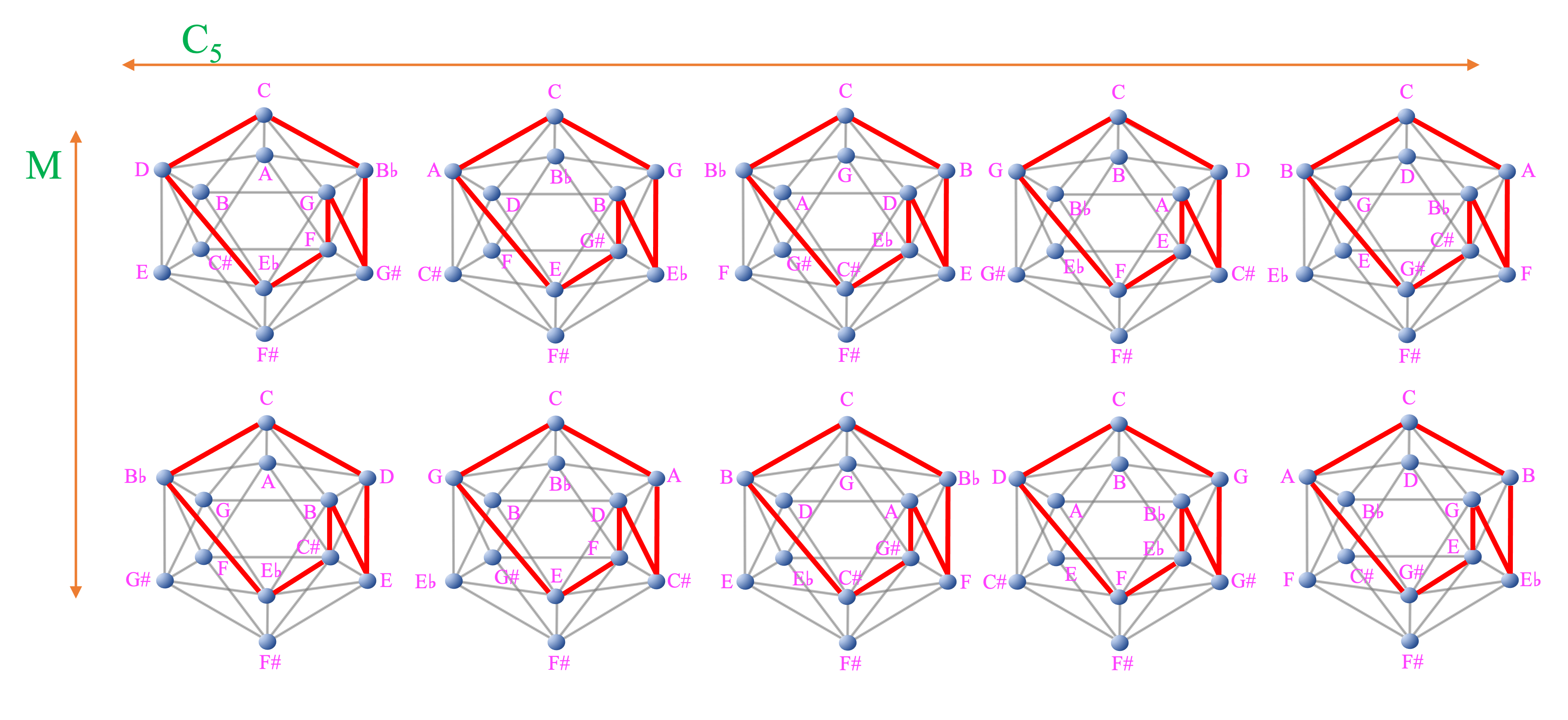}}}\hspace{5pt}
\caption{The first generalization of the $C$ minor scale by the type ${\rm 1^*}$ exceptional musical icosahedron.} \label{3_generalization_minor1_1}
\end{figure}

\begin{figure}[H]
\centering
{%
\resizebox*{14cm}{!}{\includegraphics{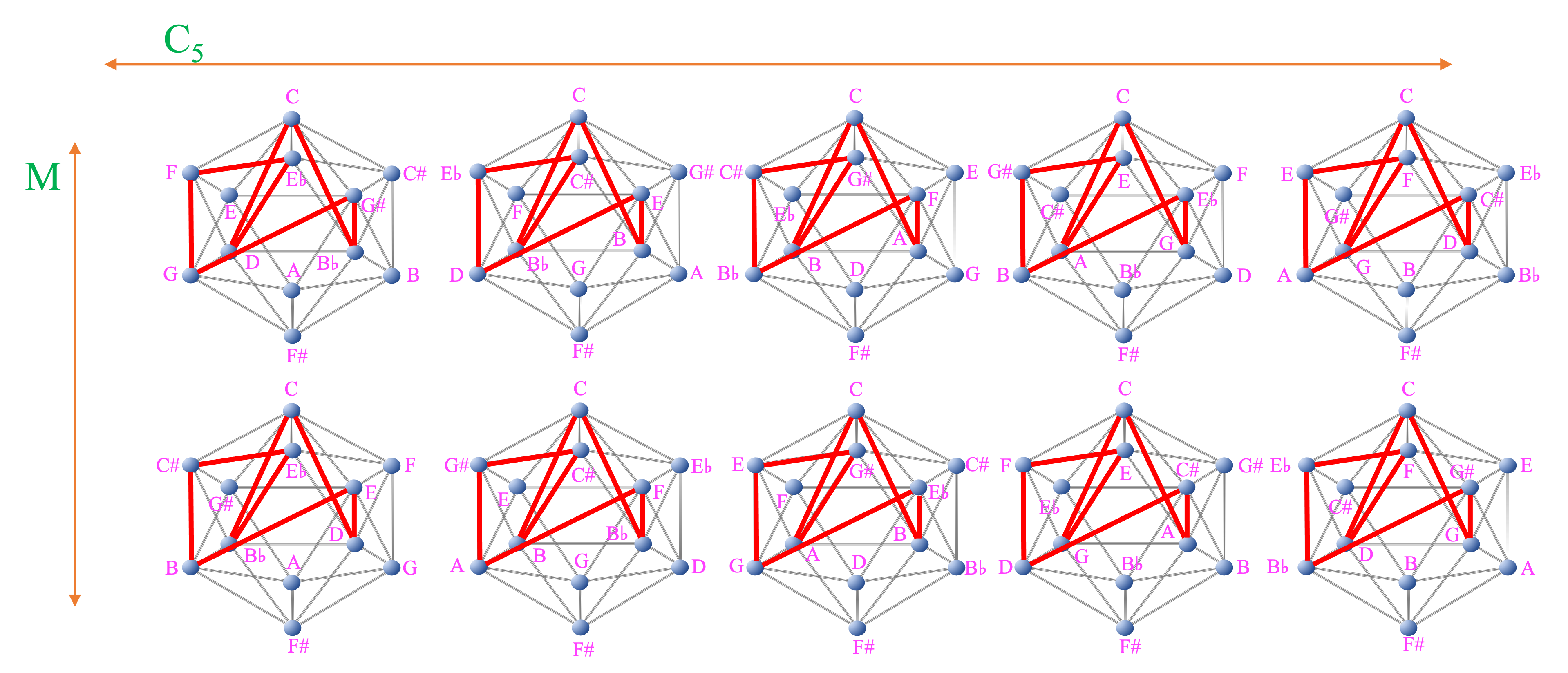}}}\hspace{5pt}
\caption{The second generalization of the $C$ minor scale by the type ${\rm 1^*}$ exceptional musical icosahedron.} \label{3_generalization_minor1_2}
\end{figure}

\begin{figure}[H]
\centering
{%
\resizebox*{14cm}{!}{\includegraphics{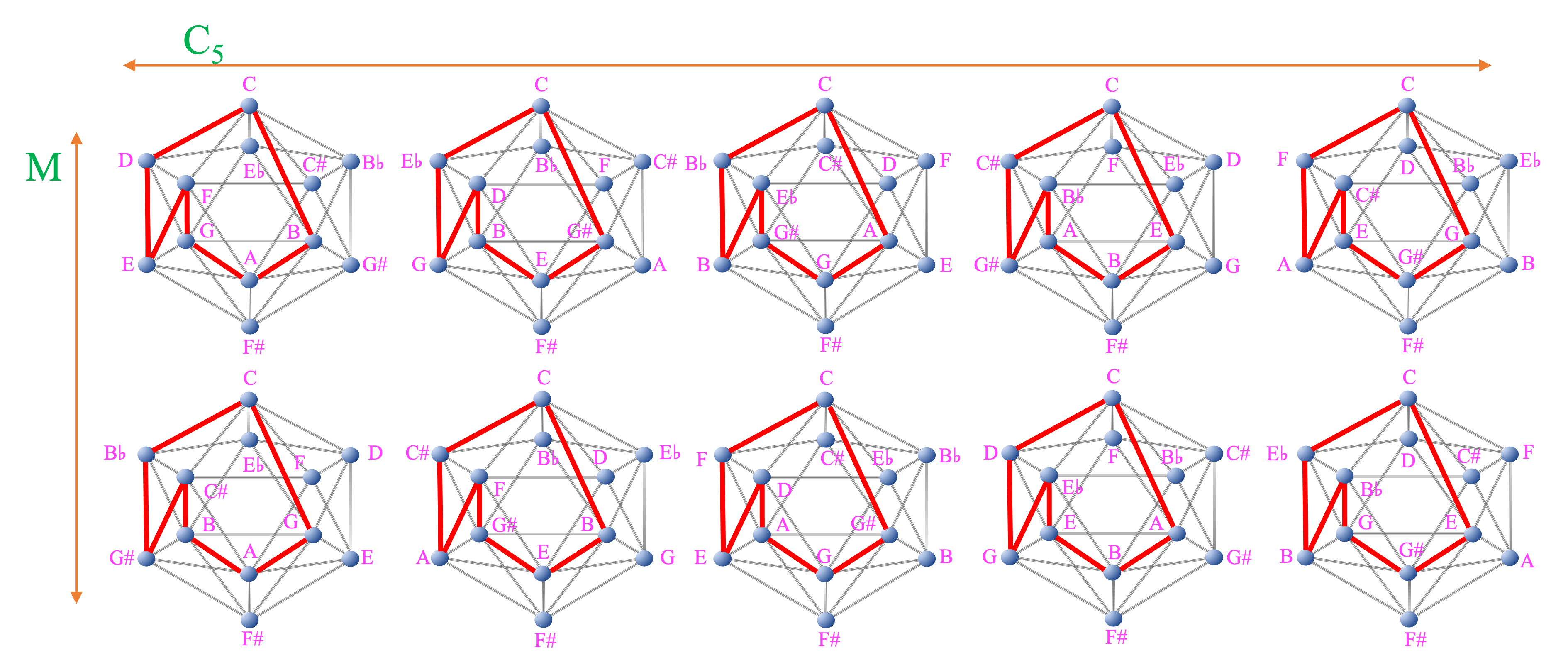}}}\hspace{5pt}
\caption{The first generalization of the $C$ major scale by the type ${\rm 2^*}$ exceptional musical icosahedron.} \label{3_generalization_major2_1}
\end{figure}

\begin{figure}[H]
\centering
{%
\resizebox*{14cm}{!}{\includegraphics{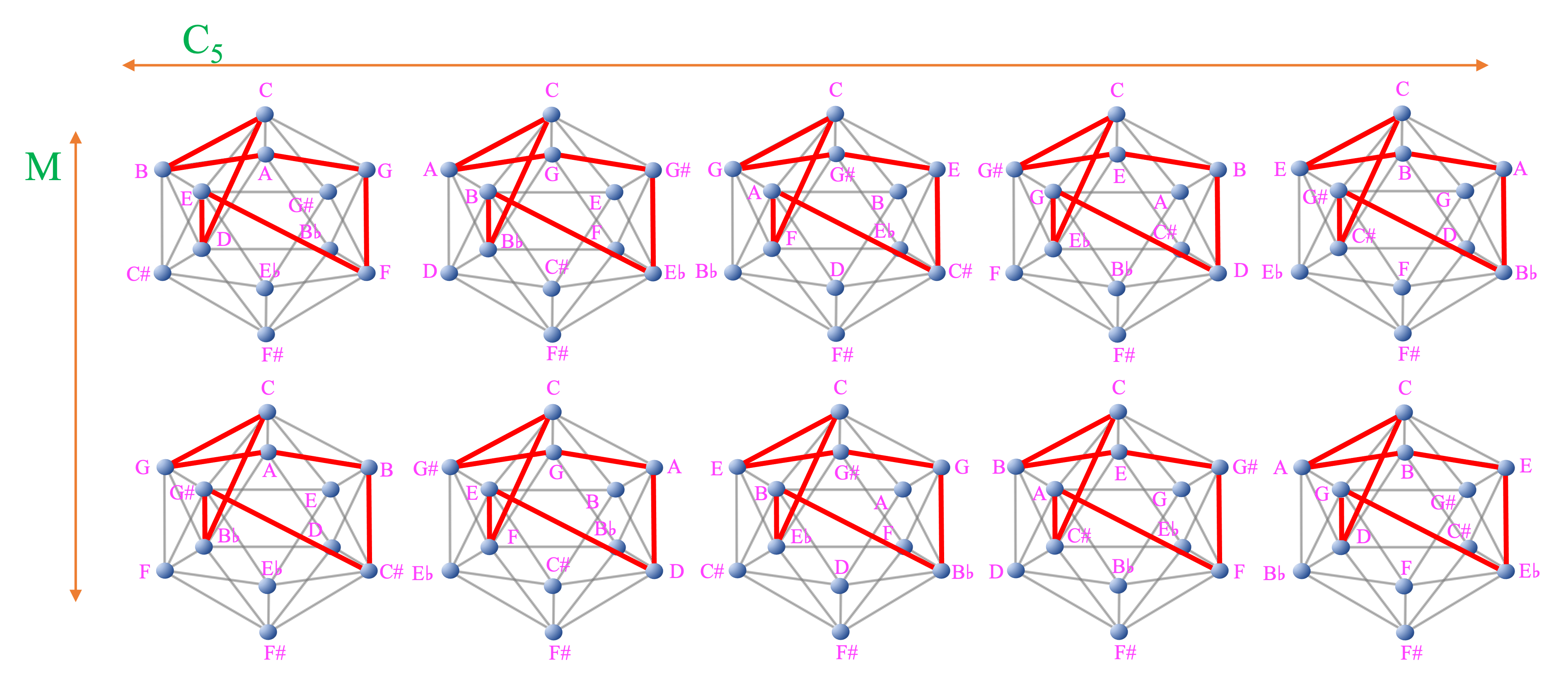}}}\hspace{5pt}
\caption{The second generalization of the $C$ major scale by the type ${\rm 2^*}$ exceptional musical icosahedron.} \label{3_generalization_major2_2}
\end{figure}

\begin{figure}[H]
\centering
{%
\resizebox*{14cm}{!}{\includegraphics{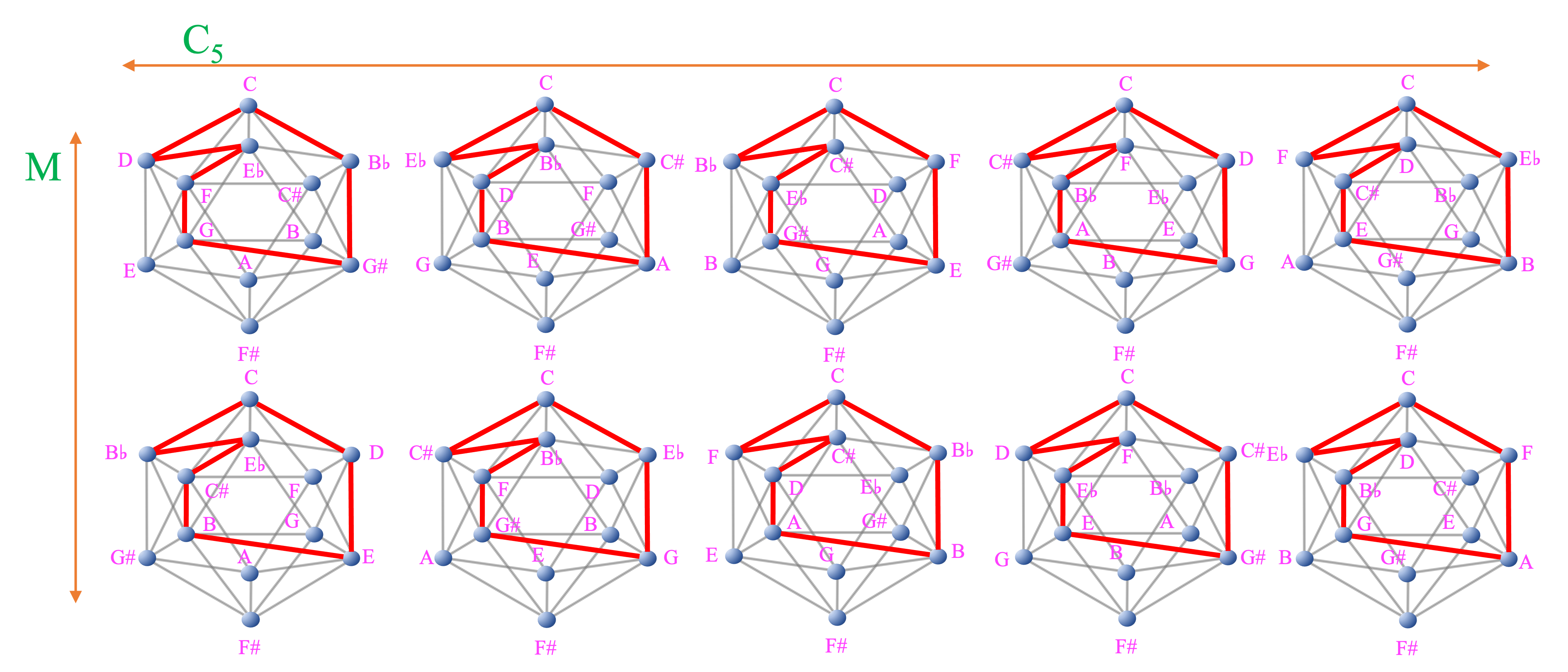}}}\hspace{5pt}
\caption{The first generalization of the $C$ minor scale by the type ${\rm 2^*}$ exceptional musical icosahedron.} \label{3_generalization_minor2_1}
\end{figure}

\begin{figure}[H]
\centering
{%
\resizebox*{14cm}{!}{\includegraphics{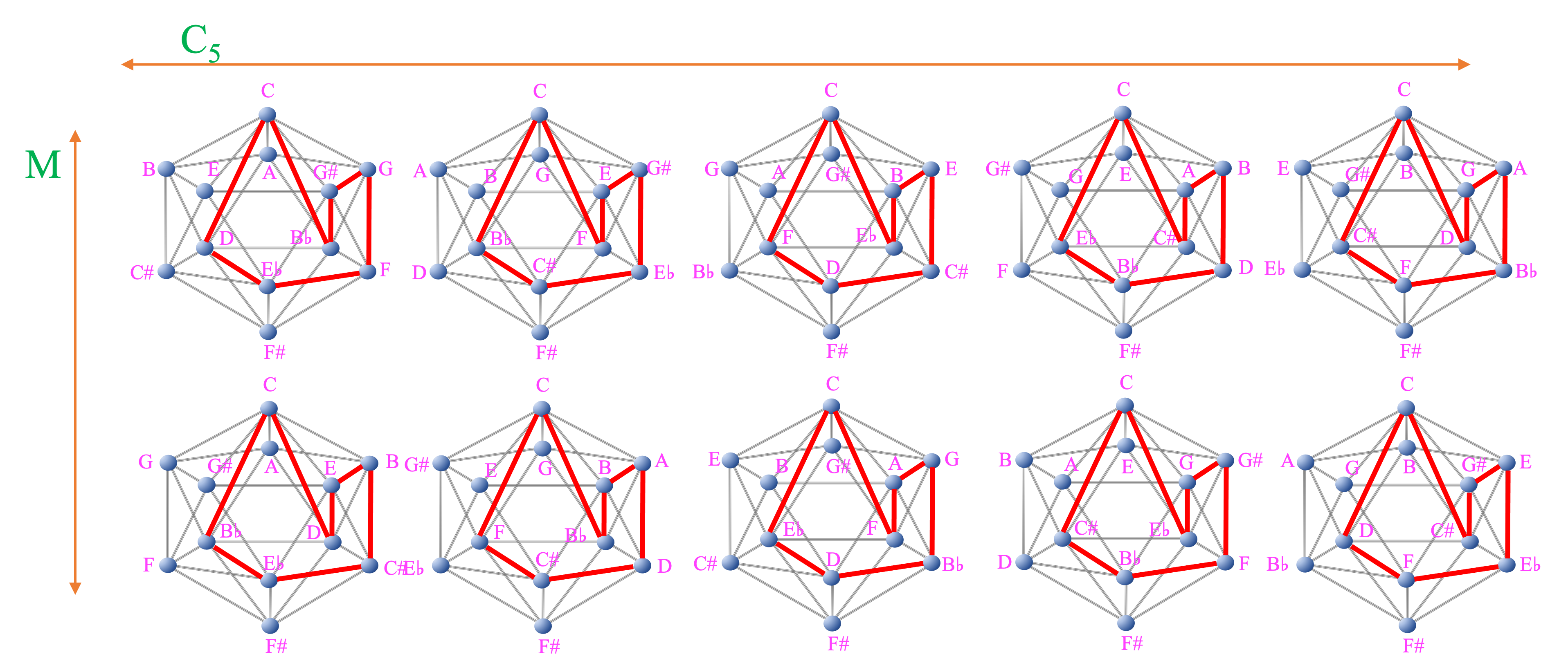}}}\hspace{5pt}
\caption{The second generalization of the $C$ minor scale by the type ${\rm 2^*}$ exceptional musical icosahedron.} \label{3_generalization_minor2_2}
\end{figure}

Anyway, one can obtain the generalizations of the major/minor scales by considering all the symmetry transformation of the regular icosahedron to the $C$-major/minor scale in all the types of the exceptional musical icosahedra (Fig.~\ref{3_generalization_major1_1}, \ref{3_generalization_major1_2}, \ref{3_generalization_minor1_1}, \ref{3_generalization_minor1_2}, \ref{3_generalization_major2_1}, \ref{3_generalization_major2_2}, \ref{3_generalization_minor2_1}, \ref{3_generalization_minor2_2}). Then, one can determine the tonic of the generalized major/minor scales by drawing lines corresponding to the $C$-major/minor triad on the left-top musical icosahedron in each figure and drawing the same lines on the other musical icosahedra.


\newpage
\section{Conclusion}

Since antiquity, Western music has always tried to base its foundation on mathematics. Until the end of the Renaissance period, music theory always involved a discussion of interval theory and acoustic theory of sound from the perspective of proportion (e.g. the theoretical writings by Boethius, Guido of Alezzo, Franco of Cologne, and Vincenzo Galilei). However, after the beginning of the Baroque period, especially after the 18th century, human perception took priority over an arithmetic truth in music; it made the connection between mathematics and music weaker\footnote{For example, it is typically indicated in a title of piece-collections for a practice of the piano published in the middle of the 18th century (Hartung 1749); ``The demonstrative music theory. It is constructed on the true principles, freed from arithmetical difficulties, but alternated by (an explanation on) harmony… [Die \emph{demonstrativische} theoria musica. Auf ihre wahre \emph{Principia} gebauet, von vielen \emph{arithmetischen Subtilitaeten} befreyet, hingegen die Abwechselung deren \emph{Harmonien}...]"}. Atonal musical experiments, which were especially prevalent in the 20th century and can be represented by the twelve-tone technique, again gave mathematics a place in music theory. In the sense of the relationship between mathematics and music, it seems likely that it goes back to the ancient or medieval age in terms of its spirit. However, there is a big difference between them: for example, Lorenz Christoph Mitlzer (1711-1778), who strongly asserted the importance of mathematics inmusic even in the middle of the 18th century and thus might be considered to be an enlightened musician witha mind of ``renaissance", believed that ``the bases of harmony and pure composition are invariable in music [die Gr\"unde der Harmonie und der reinen Composition unver\"anderlich in der Musik]" (1742, Vorrede): these ``Gr\"unde" were for Mitzler no other than the bases that can be derived from mathematics. For Mitzler, or for the musicians who possessed mathematical-based world of views until the 18th century, music based on mathematical truth is always invariable and beautiful. On the other hand, it can be interpreted that atonal music and music theories, which are derived from it or are illustrated for contributing to it in the last century, generally used mathematics as a measure for those objective but artificial systematizations; they do use mathematics, but they do not start from the mathematical truth, which – according to the ancient and medieval spirit allows music to exist in pure and beautiful form. In terms of the mathematical starting point, our idea should therefore be understood as a music theory that was born from the traditional view of music on mathematics.

Now, the results presented in this paper should briefly be summarized.

First, we introduced the concept of a musical icosahedron: the regular icosahedron each of whose vertices has one of the 12 tones without duplication, and the neighboring/distant condition that topology of scales is reflected in the icosahedron. We found that there are only four musical icosahedra satisfying the neighboring condition for the chromatic scale and the neighboring condition for one of the whole tone scales and the distant condition for the other whole tone scale. We named them chromatic/whole tone musical icosahedra. We remark that a transformation linking the type 1 (the type 2) and the type 4 (the type 3) gives rise to the analogy of the spin in quantum physics or representation theory. We also remark that the chromatic/whole tone musical icosahedra have the hexagon-icosahedron symmetry (an operation of raising all the tones of a given scale by two semitones corresponds to a symmetry transformation of the regular icosahedron). The chromatic/whole tone musical icosahedra lead to the dualities between the major triads and the minor triads, and the major scales and the minor scales, the fundamental triads for the hexatonic major/minor scales, and the duality between the major/minor scales and the Gregorian modes.

Then, we focused on other musical icosahedra that have the hexagon-icosahedron symmetry.

We found that there are only four musical icosahedra satisfying the neighboring condition for the Pythagorean chain (succession of the fifth based on $C$) and  the neighboring condition for one of the whole tone scales and the distant condition for the other whole tone scale. We named them chromatic/whole tone musical icosahedra. Then, we derived dualities between the chromatic scale and the Pythagorean chain by the chromatic/whole tone musical icosahedra and the Pythagorean/whole tone musical icosahedra. First, we showed the Pythagorean chain is constructed by only the middle lines in the chromatic/whole tone musical icosahedra and the chromatic scale is constructed by only the middle lines in the Pythagorean/whole tone musical icosahedra. Next, we found some similarities between the chromatic/whole tone musical icosahedra and the Pythagorean/whole tone musical icosahedra: the hexagon-icosahedron symmetry, tritone and opposite vertex, transposition and symmetry. Also, we showed the golden gnomons determine uniquely the major/minor triads and the fundamental triads for the hexatonic major/minor scales in the Pythagorean/whole tone musical icosahedron. This should be regarded as the dual theorem of the Golden theorem of the chromatic/whole tone musical icosahedra. In addition, we generalized the major/minor triads and scales through the Pythagorean/whole tone musical icosahedra and the symmetry of the regular icosahedron. We found that the type 1' Pythagorean/whole tone musical icosahedron is compatible with the $C\sharp$-minor triad while the chromatic/whole tone musical icosahedron is compatible with the $C$-major triad although both of them satisfy the neighboring condition for the whole tone scale including $C$. The dualities between the chromatic scale and the Pythagorean chain that are found in this paper by using musical icosahedra cannot be derived by the simple algebraic relation among 12 tones ($D-C=2$ etc.), and can be derived by combining it and mathematical properties of the regular icosahedron.

In addition, we dealt with four other musical icosahedra that have the hexagon-icosahedron symmetry, exceptional musical icosahedra. Major triads whose base tone is $C$, $D$, $E$, $F\sharp$, $G\sharp$, or $B\flat$ ($C\sharp$, $E\flat$, $F$, $G$, $A$, or $B$) and minor triads whose base tone is $C\sharp$, $E\flat$, $F$, $G$, $A$, or $B$ ($C$, $D$, $E$, $F\sharp$, $G\sharp$, or $B\flat$) are set on the golden triangles/gnomons (the golden gnomons/triangles) on them. These exceptional musical icosahedra may be applied to harmonic analysis because one can deal with a combination of some of the major/minor triads by some of the golden triangles and the golden gnomons. We also showed generalizations of the major/minor triads and scales through the exceptional musical icosahedra and the symmetry of the regular icosahedron.

Using these generalized scales may be a new way to compose music that has an atmosphere the major/minor scales do not possess. Composing this kind of music is one of the major trends in classical music from the Romantic age to today, accomplished variously through the use of the Gregorian modes, whole tone scales, twelve-tone music, tone cluster, graphic notation etc. The generalized scales proposed in this paper are not meant to deny the major/minor scales or create chaotic music. Rather, they extend the major/minor scales systematically through symmetries of the regular icosahedra. Because there are exceedingly many possible scales, this extension of the major/minor scales may be a useful tool in composing music, creating novel musical landscapes. It is also interesting to visualize music by using these musical icosahedra because beautiful harmonies (the major/minor triads) are represented by a beautiful figure, the golden triangles/gnomons.

There are many works we have not dealt with in this paper. For example, although we dealt with topology of some of the Messiaen's mode of limited transposition: the chromatic scale and the whole tone scales, other Messiaen’s modes may also lead to interesting musical icosahedra. In addition, we dealt with only one polyhedron, the icosahedron, and other polyhedra in any dimension may also connect with music beautifully. In this paper, we only dealt with branches of musical icosahedra directly connecting to the existing musical framework. However, in order to construct a mathematically beautiful theory, we should generalize our theory proposed in this paper and find the relationship between it and existing mathematical theories.

We hope our study stimulates researchers of music and mathematics, composers, musicians, and all other people who are interested in music and/or mathematics, and that it may inspire some breakthrough in mathematics, music, or other fields.

\newpage
\begin{appendix}

\section{List of scales}
\subsection{Major scales}
The $X$-major scale is given by {$X$, $X+2$, $X+4$, $X+5$, $X+7$, $X+9$, $X+11$} for any tone X. Explicitly, the major scales are represented as follows.
\\
\\
\indent
(i) $C$, $D$, $E$, $F$, $G$, $A$, $B$, $C$ ($C$-major),

(ii) $C\sharp$, $D\sharp$, $F$, $F\sharp$, $G\sharp$, $A\sharp$, $C$, $C\sharp$ ($C\sharp$-major),

(iii) $D$, $E$, $F\sharp$, $G$, $A$, $B$, $C\sharp$, $D$ ($D$-major),

(iv) $E\flat$, $F$, $G$, $A\flat$, $B\flat$, $C$, $D$, $E\flat$ ($E\flat$-major),

(v) $E$, $F\sharp$, $G\sharp$, $A$, $B$, $C\sharp$, $D\sharp$, $E$ ($E$-major),

(vi) $F$, $G$, $A$, $B\flat$, $C$, $D$, $E$, $F$ ($F$-major),

(vii) $F\sharp$, $G\sharp$, $A\sharp$, $B$, $C\sharp$, $D\sharp$, $F$, $F\sharp$ ($F\sharp$-major) ,

(viii) $G$, $A$, $B$, $C$, $D$, $E$, $F\sharp$, $G$ ($G$-major),

(ix) $A\flat$, $B\flat$, $C$, $D\flat$, $E\flat$, $F$, $G$, $A\flat$ ($A\flat$-major),

(x) $A$, $B$, $C\sharp$, $D$, $E$, $F\sharp$, $G\sharp$, $A$ ($A$-major),

(xi) $B\flat$, $C$, $D$, $E\flat$, $F$, $G$, $A$, $B\flat$ ($B\flat$-major),

(xii) $B$, $C\sharp$, $D\sharp$, $E$, $F\sharp$, $G\sharp$, $A\sharp$, $B$ ($B$-major).

\subsection{Minor scales}
The $X$-minor scale is given by {$X$, $X+2$, $X+3$, $X+5$, $X+7$, $X+8$, $X+10$}. Explicitly, the minor scales are represented as follows.
\\
\\
\indent
(i) $C$, $D$, $E\flat$, $F$, $G$, $A\flat$, $B\flat$, $C$ ($C$-minor),

(ii) $C\sharp$, $D\sharp$, $E$, $F\sharp$, $G\sharp$, $A$, $B$, $C\sharp$ ($C\sharp$-minor),

(iii) $D$, $E$, $F$, $G$, $A$, $B\flat$, $C$, $D$ ($D$-minor),

(iv) $E\flat$, $F$, $G\flat$, $A\flat$, $B\flat$, $B$, $D\flat$, $E\flat$ ($E\flat$-minor),

(v) $E$, $F\sharp$, $G$, $A$, $B$, $C$, $D$, $E$ ($E$-minor),

(vi) $F$, $G$, $A\flat$, $B\flat$, $C$, $D\flat$, $E\flat$, $F$ ($F$-minor),

(vii) $F\sharp$, $G\sharp$, $A$, $B$, $C\sharp$, $D$, $E$, $F\sharp$ ($F\sharp$-minor),

(viii) $G$, $A$, $B\flat$, $C$, $D$, $E\flat$, $F$, $G$ ($G$-minor),

(ix) $G\sharp$, $A\sharp$, $B$, $C\sharp$, $D\sharp$, $E$, $F\sharp$, $G\sharp$ ($G\sharp$-minor),

(x) $A$, $B$, $C$, $D$, $E$, $F$, $G$, $A$ ($A$-minor),

(xi) $B\flat$, $C$, $D\flat$, $E\flat$, $F$, $G\flat$, $A\flat$, $B\flat$ ($B\flat$-minor),

(xii) $B$, $C\sharp$, $D$, $E$, $F\sharp$, $G$, $A$, $B$ ($B$-minor).

\subsection{Hexatonic major scales}
The hexatonic major/minor scales and pentatonic major/minor scales are fundamental scales constructed by six (five) tones and obtained by removing one tone (two tones) from the major/minor scales (MacCabe, 2015). The hexatonic $X$-major scale is given by $X$, $X+2$, $X+4$, $X+5$, $X+7$, $X+9$. Explicitly, the hexatonic major scales are represented as follows.
\\
\\
\indent
(i) $C$, $D$, $E$, $F$, $G$, $A$, $C$,

(ii) $C\sharp$, $D\sharp$, $F$, $F\sharp$, $G\sharp$, $A\sharp$, $C\sharp$,

(iii) $D$, $E$, $F\sharp$, $G$, $A$, $B$, $D$,

(iv) $E\flat$, $F$, $G$, $A\flat$, $B\flat$, $C$, $E\flat$,

(v) $F$, $G$, $A$, $B\flat$, $C$, $D$, $F$,

(vi) $F\sharp$, $G\sharp$, $B\flat$, $B$, $C\sharp$, $E\flat$, $F\sharp$,

(viii) $G$, $A$, $B$, $C$, $D$, $E$, $G$,

(ix) $A\flat$, $B\flat$, $C$, $D\flat$, $E\flat$, $F$, $A\flat$,

(x) $A$, $B$, $C\sharp$, $D$, $E$, $F\sharp$, $A$,

(xi) $B\flat$, $C$, $D$, $E\flat$, $F$, $G$, $B\flat$,

(xii) $B$, $C\sharp$, $D\sharp$, $E$, $F\sharp$, $G\sharp$, $B$.

\subsection{Hexatonic minor scales}
The $X$-hexatonic minor scale is given by $X$, $X+2$, $X+3$, $X+5$, $X+7$, $X+10$. Explicitly, the hexatonic $X$-minor scales are represented as follows.
\\
\\
\indent
(i) $C$, $D$, $E\flat$, $F$, $G$, $B\flat$, $C$,

(ii) $C\sharp$, $D\sharp$, $E$, $F\sharp$, $G\sharp$, $B$, $C\sharp$,

(iii) $D$, $E$, $F$, $G$, $A$, $C$, $D$,

(iv) $E\flat$, $F$, $G\flat$, $A\flat$, $B\flat$, $D\flat$, $E\flat$,

(v) $E$, $F\sharp$, $G$, $A$, $B$, $D$, $E$,

(vi) $F$, $G$, $A\flat$, $B\flat$, $C$, $E\flat$, $F$,

(vii) $F\sharp$, $G\sharp$, $A$, $B$, $C\sharp$, $E$, $F\sharp$,

(viii) $G$, $A$, $B\flat$, $C$, $D$, $F$, $G$,

(ix) $G\sharp$, $A\sharp$, $B$, $C\sharp$, $D\sharp$, $F\sharp$, $G\sharp$,

(x) $A$, $B$, $C$, $D$, $E$, $G$, $A$,

(xi) $B\flat$, $C$, $D\flat$, $E\flat$, $F$, $A\flat$, $B\flat$,

(xii) $B$, $C\sharp$, $D$, $E$, $F\sharp$, $A$, $B$.

\subsection{Pentatonic major scales}
The pentatonic $X$-major scale is given by $X$, $X+2$, $X+4$, $X+7$, $X+9$. Explicitly, the pentatonic major scales are represented as follows.
\\
\\
\indent
(i) $C$, $D$, $E$, $G$, $A$, $C$,

(ii) $C\sharp$, $D\sharp$, $F$, $G\sharp$, $B\flat$, $C\sharp$,

(iii) $D$, $E$, $F\sharp$, $A$, $B$, $D$,

(iv) $E\flat$, $F$, $G$, $B\flat$, $C$, $E\flat$,

(v) $F$, $G$, $A$, $C$, $D$, $F$,

(vi) $F\sharp$, $G\sharp$, $A\sharp$, $C\sharp$, $D\sharp$, $F\sharp$,

(viii) $G$, $A$, $B$, $D$, $E$, $G$,

(ix) $A\flat$, $B\flat$, $C$, $E\flat$, $F$, $A\flat$,

(x) $A$, $B$, $C\sharp$, $E$, $F\sharp$, $A$,

(xi) $B\flat$, $C$, $D$, $F$, $G$, B$\flat$,

(xii) $B$, $C\sharp$, $D\sharp$, $F\sharp$, $G\sharp$, $B$.

\subsection{Pentatonic minor scales}
The pentatonic $X$-minor scale is given by $X$, $X+3$, $X+5$, $X+7$, $X+10$. Explicitly, the pentatonic minor scales are represented as follows.
\\
\\
\indent
(i) $C$, $E\flat$, $F$, $G$, $B\flat$, $C$,

(ii) $C\sharp$, $E$, $F\sharp$, $G\sharp$, $B$, $C\sharp$,

(iii) $D$, $F$, $G$, $A$, $C$, $D$,

(iv) $E\flat$, $G\flat$, $A\flat$, $B\flat$, $D\flat$, $E\flat$,

(v) $E$, $G$, $A$, $B$, $D$, $E$,

(vi) $F$, $A\flat$, $B\flat$, $C$, $E\flat$, $F$,

(vii) $F\sharp$, $A$, $B$, $C\sharp$, $E$, $F\sharp$,

(viii) $G$, $B\flat$, $C$, $D$, $F$, $G$,

(ix) $G\sharp$, $B$, $C\sharp$, $D\sharp$, $F\sharp$, $G\sharp$,

(x) $A$, $C$, $D$, $E$, $G$, $A$,

(xi) $B\flat$, $D\flat$, $E\flat$, $F$, $A\flat$, $B\flat$,

(xii) $B$, $D$, $E$, $F\sharp$, $A$, $B$. 
\\
\\
\indent
We remark that the major/minor scales each include one pair of tritones, e.g. the $C$-major scale includes $F$ and $B$ and the $C$-minor scale includes $D$ and $A\flat$. Then, the hexatonic major/minor scale is obtained by removing one of them, e.g. the hexatonic $C$-major scale can be obtained by removing $B$ from the $C$-major scale and the hexatonic $C$-minor scale can be obtained by removing $A\flat$ from the $C$-minor scale. In addition, the pentatonic major/minor scale is obtained by removing the two tones constructing a tritone, e.g. the pentatonic $C$-major scale can be obtained by removing $F$ and $B$ from the $C$-major scale and the pentatonic $C$-minor scale can be obtained by removing $D$ and $A\flat$ from the $C$-minor scale.

\subsection{Gregorian modes}
The Gregorian modes are represented as follows.
\\
\\
\indent
(i) $D(f)$, $E$, $F$, $G$, $A$, $B$, $C$, $D$ (Dorian),

(ii) $A$, $B$, $C$, $D(f)$, $E$, $F$, $G$, $A$ (Hypodorian),

(iii) $E(f)$, $F$, $G$, $A$, $B$, $C$, $D$, $E$ (Phrygian),

(iv) $B$, $C$, $D$, $E(f)$, $F$, $G$, $A$, $B$ (Hypophrigian),

(v) $F(f)$, $G$, $A$, $B$, $C$, $D$, $E$, $F$ (Lydian),

(vi) $C$, $D$, $E$, $F(f)$, $G$, $A$, $B$, $C$ (Hypolydian),

(vii) $G(f)$, $A$, $B$, $C$, $D$, $E$, $F$, $G$ (Myxolydian),

(viii) $D$, $E$, $F$, $G(f)$, $A$, $B$, $C$, $D$ (Hipomyxolydian),
\\
\\
\indent
where $f$ mean the final tone of each mode. In this paper, we deal with the Gregorian modes and all the transpositions of them. We define $X$-Dorian ($X$-Phrygian, $X$-Lydian, $X$-Mixolydian) as the scale whose base tone is $X$ and is obtained by transposing Dorain (Phrygian, Lydian, Mixolydian). They are represented as follows. 
\\
\\
\indent
$X$, $X+2$, $X+3$, $X+5$, $X+7$, $X+9$, $X+10$, $X+12$ ($X$-Dorian),
\\
\indent
$X$, $X+1$, $X+3$, $X+5$, $X+7$, $X+8$, $X+10$, $X+12$ ($X$-Phrygian),
\\
\indent
$X$, $X+2$, $X+4$, $X+6$, $X+7$, $X+9$, $X+11$, $X+12$ ($X$-Lydian),
\\
\indent
$X$, $X+2$, $X+4$, $X+5$, $X+7$, $X+9$, $X+11$, $X+12$ ($X$-Mixolydian).
\\
\\
Then, we analyze the structure of the above modes. The $X$-Dorain can be obtained by lowering the sixth tone of the $X$-minor scale by one semitone, and the $X$-Phrygian can be obtained by raising the second tone of the $X$-minor scale by one semitone, and the $X$-Lydian can be obtained by raising the fourth tone of the $X$-major by one semitone, and the $X$-Myxolydian can be obtained by lowering the seventh tone of the $X$-major by one semitone.

\subsection{Chromatic Scale and Whole tone scales}
The chromatic scale is the following scale.
\\
\\
\indent
$C$, $C\sharp$, $D$, $E\flat$, $E$, $F$, $F\sharp$, $G$, $G\sharp$, $A$, $B\flat$, $B$
\\
\\
\indent
Note that all the intervals in the chromatic scale are one semitone.
\\
\indent
The whole tone scales are the following two scales.
\\
\\
\indent
$C$, $D$, $E$, $F\sharp$, $G\sharp$, $B\flat$,
\\
\indent
$C\sharp$, $E\flat$, $F$, $G$, $A$, $B$.
\\
\\
\indent
Note that all the intervals in the whole tone scale are one whole tone (two semitones). The second one is obtained by transposing the first one by one whole tone (two semitones). Messiaen extended these scales by using the concept of transposition. There are modes that can be transposed in 1 way, 2 ways, 3 ways, 4 ways, 6 ways and 12 ways. The chromatic scale is the only scale that can be transposed in 1 way, and the whole tone scales are the only scales that can be transposed in 2 ways. The major/minor scales or the Gregorian modes can be transposed in 12 ways. Although the modes of limited transposition we deal with mainly in this paper are only the chromatic scale and whole tone scales, other modes of limited transposition, e.g. $C$, $C\sharp$, $F$, $F\sharp$, $G$, $G\sharp$, $B\flat$, are also interesting and we will deal with those modes in another paper.

\subsection{First generalization of major scales} 
In this subsection, we show the first generalization of the major scales and the three generalized major/minor triads constructing them. Note that the following scales are rearranged to become an ascending scale because these arrangements may be useful for composition or analysis of music because a scale is usually considered as a list of tones that is ascending or descending. Theoretically, the original arrangements should be respected because, except for the $C$-major scale, the following scales correspond to figures that are different from the figure of the $C$-major scale on the type 1 chromatic/whole tone musical icosahedron.
\\
\\
\indent
$C$-major: $C$, $D$, $E$, $F$, $G$, $A$, $B$ ($CEG$, $FAC$, $GBD$),
\\
\indent
$C_1$-major: $C$, $E\flat$, $E$, $F$, $G$, $G\sharp$, $B\flat$ ($CGG\sharp$, $EFC$, $G\sharp B\flat E\flat$),
\\
\indent
$C_2$-major: $C$, $C\sharp$, $E$, $G$, $G\sharp$, $A$, $B$ ($CG\sharp A$, $GEC$, $AC\sharp B$),
\\
\indent
$C_3$-major: $C$, $D$, $F$, $G$, $G\sharp$, $A$, $B\flat$ ($CAF$, $G\sharp GC$, $FDB\flat$),
\\
\indent
$C_4$-major: $C$, $C\sharp$, $E\flat$, $E$, $F$, $G\sharp$, $A$ ($CFE$, $G\sharp AC$, $EE\flat C\sharp$),
\\
\indent
$MC$-major: $C$, $E\flat$, $F$, $G$, $G\sharp$, $A$, $B\flat$ ($CG\sharp G$, $AFC$, $GE\flat B\flat$),
\\
\indent
$MC_1$-major: $C$, $C\sharp$, $E$, $F$, $G\sharp$, $A$, $B$ ($CAG\sharp$, $EFC$, $G\sharp B C\sharp$),
\\
\indent
$MC_2$-major: $C$, $D$, $E$, $F$, $G$, $A$, $B\flat$ ($CFA$, $EGC$, $AB\flat D$),
\\
\indent
$MC_3$-major: $C$, $C\sharp$, $E\flat$, $E$, $F$, $G$, $G\sharp$ ($CEF$, $GG\sharp C$, $FC\sharp E\flat$),
\\
\indent
$MC_4$-major: $C$, $D$, $E$, $G$, $G\sharp$, $A$, $B$ ($CGE$, $G\sharp AC$, $EDB$),
\\
\indent
and the scales obtained by raising all the tones of each of the above scales by $n$ semitones ($1\leq n\leq12$).

\subsection{First generalization of minor scales}
In this subsection, we show the first generalization of the minor scales and the three generalized major/minor triads constructing them. Note that the following scales are rearranged to become an ascending scale.  
\\
\\
\indent
$C$-minor: $C$, $D$, $E\flat$, $F$, $G$, $G\sharp$, $B\flat$ ($CE\flat G$, $FG\sharp C$, $GB\flat D$),
\\
\indent
$C_1$-minor: $C$, $E$, $F$, $G$, $G\sharp$, $A$, $B$ ($CGG\sharp$, $EFC$, $G\sharp AB$),
\\
\indent
$C_2$-minor: $C$, $C\sharp$, $E\flat$, $E$, $F$, $G\sharp$, $B\flat$ ($CG\sharp F$, $E\flat EC$, $FC\sharp B\flat$),
\\
\indent
$C_3$-minor: $C$, $D$, $E\flat$, $E$, $F$, $G$, $A$ ($CFE$, $GE\flat C$, $EDA$),
\\
\indent
$C_4$-minor: $C$, $C\sharp$, $E\flat$, $E$, $G$, $G\sharp$, $B$ ($CEE\flat$, $G\sharp GC$, $E\flat BC\sharp$),
\\
\indent
$MC$-minor: $C$, $D$, $E\flat$, $E$, $F$, $G$, $B\flat$ ($CGE\flat$, $FEC$, $E\flat DB\flat$), 
\\
\indent
$MC_1$-minor: $C$, $E\flat$, $E$, $G$, $G\sharp$, $A$, $B$ ($CG\sharp G$, $EE\flat C$, $GBA$), 
\\
\indent
$MC_2$-minor: $C$, $C\sharp$, $E\flat$, $F$, $G$, $G\sharp$, $B\flat$ ($CFG\sharp$, $E\flat GC$, $G\sharp B\flat C\sharp$),
\\
\indent
$MC_3$-minor: $C$, $D$, $E$, $F$, $G$, $G\sharp$, $A$ ($CEF$, $GG\sharp C$, $FAD$),
\\
\indent
$MC_4$-minor: $C$, $C\sharp$, $E\flat$, $E$, $F$, $G\sharp$, $B$ ($CE\flat E$, $G\sharp FC$, $EC\sharp B$),
\\
\indent
and the scales obtained by raising all the tones of each of the above scales by $n$ semitones ($1\leq n\leq12$)

\subsection{Second generalization of major scales}
In this subsection, we show the second generalization of the major scales and the three generalized major/minor triads constructing them. Note that the following scales are rearranged to become an ascending scale.  
\\
\\
\indent
$C$-major: $C$, $D$, $E$, $F$, $G$, $A$, $B$ ($CEG$, $FAC$, $GBD$),
\\
\indent
$C_1'$-major: $C$, $C\sharp$, $E\flat$, $F$, $G\sharp$, $A$, $B\flat$ ($CC\sharp F$, $E\flat G\sharp C$, $FAB\flat$),
\\
\indent
$C_2'$-major: $C$, $D$, $E\flat$, $E$, $G$, $G\sharp$, $B$ ($CBE\flat$, $DEC$, $E\flat G\sharp G$),
\\
\indent
$C_3'$-major: $C$, $C\sharp$, $D$, $E$, $F$, $A$, $B\flat$ ($CAD$, $B\flat C\sharp C$, $DEF$),
\\
\indent
$C_4'$-major: $C$, $C\sharp$, $E\flat$, $G$, $G\sharp$, $B\flat$, $B$ ($CG\sharp B\flat$, $GBC$, $B\flat C\sharp E\flat$),
\\
\indent
$MC'$-major: $C$, $C\sharp$, $E\flat$, $F$, $G\sharp$, $B\flat$, $B$ ($CG\sharp E\flat$, $FC\sharp C$, $E\flat BB\flat$),
\\
\indent
$MC_1'$-major: $C$, $D$, $E\flat$, $E$, $G$, $A$, $B$ ($CED$, $E\flat BC$, $DAG$),
\\
\indent
$MC_2'$-major: $C$, $C\sharp$, $D$, $F$, $G\sharp$, $A$, $B\flat$ ($CG\sharp B\flat$, $DAC$, $B\flat G\sharp F$), 
\\
\indent
$MC_3'$-major: $C$, $E\flat$, $E$, $G$, $G\sharp$, $B\flat$, $B$ ($CBG$, $B\flat G\sharp C$, $GEE\flat$),
\\
\indent
$MC_4'$-major: $C$, $C\sharp$, $D$, $E$, $F$, $G$, $A$ ($CAF$, $GEC$, $FC\sharp D$),
\\
\indent
and the scales obtained by raising all the tones of each of the above scales by $n$ semitones ($1\leq n\leq12$)

\subsection{Second generalization of minor scales}
In this subsection, we show the second generalization of the minor scales and the three generalized major/minor triads constructing them. Note that the following scales are rearranged to become an ascending scale.  
\\
\\
\indent
$C$-minor: $C$, $D$, $E\flat$, $F$, $G$, $A\flat$, $B\flat$ ($CE\flat G$, $FA\flat C$, $GB\flat D$),
\\
\indent
$C_1'$-minor: $C$, $D$, $E$, $F$, $G$, $A$, $B$ ($CEA$, $GBC$, $ADF$),
\\
\indent
$C_2'$-minor: $C$, $C\sharp$, $F$, $G$, $G\sharp$, $A$, $B\flat$ ($CG\sharp B\flat$, $AC\sharp C$, $B\flat F G$),
\\
\indent
$C_3'$-minor: $C$, $D$, $E\flat$, $G$, $A$, $B\flat$, $B$ ($CBD$, $B\flat E\flat C$, $DGA$),
\\
\indent
$C_4'$-minor: $C$, $C\sharp$, $D$, $E$, $F$, $A$, $B\flat$ ($CC\sharp F$, $DEC$, $FAB\flat$),
\\
\indent
$MC'$-minor: $C$, $D$, $E$, $G$, $A$, $B\flat$, $B$ ($CBG$, $AEC$, $GDB\flat$),
\\
\indent
$MC_1'$-minor: $C$, $C\sharp$, $D$, $F$, $G\sharp$, $A$, $B\flat$ ($CC\sharp A$, $B\flat G\sharp C$, $AFD$),
\\
\indent
$MC_2'$-minor: $C$, $D$, $E\flat$, $F$, $G$, $B\flat$, $B$ ($CE\flat B\flat$, $DBC$, $B\flat GF$),
\\
\indent
$MC_3'$-minor: $C$, $C\sharp$, $D$, $E$, $F$, $G$, $A$ ($CED$, $FC\sharp C$, $DAG$),
\\
\indent
$MC_4'$-minor: $C$, $E\flat$, $F$, $G$, $G\sharp$, $A$, $B\flat$ ($CG\sharp F$, $GE\flat C$, $FB\flat A$),
\\
\indent
and the scales obtained by raising all the tones of each of the above scales by $n$ semitones ($1\leq n\leq12$)

\end{appendix}

Archer-Hind, R. D. (1888). \emph{The Timaeus of Plato}, London: Macmillan.

Benary, P. (1995). Dur und Moll. In \emph{Musik in Geschichte und Gegenwart} (2nd ed., col. 1593-1597),
\\
\qquad\quad Kassel: B\"{a}renreiter-Verlag.

Coxeter, H. S. M. \& Rouce Ball, W. W. (2010). \emph{Mathematical Recreations and Essays} (13th. ed.). 
\\
\qquad\quad
New York, NY: Dover Publications.

Cromwell, P. R. (1999). \emph{Polyhedra}, Cambridge, UK: Cambridge University Press.

Hartung, P. Ch. (1749). \emph{Musicus theoretico-practicus}, N\"{u}rnberg.

Hawking, S. (2007). \emph{God Created The Integers: The Mathematical Breakthroughs that Changed History},
\\
\qquad\quad
Philadelphia: Running Press Adult.

JAXA EORC (2019, October 15). AMSR AMSR-E,
\\
\qquad\quad
Retrieved from https://sharaku.eorc.jaxa.jp/AMSR/giftbox/s\_ball/index\_j.html.

Kepler, J. (2009). \emph{Mysterium Cosmographicum} [Cosmic Mystery] (S. Ohtsuki \& M. Kishimoto, Trans.). 
\\
\qquad\quad
Tokyo: Kohsakusha.

Lindley, M. (2001). "Pythagorean intonation." 
\\
\qquad\quad
${\it The  \ New \ Grove \  Dictionary \ of \ Music \ and \ Musicians}$. 2nd ed. London: Macmillan.

MacCabe, L. (2015). \emph{You Can Teach Yourself Songwriting} (You Can Teach Yourself Series), Fenton, 
\\
\qquad\quad
Mel Bay Publications.

Meisner, G. B. \& Araujo, R. (2018). \emph{The Golden Ratio: The Divine Beauty of Mathematics},
\\
\qquad\quad
New York, NY: Race Point Publishing.

Mitzler, Lorenz Christoph. Trans. (1742). \emph{radus ad Parnassum, oder Anf\"{u}hrung zur regelm\"{a}ssigen 
\\
\qquad\quad
Composition, aus dem Lateinischen ins Deutsche \"{u}bersetzt, und mit Anmerkungen versehen}, Leipzig.

Polyhedron Club (2019a, November 23). Polyhedron-Temari collection,

Retrieved from http://polyhedra.cocolog-nifty.com/blog/2012/11/post-a6c3.html.

Polyhedron Club (2019b, November 23). Polyhedron-paper-craft, 
\\
\qquad\quad
Retrieved from http://polyhedra.cocolog-nifty.com/blog/2017/11/post-dcc9.html.

Schein, S. \& Gayed, M. (2014). Fourth class of convex equilateral polyhedron with polyhedral symmetry 
\\
\qquad\quad
related to fullerenes and viruses, \emph{Proceedings of the National Academy of Sciences}, \it{111}\rm{(8)}, 2920-2925.

Shi C. Y. (2002). On the gong diao of Chinese Traditional Music Theory, \emph{Musicology in China}, \it{4}, 109-116.

\end{document}